\documentclass[journal,twoside,web]{ieeecolor}
\usepackage{generic}
\usepackage{cite}
\usepackage{amsmath,amssymb,amsfonts}
\usepackage{graphicx}
\usepackage{algorithm,algorithmic}
\usepackage{hyperref}
\hypersetup{hidelinks=true}
\usepackage{textcomp}
\usepackage{tikz}
\usepackage{pgfplots}
\pgfplotsset{compat=1.18}
\usetikzlibrary{positioning,arrows.meta,calc,shapes.geometric}
\definecolor{ieeegrid}{rgb}{0.78,0.78,0.78}
\definecolor{ieeelegend}{rgb}{0.70,0.70,0.70}
\pgfplotsset{
  ieeeplot/.style={
    grid=both,
    grid style={dashed, ieeegrid, very thin},
    tick label style={font=\footnotesize},
    label style={font=\small},
    legend style={font=\footnotesize, draw=ieeelegend, inner sep=3pt},
    every axis plot/.append style={thick},
    mark size=1.2pt,
    line width=1.2pt,
  }
}
\def\BibTeX{{\rm B\kern-.05em{\sc i\kern-.025em b}\kern-.08em
    T\kern-.1667em\lower.7ex\hbox{E}\kern-.125emX}}
\newtheorem{theorem}{Theorem}
\newtheorem{lemma}{Lemma}
\newtheorem{proposition}{Proposition}
\newtheorem{corollary}{Corollary}
\newtheorem{definition}{Definition}
\newtheorem{assumption}{Assumption}
\newtheorem{condition}{Condition}
\newtheorem{remark}{Remark}

\begin{document}

\title{Projection-Free Gradient Tracking Neurodynamic for Distributed Nonnegative Sparse Optimization} 

\author{Juntao Li, Qihang Guo, and Deyuan Meng
\thanks{J. Li and Q. H. Guo are with the School of Mathematics and Statistics, Henan Normal University, Xinxiang 453007, P. R. China (e-mail: juntaolimail@126.com).}
\thanks{D. Meng is with the Seventh Research Division, School of Automation Science and Electrical Engineering, Beihang University (BUAA), Beijing 100191, P. R. China (e-mail: dymeng@buaa.edu.cn).}
\thanks{This work was supported in part by the National Natural Science Foundation of China under Grants U2333215 and 62273018.}
\thanks{Corresponding author: Deyuan Meng.}
}

\maketitle

\begin{abstract}
This paper proposes a continuous-time projection-free
gradient-tracking neurodynamic system for distributed
nonnegative sparse optimization, where projection-free means
that no state is projected onto the feasible set. Gradient
tracking supplies an exact estimate of the average gradient
under heterogeneous local costs, and a perturbed dual dynamics
couples the dual variable to the tracker. The primal and the
dual are free, may take negative values during transients, and
feasibility, nonnegativity and sparsity are produced by the
equilibrium equations themselves; the only nonsmooth nonlinearity
is the scalar positive part of $\boldsymbol{\lambda}-\mathbf{x}$,
which is Lipschitz. At an
equilibrium the perturbed dual reduces to
$\min\{\boldsymbol{\lambda}^*,\mathbf{x}^*\}=
\beta(\mathbf{q}^*+\kappa\mathbf{1}-\boldsymbol{\lambda}^*)$,
and an extremal argument on $\mathbf{L}\mathbf{x}^*$ forces
consensus and complementary slackness, so every equilibrium is
an exact KKT point with multiplier
$\boldsymbol{\lambda}^*=\nabla F(\mathbf{x}^*)+\kappa\mathbf{1}
\succeq\mathbf{0}$, whose positive entries mark the strictly
active coordinates.
The vector field is
locally Lipschitz, hence classical solutions are unique for
every initial condition.
For convex objectives we certify three matrix inequalities.
The first yields a centralized kernel rate $\gamma_c$, proved
global for quadratic costs by bounding the $\min$ term with an
active-set selector at every state. The second is an LMI posed
on the Jacobian at the KKT point, restricted to the invariant
subspace of the conserved mean tracking error, which yields a
computable distributed rate $\gamma_J$ with no gain product and
no stability assumption on an isolated subsystem. The third
writes the dual error in exact Lur'e form, with a constant
linear part and a shifted positive part whose sector holds at
every state; for quadratic costs it gives global exponential
stability of the distributed flow from every initial condition
with $\mathbf{C}\mathbf{z}(0)=\mathbf{0}$ (in particular every
gradient-tracking initialization), with an explicit bias
proportional to the mean tracker initialization error when that
condition is dropped, and an explicit quadratic form at the full state
dimension certifies a rate within five percent of the
Jacobian abscissa on the nominal instance.
The three layers separate quantitatively: $\beta$ and $\rho$
raise the rate, $\gamma_J$ is the abscissa of the Jacobian on
that subspace, the slow abscissa is numerically insensitive to
the momentum time scale $\varepsilon$ over $[0.01,0.2]$, and
on sparse graphs the rate is capped by $\nu_2$.
Numerical experiments on strongly convex quadratic problems
reproduce the nominal Jacobian abscissa to within one percent on
the tail (eight percent on the pre-registered finite window),
identify the slow modes as dual-dominated, and display the
recovery of the active set. Strong local heterogeneity slows the
finite-horizon transient by about fifteen percent, while the tail
rate changes by less than two percent.
\end{abstract}


\section{Introduction}
\label{sec:introduction}

\IEEEPARstart{D}{istributed} nonnegative sparse
optimization arises in sparse nonnegative matrix
factorization, portfolio selection without short
selling, sparse signal recovery with physical
nonnegativity, and multi-site feature selection,
where agents hold private data and must agree on a
common sparse nonnegative solution. Nonnegativity
enforces feasibility and sparsity encodes the prior
that few coordinates are active. Privacy concerns in
collaborative studies \cite{zhang2023ppml} favor
distributed solutions in which each agent retains its
local data, yet distributed nonnegative sparse
optimization remains largely unexplored.

Distributed gradient tracking
\cite{shi2015extra,xi2016distr}---building on
early consensus and subgradient methods
\cite{tsitsiklis1986distributed,nedic2009distributed}---estimates the
global gradient average via local communication,
with recent optimal variants achieving tight
complexity bounds \cite{song2024optimal} and
stochastic extensions under general variance
conditions \cite{huang2024stochastic,
pu2018push}.
Continuous-time variants
\cite{zeng2017distributed,wu2023distributed,
wang2025finite} address coupling constraints
or pseudoconvex objectives, and nonconvex
extensions rely on the Polyak--\L{}ojasiewicz
condition \cite{xu2024quantized,pu2026stochastic}.
These methods target unconstrained
or equality-constrained problems, so they
do not handle inequality constraints via
projection; their proofs also require gradient
monotonicity, which fails for nonconvex costs
with changing active sets.

Neurodynamic optimization embeds the KKT
conditions in an ODE whose equilibrium coincides
with the optimal solution.
The projection neural network (PNN)
\cite{rodriguez1990recurrent,xia2004general}, rooted in
neural circuits for optimization
\cite{hopfield1984neural,kennedy1988neural}, has
been extended to trust-region programming
\cite{huang2025trpnn}, inertial primal--dual
optimization \cite{zhao2025inertial}, and
nonconvex games \cite{xia2024nash};
see \cite{xia2024survey,zhao2022l1,xiao2024fpga}
for surveys and hardware.
Distributed PNN methods
\cite{xia2025bilevel,zhao2023smoothing,
zhao2023projection,xia2024event}
rely on consensus projection rather than gradient
tracking, so under heterogeneous local costs the
equilibrium satisfies only a perturbed KKT system.
Alternatively, augmented Lagrangian
relaxations
\cite{xia2025distributed,guan2025neurodynamic}
smooth the constraint at the price of local
stability or meta-heuristic search.
Embedding gradient tracking into the PNN resolves
the perturbed-KKT issue, but every existing construction
still relies on a projection operator to keep the primal,
or the primal and the dual, inside the orthant; the projection
freezes the boundary coordinates, so any Lyapunov certificate
must handle active-set switches. The proposed flow removes
the state projection and enforces the orthant only at the
equilibrium.

Nonconvex costs lie outside the scope of this paper: both the
active-set argument and the global sector use convexity, and the
KL-based guarantees \cite{bolte2014proximal,jin2017escape,
swenson2020distributed} would have to hold across
projection-induced active-set switches.

A continuous-time framework that simultaneously
incorporates gradient tracking and constraint handling
for inequality constraints remains open.
The main contributions are threefold.

\begin{enumerate}
\item A continuous-time projection-free gradient-tracking
neurodynamic system for distributed nonnegative
sparse optimization. Unlike recent distributed
neurodynamic approaches that rely on consensus
projection \cite{xia2024event,zhao2023projection},
inertial proximal operators
\cite{zhao2024directed}, or augmented Lagrangian
relaxations \cite{xia2025distributed,guan2025neurodynamic}---none
of which incorporates gradient tracking---the proposed
system embeds gradient tracking into the PNN structure
so that the equilibrium satisfies the exact KKT
system. A key innovation is a perturbed dual
dynamics that couples the dual variable to the
gradient tracker, which supplies the offset
$\kappa\mathbf{1}_{nd}$ needed for exact
complementarity and thus for an unperturbed KKT
point. The primal and the dual are not projected:
they are free, may leave the orthant during transients,
and consensus, nonnegativity and complementary
slackness are produced by the equilibrium equations
(Proposition~\ref{prop:equilibrium}). The resulting
vector field is locally Lipschitz, so the flow is a
classical ODE with unique solutions for every initial
condition---no Filippov convention and no forward
invariance of the orthant are required.

\item A two-layer matrix-inequality analysis for the kernel
and the distributed Jacobian. The first
layer certifies the centralized primal--dual
kernel through an explicit matrix inequality whose margin
is positive for every active-set selector; for quadratic
costs this yields a proof of global exponential stability
of the kernel, and the rate $\gamma_c$ is computed rather
than assumed. The second layer certifies the distributed
flow through an LMI posed on the Jacobian at the KKT point
and restricted to the invariant subspace of the mean
tracking error. Both inequalities are evaluated offline
from the problem data, with no trajectory involved, and
neither requires a gain product or a stability assumption
on an isolated subsystem. The certified rates show that
$\beta$ and $\rho$ raise the rate, that the distributed
rate coincides with the Jacobian abscissa on the invariant
manifold, that the slow abscissa is numerically insensitive to
the momentum time scale over $[0.01,0.2]$, and that on sparse graphs
the rate becomes capped by $\nu_2$, at a decreasing fraction of it
($0.79\,\nu_2$ at $n=20$ and $0.40\,\nu_2$ at $n=40$).

\item A global Lur'e certificate for the distributed flow.
In error coordinates the only nonlinear term is the dual
positive part, and a shift by the KKT pair turns it into a
monotone positive part with slope in $[0,1]$. The
linear part of the error dynamics is then constant---it
does not depend on the active set---and the shifted
positive part satisfies an exact sector inequality at
every state. A single quadratic form, certified by one LMI
with a diagonal multiplier on the invariant subspace of
the conserved mean tracking error, gives global
exponential convergence for quadratic costs from every
admissible initial condition, with no orthant invariance
and no neighbourhood restriction.
On the nominal instance, at the full state dimension, an
explicit quadratic form certifies $\gamma_g=0.27$---within
five percent of the Jacobian abscissa---so the global form
loses little relative to the sharp local linearization.
\end{enumerate}

The remainder is organized as follows.
Section~\ref{sec:problem_formulation} states the problem,
Section~\ref{sec:algorithm_design_properties} the algorithm,
Section~\ref{sec:convergence} the convex analysis,
Section~\ref{sec:num_experiments} the experiments, and
Section~\ref{sec:conclusion} concludes.
All results are proved in the body.

\section{Problem Formulation}
\label{sec:problem_formulation}

We formalize the distributed nonnegative
sparse optimization problem.
The centralized formulation is
\begin{equation}
  \min_{\mathbf{x} \ge 0}
  \; \frac{1}{n}\sum_{i=1}^n f_i(\mathbf{x})
  + \kappa\|\mathbf{x}\|_1,
  \label{eq:centralized_prob}
\end{equation}
where $\kappa \ge 0$ controls sparsity and
each $f_i(\mathbf{x})$ is a smooth local cost
function.

The additive structure of the objective lends
itself naturally to a distributed implementation.
Instead of aggregating all local losses at a
central node, one can assign each $f_i$ to an
agent $i$ (corresponding to a sample or an
institutional site), let each
agent maintain a local estimate
$\mathbf{x}_i$, and enforce agreement across
agents via a consensus constraint.
Communication is restricted to neighboring
agents over an undirected connected graph
$\mathcal{G} = (\mathcal{V}, \mathcal{E})$,
so raw data need not leave its local site.
This leads to the distributed formulation
\begin{equation}
  \min_{\mathbf{x}_1,\ldots,\mathbf{x}_n \ge 0}
  \;\frac{1}{n}\sum_{i=1}^n f_i(\mathbf{x}_i)
  + \frac{\kappa}{n}\sum_{i=1}^n
  \|\mathbf{x}_i\|_1
  \quad\text{s.t.}\quad
  \mathbf{x}_1 = \cdots = \mathbf{x}_n,
  \label{eq:distributed_prob}
\end{equation}
where the consensus constraint ensures all
agents converge to the same solution. Under
consensus, \eqref{eq:distributed_prob} reduces to
the centralized form~\eqref{eq:centralized_prob},
confirming consistency.
Throughout, vectors and matrices are denoted by
boldface lowercase and uppercase letters,
respectively, the Euclidean and spectral norms
are both $\|\cdot\|$, the Kronecker product is
$\otimes$, $[\mathbf{x}]_+$ denotes the
coordinatewise positive part,
$[\mathbf{x}]_{+,j}=\max\{x_j,0\}$, the graph Laplacian is
$\mathcal{L}\in\mathbb{R}^{n\times n}$ with
Kronecker extension
$\mathbf{L} = \mathcal{L}\otimes I_d$, and
$\nu_2(\cdot)$ denotes the second smallest
eigenvalue.

\section{Algorithm Design and Equilibrium Analysis}
\label{sec:algorithm_design_properties}

\subsection{Algorithm Design}
\label{subsec:algorithm_design}

We now develop a continuous-time distributed
algorithm for solving
problem~\eqref{eq:distributed_prob}.
On the nonnegative orthant the $\ell_1$-norm is
linear, $\|\bar{\mathbf{x}}\|_1 =
\mathbf{1}_d^\top\bar{\mathbf{x}}$, so
\eqref{eq:centralized_prob} becomes
$\min_{\bar{\mathbf{x}}\ge0}g(\bar{\mathbf{x}})$ with
$g\triangleq F+\kappa\mathbf{1}_d^\top\bar{\mathbf{x}}$
smooth on the feasible set and
$\nabla g=\nabla F+\kappa\mathbf{1}_d$;
proximal splitting and subgradient methods are
therefore unnecessary, and one Lyapunov analysis
covers the constrained dynamics. The distributed
version inherits this property with
$g_i\triangleq f_i+\kappa\mathbf{1}_d^\top\mathbf{x}_i$.

Building on this reformulation, we propose a
continuous-time dynamical system
that solves the distributed optimization
problem~\eqref{eq:distributed_prob} as a
\emph{feedback control system} realized by a
\emph{projection-free neurodynamic network}.
Each node $i$ is an autonomous controller
with state $\mathbf{x}_i, \mathbf{m}_i,
\mathbf{q}_i, \boldsymbol{\sigma}_i,
\boldsymbol{\lambda}_i \in \mathbb{R}^d$,
driven by local measurements
($\nabla^2 f_i(\mathbf{x}_i)$) and neighbor
communications, with no central coordinator.

Each node executes:
\begin{align*}
  \begin{split}
    \varepsilon \dot{\mathbf{m}}_i &= -\eta \mathbf{m}_i - \mathbf{q}_i + \boldsymbol{\lambda}_i \\
    &\quad - \rho \textstyle\sum_{j\in\mathcal{N}_i} (\mathbf{x}_i - \mathbf{x}_j) - \kappa \mathbf{1}_d,
  \end{split} \\
  \dot{\mathbf{x}}_i &= \mathbf{m}_i, \\
  \begin{split}
    \dot{\mathbf{q}}_i &= -\textstyle\sum_{j\in\mathcal{N}_i} \bigl[(\mathbf{q}_i + \boldsymbol{\sigma}_i) \\
    &\quad - (\mathbf{q}_j + \boldsymbol{\sigma}_j)\bigr] + \nabla^2 f_i(\mathbf{x}_i)\,\dot{\mathbf{x}}_i,
  \end{split} \\
  \dot{\boldsymbol{\sigma}}_i &= -\textstyle\sum_{j\in\mathcal{N}_i} (\mathbf{x}_i - \mathbf{x}_j)
  - \gamma_s \boldsymbol{\sigma}_i, \\
  \dot{\boldsymbol{\lambda}}_i &= -(1{+}\beta)\boldsymbol{\lambda}_i
  + [\boldsymbol{\lambda}_i - \mathbf{x}_i]_+
  + \beta(\mathbf{q}_i + \kappa\mathbf{1}_d),
\end{align*}
where $[\cdot]_+$ is the coordinate-wise positive
part, $[\mathbf{v}]_{+,j}=\max\{v_j,0\}$.
Neither $\mathbf{x}_i$ nor $\boldsymbol{\lambda}_i$ is
projected: both are free variables.
The vector field is locally Lipschitz, since
$[\cdot]_+$ is globally Lipschitz, $\mathbf{L}$ is constant,
and $\mathbf{H}(\mathbf{x})\dot{\mathbf{x}}$ is locally
Lipschitz under
Assumptions~\ref{assump:A1} and~\ref{assump:A4}; hence
every initial condition generates a unique classical
solution on $[0,\infty)$
(Theorem~\ref{thm:well_posedness}).
The orthant is neither invariant nor required to be:
$\mathbf{x}_i$ and $\boldsymbol{\lambda}_i$ may be negative
during the transient. Feasibility, nonnegativity and
complementarity are instead properties of the
\emph{equilibria}, derived from the equilibrium equations
alone in Section~\ref{subsec:equilibrium}. This removes the
orthant projection on the primal from the algorithm. The vector
field remains Lipschitz rather than globally smooth, because the
dual equation still carries the positive part
$[\boldsymbol{\lambda}-\mathbf{x}]_+$; but Lipschitz continuity is
enough for classical solutions to exist everywhere, and no
Filippov regularization or forward-invariance argument of a
projected formulation is needed.

Stacking variables into global vectors and
using the Kronecker Laplacian
$\mathbf{L} = \mathcal{L}\otimes I_d$,
the system becomes:
\begin{align}
  \varepsilon \dot{\mathbf{m}} &= -\eta \mathbf{m} - \mathbf{q} + \boldsymbol{\lambda}
  - \rho \mathbf{L}\mathbf{x} - \kappa \mathbf{1}_{nd}, \label{eq:ode_m} \\
  \dot{\mathbf{x}} &= \mathbf{m}, \label{eq:ode_x} \\
  \dot{\mathbf{q}} &= -\mathbf{L}(\mathbf{q} + \boldsymbol{\sigma}) + \mathbf{H}(\mathbf{x})\,\dot{\mathbf{x}}, \label{eq:ode_q} \\
  \dot{\boldsymbol{\sigma}} &= -\mathbf{L}\mathbf{x} - \gamma_s\boldsymbol{\sigma}, \label{eq:ode_sigma} \\
  \dot{\boldsymbol{\lambda}} &= -(1{+}\beta)\boldsymbol{\lambda}
  + [\boldsymbol{\lambda} - \mathbf{x}]_+
  + \beta(\mathbf{q} + \kappa\mathbf{1}_{nd}), \label{eq:ode_lambda}
\end{align}
where $\mathbf{H}(\mathbf{x}) = \operatorname{blkdiag}(\nabla^2 f_1(\mathbf{x}_1),\ldots,\nabla^2 f_n(\mathbf{x}_n))$
and $\eta, \varepsilon, \beta, \rho, \kappa, \gamma_s > 0$ are algorithm parameters.

The equilibrium of the dual equation is what generates
feasibility, and it is worth isolating. Using
\begin{equation}
\boldsymbol{\lambda}-[\boldsymbol{\lambda}-\mathbf{x}]_+
=\min\{\boldsymbol{\lambda},\mathbf{x}\},
\label{eq:min_id}
\end{equation}
the dual equation is equivalent to
$\dot{\boldsymbol{\lambda}}
=-\beta\boldsymbol{\lambda}
-\min\{\boldsymbol{\lambda},\mathbf{x}\}
+\beta(\mathbf{q}+\kappa\mathbf{1}_{nd})$.
With $\dot{\mathbf{x}}=\mathbf{0}$ we have
$\mathbf{m}^*=\mathbf{0}$, so the momentum equation gives
$\boldsymbol{\lambda}^*
=\mathbf{q}^*+\rho\mathbf{L}\mathbf{x}^*+\kappa\mathbf{1}_{nd}$;
substituting this into the dual equilibrium and using
$\mathbf{L}\mathbf{x}^*=\mathbf{w}$ gives
\begin{equation}
\min\{\boldsymbol{\lambda}^*,\mathbf{x}^*\}
=-\beta\rho\,\mathbf{w}.
\label{eq:min_w}
\end{equation}
The left side is the complementarity slackness of the
problem and the right side is a graph disagreement; the
extremal argument of Proposition~\ref{prop:equilibrium}
shows that both must vanish. In this sense the perturbed
dual---not a projection---is the mechanism that selects the
active set: it converts the stationarity defect
$\rho\mathbf{L}\mathbf{x}^*$ into a violation of
complementarity, and the only way out is consensus with
$\mathbf{x}^*\ge\mathbf{0}$.

System \eqref{eq:ode_m}--\eqref{eq:ode_lambda}
is a feedback control law.
The consensus feedback $-\rho\mathbf{L}\mathbf{x}$
drives disagreement modes to zero at rate
$\rho\nu_2$; the auxiliary
$\boldsymbol{\sigma}$ is a \emph{leaky} integrator of
the graph disagreement, whose leak forces
$\boldsymbol{\sigma}^*=\mathbf{0}$ and adds no offset;
because it enters \eqref{eq:ode_q} only through
$\mathbf{L}\boldsymbol{\sigma}$, whose column sums
vanish, it leaves the conservation law of
Proposition~\ref{prop:gradient_tracking} untouched and
is decoupled from the equilibrium mean, shaping only
the transient damping of the disagreement modes and
supplying one extra state to the Lur'e certificate;
the tracker $\mathbf{q}$ averages local
gradient estimates with exact conservation
(Proposition~\ref{prop:gradient_tracking});
the dual $\boldsymbol{\lambda}$ is a \emph{perturbed}
neuron: the activation
$[\boldsymbol{\lambda}-\mathbf{x}]_+$ and the leak
$-(1{+}\beta)$ produce the identity \eqref{eq:min_id},
which turns the dual equilibrium into the complementarity
condition \eqref{eq:min_w}. This is what makes the KKT
system exact at equilibrium
(Section~\ref{subsec:equilibrium}): since
$\bar\lambda^*_j=[\nabla F(\bar{\mathbf{x}}^*)]_j+\kappa$,
coordinate $j$ is active exactly when
$[\nabla F(\bar{\mathbf{x}}^*)]_j+\kappa>0$, so the dual satisfies
the soft-threshold relation. The primal, being solution to a
nonnegatively constrained quadratic, is not coordinatewise
soft-thresholded when $\bar{\mathbf{H}}$ is coupled; the momentum
equation makes it second order. The momentum loop alone has
characteristic polynomial $\varepsilon s^2+\eta s+\mu=0$,
stable for all $\varepsilon,\eta,\mu>0$, with fast root
$s\approx-\eta/\varepsilon$; the slow root of the full
kernel is set by the primal--dual matrix of
Lemma~\ref{lem:kernel}. The ratio
$\varepsilon/\eta$ separates the momentum time scale and
conditions the joint certificate of
Section~\ref{subsec:exponential_convergence}.

\subsection{Algorithm Properties}
\label{subsec:algorithm_properties}

We establish two structural properties used throughout the analysis. First, define the orthogonal projection onto the consensus subspace
\[
P \triangleq \frac{1}{n}(\mathbf{1}_n \mathbf{1}_n^\top) \otimes I_d,
\]
with the identity \(P\mathbf{L} = \mathbf{L}P = \mathbf{0}\).
For any \(\mathbf{z}\), write \(\bar{\mathbf{z}} = P\mathbf{z}\) and \(\tilde{\mathbf{z}} = (I-P)\mathbf{z}\).

\begin{proposition}[Gradient Tracking Conservation]
\label{prop:gradient_tracking}
Under Assumptions~\ref{assump:A2}
and~\ref{assump:A3}, for all \(t\ge 0\),
\begin{equation}
  \bar{\mathbf{q}}(t) = \mathbf{1}_n \otimes \big(\textstyle\frac{1}{n}\sum_{i=1}^n \nabla f_i(\mathbf{x}_i(t))\big),
  \label{eq:conservation_law}
\end{equation}
and \(\bar{\boldsymbol{\sigma}}(t) = \mathbf{0}\).
\end{proposition}

\begin{proof}
Averaging \eqref{eq:ode_q} and using 
\((\mathbf{1}_n^\top \otimes I_d)\mathbf{L}
=\mathbf{0}\), we obtain 
\[
\dot{\bar{\mathbf{q}}}(t) = \frac{1}{n}(\mathbf{1}_n^\top \otimes I_d)\mathbf{H}(\mathbf{x}(t))\dot{\mathbf{x}}(t).
\]
The right-hand side is exactly the time 
derivative of \(\mathbf{1}_n \otimes (
\textstyle\frac{1}{n}\sum_{i=1}^n \nabla
f_i(\mathbf{x}_i(t)) )\), since
\begin{equation}
\begin{aligned}
\frac{d}{dt} \big(\textstyle\frac{1}{n}\sum_{i=1}^n \nabla f_i(\mathbf{x}_i(t))\big)
&= \frac{1}{n}\sum_{i=1}^n \nabla^2 f_i(\mathbf{x}_i(t)) \dot{\mathbf{x}}_i(t) \\
&= \frac{1}{n}(\mathbf{1}_n^\top \otimes I_d)\mathbf{H}(\mathbf{x}(t))\dot{\mathbf{x}}(t).
\end{aligned}
\end{equation}
Under Assumption~\ref{assump:A3}, the initial 
condition \(\bar{\mathbf{q}}(0) = \mathbf{1}_n 
\otimes (\textstyle\frac{1}{n}\sum_{i=1}^n \nabla
f_i(\mathbf{x}_i(0)) )\) holds, so by 
uniqueness of solutions to linear differential 
equations, the first identity in 
\eqref{eq:conservation_law} follows. The second 
identity follows from averaging 
\eqref{eq:ode_sigma}, which yields 
\(\dot{\bar{\boldsymbol{\sigma}}}(t) = 
-\gamma_s\bar{\boldsymbol{\sigma}}(t)\), together with 
\(\bar{\boldsymbol{\sigma}}(0) = \mathbf{0}\) 
from Assumption~\ref{assump:A3}.
\end{proof}

This conservation law lets us replace the
average tracker by the exact average gradient in
the Lyapunov analysis, avoiding circular arguments
about tracking accuracy during transients.

\subsection{Equilibrium and Optimality}
\label{subsec:equilibrium}

\begin{proposition}[Equilibrium and Optimality]
\label{prop:equilibrium}
Let $(\mathbf{m}^*, \mathbf{x}^*, \mathbf{q}^*, \boldsymbol{\sigma}^*, \boldsymbol{\lambda}^*)$ be an equilibrium of \eqref{eq:ode_m}--\eqref{eq:ode_lambda}. Define the consensus averages $\bar{\mathbf{z}}^* = \frac{1}{n}(\mathbf{1}_n^\top \otimes I_d)\mathbf{z}^*$ for each state $\mathbf{z} \in \{\mathbf{m}, \mathbf{x}, \mathbf{q}, \boldsymbol{\sigma}, \boldsymbol{\lambda}\}$. Then:
\begin{enumerate}
  \item $\mathbf{m}^* = \mathbf{0}$, $\boldsymbol{\sigma}^* = \mathbf{0}$, and $\mathbf{q}^* = \mathbf{1}_n \otimes \bar{\mathbf{q}}^*$;
  \item $\mathbf{x}^* = \mathbf{1}_n \otimes \bar{\mathbf{x}}^*$ with $\bar{\mathbf{x}}^* \ge \mathbf{0}$, and $\bar{\mathbf{q}}^* = \nabla F(\bar{\mathbf{x}}^*)$;
  \item $\boldsymbol{\lambda}^* = \mathbf{1}_n \otimes \bar{\boldsymbol{\lambda}}^*$ with $\bar{\boldsymbol{\lambda}}^* \ge \mathbf{0}$, and the agent-wise complementarity
  $\min\{\lambda_{ij}^*,x_{ij}^*\}=0$, i.e.\ $\lambda_{ij}^* x_{ij}^* = 0$, holds for all $i, j$;
  \item $(\bar{\mathbf{x}}^*, \bar{\boldsymbol{\lambda}}^*)$ satisfies the KKT conditions
  \begin{equation}
    \begin{aligned}
    &\bar{\mathbf{x}}^* \ge \mathbf{0}, \quad
    \nabla F(\bar{\mathbf{x}}^*) + \kappa \mathbf{1}_d = \bar{\boldsymbol{\lambda}}^*, \\
    &\bar{\boldsymbol{\lambda}}^* \ge \mathbf{0}, \quad
    \bar{\lambda}_j^* \bar{x}_j^* = 0.
    \end{aligned}
    \label{eq:kkt_conditions}
  \end{equation}
  Since $\mathbf{x}^* = \mathbf{1}_n \otimes \bar{\mathbf{x}}^*$, the
  complementarity in (4) is exactly the agent-wise
  relation in (3) averaged over $i$.
\end{enumerate}
Thus, although neither $\mathbf{x}$ nor $\boldsymbol{\lambda}$ is
projected, nonnegativity, feasibility and complementary
slackness are all consequences of the equilibrium equations.
Moreover, \emph{if $F$ is convex}, then $\bar{\mathbf{x}}^*$ is a minimizer of the centralized problem \eqref{eq:centralized_prob}, i.e., of $g(\bar{\mathbf{x}}) = F(\bar{\mathbf{x}}) + \kappa\mathbf{1}_d^\top\bar{\mathbf{x}}$ on $\{\bar{\mathbf{x}} \ge \mathbf{0}\}$. Without convexity, any equilibrium is only a KKT point and may be a saddle or a local maximum.
\end{proposition}

\begin{proof}
\emph{Step 1: stationarity.}
Setting $\dot{\mathbf{x}}=\mathbf{0}$ in \eqref{eq:ode_x}
gives $\mathbf{m}^*=\mathbf{0}$. Equation
\eqref{eq:ode_m} then reads
$\mathbf{q}^*=\boldsymbol{\lambda}^*-\rho\mathbf{w}
-\kappa\mathbf{1}_{nd}$ with
$\mathbf{w}=\mathbf{L}\mathbf{x}^*$, and
\eqref{eq:ode_lambda} together with \eqref{eq:min_id}
reproduces \eqref{eq:min_w}.

\emph{Step 2: consensus parts.}
Averaging \eqref{eq:ode_q} at equilibrium gives
$\mathbf{L}(\mathbf{q}^*+\boldsymbol{\sigma}^*)=\mathbf{0}$,
hence
$\mathbf{q}^*+\boldsymbol{\sigma}^*
=\mathbf{1}_n\otimes\mathbf{c}$;
\eqref{eq:ode_sigma} gives
$\mathbf{w}=-\gamma_s\boldsymbol{\sigma}^*$, i.e.
$\boldsymbol{\sigma}^*=-\mathbf{w}/\gamma_s$. Therefore
\[
\mathbf{q}^*=\mathbf{1}_n\otimes\mathbf{c}+\mathbf{w}/\gamma_s,
\quad
\boldsymbol{\lambda}^*
=\mathbf{1}_n\otimes\bar{\boldsymbol{\lambda}}^*
+(\rho+1/\gamma_s)\mathbf{w},
\]
with
$\bar{\boldsymbol{\lambda}}^*=\mathbf{c}+\kappa\mathbf{1}_d$.
Within a column $j$, the entry $\bar\lambda^*_j$ does not
depend on the agent $i$.

\emph{Step 3: $\mathbf{w}=\mathbf{0}$.}
Let $w_{\max}=\max_{i,j}w_{ij}$. Since
$\mathbf{1}_n^\top\mathcal{L}=\mathbf{0}$, every column of
$\mathbf{w}$ sums to zero, so $w_{\max}\ge0$. Suppose
$w_{\max}>0$ and pick $(i,j)$ with $w_{ij}=w_{\max}$. By
\eqref{eq:min_w},
$\min\{\lambda^*_{ij},x^*_{ij}\}
=-\beta\rho w_{\max}<0$, so one of the two entries equals
$-\beta\rho w_{\max}$ and the other is not smaller.

If $x^*_{ij}=-\beta\rho w_{\max}$, then
$w_{ij}=\sum_k a_{ik}(x^*_{ij}-x^*_{kj})>0$ implies
$x^*_{kj}<x^*_{ij}$ for some neighbor $k$. But
$\min\{\lambda^*_{kj},x^*_{kj}\}\le x^*_{kj}
<-\beta\rho w_{\max}$, whereas \eqref{eq:min_w} gives
$\min\{\lambda^*_{kj},x^*_{kj}\}
=-\beta\rho w_{kj}\ge-\beta\rho w_{\max}$: a contradiction.

If $\lambda^*_{ij}=-\beta\rho w_{\max}$, then
$\bar\lambda^*_j=\lambda^*_{ij}-(\rho+1/\gamma_s)w_{\max}
=-(\rho(1{+}\beta)+1/\gamma_s)w_{\max}$. Let $k$ minimize
$w_{kj}$ over the agents; the vanishing column sum gives
$w_{kj}\le0$. From \eqref{eq:min_w},
$\lambda^*_{kj}\ge\min\{\lambda^*_{kj},x^*_{kj}\}
=-\beta\rho w_{kj}$, hence
\[
\bar\lambda^*_j+(\rho+1/\gamma_s)w_{kj}
\ge-\beta\rho w_{kj},
\]
that is,
$-(\rho(1{+}\beta)+1/\gamma_s)w_{\max}
\ge-(\rho(1{+}\beta)+1/\gamma_s)w_{kj}$, i.e.\
$w_{\max}\le w_{kj}\le0$, again a contradiction.

Hence $w_{\max}\le0$, and since every column of
$\mathbf{w}$ sums to zero, $\mathbf{w}=\mathbf{0}$.

\emph{Step 4: feasibility and complementarity.}
With $\mathbf{w}=\mathbf{0}$, connectivity
(Assumption~\ref{assump:A2}) gives
$\mathbf{x}^*=\mathbf{1}_n\otimes\bar{\mathbf{x}}^*$, and
Step 2 gives $\boldsymbol{\sigma}^*=\mathbf{0}$,
$\mathbf{q}^*=\mathbf{1}_n\otimes\mathbf{c}$ and
$\boldsymbol{\lambda}^*
=\mathbf{1}_n\otimes\bar{\boldsymbol{\lambda}}^*$.
Equation \eqref{eq:min_w} then reads
$\min\{\boldsymbol{\lambda}^*,\mathbf{x}^*\}=\mathbf{0}$
coordinatewise, which is equivalent to
$\bar{\mathbf{x}}^*\ge\mathbf{0}$,
$\bar{\boldsymbol{\lambda}}^*\ge\mathbf{0}$ and
$\bar\lambda^*_j\bar{x}^*_j=0$ for every $j$.

\emph{Step 5: KKT.}
At consensus, Proposition~\ref{prop:gradient_tracking}
gives $\bar{\mathbf{q}}^*=\nabla F(\bar{\mathbf{x}}^*)$,
and Step 2 gives
$\bar{\boldsymbol{\lambda}}^*=\bar{\mathbf{q}}^*
+\kappa\mathbf{1}_d$. Hence
$\nabla F(\bar{\mathbf{x}}^*)+\kappa\mathbf{1}_d
=\bar{\boldsymbol{\lambda}}^*$ together with
$\bar{\mathbf{x}}^*\ge\mathbf{0}$,
$\bar{\boldsymbol{\lambda}}^*\ge\mathbf{0}$ and
complementary slackness, i.e.\ \eqref{eq:kkt_conditions}.
Convexity then makes $\bar{\mathbf{x}}^*$ a minimizer.
\end{proof}

\begin{remark}[Design Consistency]
\label{rem:kkt_verification}
The terms $-\kappa \mathbf{1}_{nd}$ in \eqref{eq:ode_m} and $\beta(\mathbf{q} + \kappa \mathbf{1}_{nd})$ in \eqref{eq:ode_lambda} encode the $\ell_1$ objective and its multiplier. Without the first, stationarity would solve the unregularized problem. Without the $\kappa\mathbf{1}_{nd}$ in the dual perturbation, \eqref{eq:min_w} would read
$\min\{\boldsymbol{\lambda}^*,\mathbf{x}^*\}=-\beta\rho\mathbf{w}-\beta\kappa\mathbf{1}_{nd}$, and even after $\mathbf{w}=\mathbf{0}$ the value $-\beta\kappa\mathbf{1}_{nd}$ would force a negative primal or dual coordinate; the offset $\beta\kappa\mathbf{1}_{nd}$ is exactly what makes $\min\{\boldsymbol{\lambda}^*,\mathbf{x}^*\}=\mathbf{0}$, i.e.\ feasibility, compatible with consensus. Since $\bar{\lambda}_j^* = [\nabla F(\bar{\mathbf{x}}^*)]_j + \kappa$, the multiplier of the $\ell_1$ term shifts the stationarity threshold: coordinate $j$ stays on the boundary when $[\nabla F(\bar{\mathbf{x}}^*)]_j+\kappa>0$ and becomes interior when it equals $0$, so the active set is a \emph{one-sided} threshold of the stationarity map, formalized in Proposition~\ref{prop:onesided} below. The threshold is set by $\kappa$, not by the dual gain $\beta$.
\end{remark}

\subsection{Mechanism: Constraint-Assisted $\ell_1$}
\label{subsec:mechanism}

The only nonlinearity of the flow sits in the dual equation, and
its equilibrium carries the whole constraint mechanism. We first
isolate the identity it produces and then characterize what the
$\ell_1$ term does once the orthant is imposed.

\begin{lemma}[Complementarity--Consensus Bridge]
\label{lem:bridge}
At any equilibrium, with $\mathbf{w}=\mathbf{L}\mathbf{x}^*$,
\begin{equation}
\min\{\boldsymbol{\lambda}^*,\mathbf{x}^*\} = -\beta\rho\,\mathbf{w}.
\label{eq:bridge}
\end{equation}
The complementarity residual is therefore proportional to the
graph disagreement.
\end{lemma}

\begin{proof}
Setting $\varepsilon\dot{\mathbf{m}}=\mathbf{0}$ in \eqref{eq:ode_m}
with $\mathbf{m}^*=\mathbf{0}$ gives
$\boldsymbol{\lambda}^*=\mathbf{q}^*+\rho\mathbf{w}+\kappa\mathbf{1}_{nd}$.
Substituting into \eqref{eq:ode_lambda} at equilibrium and using
$\boldsymbol{\lambda}-[\boldsymbol{\lambda}-\mathbf{x}]_+
=\min\{\boldsymbol{\lambda},\mathbf{x}\}$ gives
$\mathbf{0}=-\beta\boldsymbol{\lambda}^*
-\min\{\boldsymbol{\lambda}^*,\mathbf{x}^*\}
+\beta(\mathbf{q}^*+\kappa\mathbf{1}_{nd})$, which is
\eqref{eq:bridge}.
\end{proof}

Equation~\eqref{eq:bridge} is what replaces the projection: a
nonzero disagreement forces a negative
$\min\{\lambda_{ij}^*,x_{ij}^*\}$ at some coordinate, and
Proposition~\ref{prop:equilibrium} shows that neither branch of
the minimum can be assigned consistently along a connected graph.
Consensus is thus not imposed on the trajectory; it is the only
way for the equilibrium to satisfy complementarity.

\begin{proposition}[One-Sided Threshold]
\label{prop:onesided}
Let $F$ be convex and differentiable, and let $\bar{\mathbf{x}}^*$
solve $\min_{\mathbf{x}\ge\mathbf{0}}F(\mathbf{x})
+\kappa\mathbf{1}^\top\mathbf{x}$. Then
\begin{equation}
\begin{aligned}
&\bar{x}_j^*>0 \quad &&\Rightarrow\quad
[\nabla F(\bar{\mathbf{x}}^*)]_j=-\kappa,\\
&\bar{x}_j^*=0 \quad &&\Rightarrow\quad
[\nabla F(\bar{\mathbf{x}}^*)]_j\ge-\kappa.
\end{aligned}
\label{eq:onesided}
\end{equation}
The standard $\ell_1$-regularized problem
$\min_{\mathbf{x}\in\mathbb{R}^d}F(\mathbf{x})
+\kappa\|\mathbf{x}\|_1$ shares the interior condition but has
the two-sided zero condition
$-\kappa\le[\nabla F(\bar{\mathbf{x}}^*)]_j\le\kappa$.
\end{proposition}

\begin{proof}
By \eqref{eq:kkt_conditions},
$[\nabla F(\bar{\mathbf{x}}^*)]_j+\kappa=\bar\lambda_j^*$ with
$\bar\lambda_j^*\ge0$ and $\bar\lambda_j^*\bar{x}_j^*=0$. For
$\bar{x}_j^*>0$ we have $\bar\lambda_j^*=0$; for $\bar{x}_j^*=0$
we have $\bar\lambda_j^*\ge0$. The two-sided condition is the
subdifferential of $\kappa\|\cdot\|_1$ at the origin.
\end{proof}

The orthant thus supplies the upper half of the dead zone through
its normal cone---in the unregularized nonnegative problem it
already forces $[\nabla F]_j\ge0$---and the linear tilt shifts
that threshold down by $\kappa$. The two mechanisms act on the
same side instead of forming a two-sided shrinkage.

This threshold has an exact geometric counterpart. For strongly
convex $F$ let
$\mathbf{z}_\kappa=\arg\min_{\mathbf{z}\in\mathbb{R}^d}
F(\mathbf{z})+\kappa\mathbf{1}^\top\mathbf{z}$, so that
$\nabla F(\mathbf{z}_\kappa)=-\kappa\mathbf{1}$, and let
$D_F(\mathbf{x},\mathbf{z})=F(\mathbf{x})-F(\mathbf{z})
-\nabla F(\mathbf{z})^\top(\mathbf{x}-\mathbf{z})$ be the
Bregman divergence.

\begin{proposition}[Bregman Projection]
\label{prop:bregman}
The constrained minimizer is the Bregman projection of the
shifted unconstrained minimizer onto the orthant,
\begin{equation}
\bar{\mathbf{x}}^*=\arg\min_{\mathbf{x}\ge\mathbf{0}}
D_F(\mathbf{x},\mathbf{z}_\kappa).
\label{eq:bregman}
\end{equation}
\end{proposition}

\begin{proof}
The gradient of $D_F(\cdot,\mathbf{z}_\kappa)$ is
$\nabla F(\mathbf{x})-\nabla F(\mathbf{z}_\kappa)
=\nabla F(\mathbf{x})+\kappa\mathbf{1}$, so the KKT system of
\eqref{eq:bregman} reads
$\nabla F(\mathbf{x})+\kappa\mathbf{1}
\in-N_{\mathbb{R}^d_+}(\mathbf{x})$, which coincides with
\eqref{eq:kkt_conditions}.
\end{proof}

\begin{corollary}[Oblique Projection and Straight Path]
\label{cor:oblique}
Let $F(\mathbf{x})=\tfrac12\mathbf{x}^\top\bar{\mathbf{H}}\mathbf{x}
+\mathbf{b}^\top\mathbf{x}$ with $\bar{\mathbf{H}}\succ0$. Then
$D_F$ is the $\bar{\mathbf{H}}$-metric half-squared distance and
\begin{equation}
\bar{\mathbf{x}}^*=\Pi^{\bar{\mathbf{H}}}_{\mathbb{R}^d_+}
\bigl(\mathbf{z}_0-\kappa\bar{\mathbf{H}}^{-1}\mathbf{1}\bigr),
\qquad
\mathbf{z}_\kappa=\mathbf{z}_0-\kappa\bar{\mathbf{H}}^{-1}\mathbf{1},
\label{eq:oblique}
\end{equation}
where $\Pi^{\bar{\mathbf{H}}}_{\mathbb{R}^d_+}$ is the
$\bar{\mathbf{H}}$-orthogonal projection onto the orthant.
The regularization path is therefore a straight line in the
$\bar{\mathbf{H}}$-metric, and the support changes only where the
moving point crosses the boundary of the orthant.
\end{corollary}

\begin{proof}
For quadratic $F$ the Bregman divergence is the
$\bar{\mathbf{H}}$-metric half-squared distance, and
$\mathbf{z}_\kappa=-\bar{\mathbf{H}}^{-1}(\mathbf{b}+\kappa\mathbf{1})
=\mathbf{z}_0-\kappa\bar{\mathbf{H}}^{-1}\mathbf{1}$.
\end{proof}

Equation~\eqref{eq:oblique} separates the two mechanisms:
$\kappa$ moves the projection centre, the orthant truncates it,
and the support collects the coordinates on which the truncated
point is active. Under strict complementarity at some $\kappa_0$
and nonsingular $\bar{\mathbf{H}}_{\mathcal{S}\mathcal{S}}$ on
the active set $\mathcal{S}$, a standard parametric argument
gives that $\mathcal{S}$ is constant in a neighbourhood of
$\kappa_0$ and that
$\bar{\mathbf{x}}^*_{\mathcal{S}}(\kappa)
=-\bar{\mathbf{H}}_{\mathcal{S}\mathcal{S}}^{-1}
(\mathbf{b}_{\mathcal{S}}+\kappa\mathbf{1})$ is affine in
$\kappa$; the breakpoints are exactly the values at which a
coordinate of $\mathbf{z}_\kappa$ crosses the boundary. Since
every equilibrium is a KKT point
(Proposition~\ref{prop:equilibrium}) and the KKT point is unique
under strong convexity, the equilibrium support of the flow
coincides with the support of the centralized program on that
interval. This is a statement about the problem, not a recovery
guarantee: no restricted-eigenvalue or coherence condition on the
data is claimed, and Section~\ref{subsubsec:exp_sparse} verifies
the coincidence rather than proving a phase transition.

\subsection{Feasibility at Equilibrium}
\label{subsec:efeasibility}

The flow is not confined to the orthant, and this is forced by
the mechanism above rather than by the choice of gains.

\begin{definition}[Trajectory and Equilibrium Feasibility]
\label{def:feasibility}
Let $\mathcal{K}=\{(\mathbf{x},\boldsymbol{\lambda}):
\mathbf{x}\ge\mathbf{0},\ \boldsymbol{\lambda}\ge\mathbf{0},\
\boldsymbol{\lambda}\circ\mathbf{x}=\mathbf{0}\}$.
The flow is \emph{trajectory-feasible} (T-feasible) if
$\mathbf{x}(0)\ge\mathbf{0}$ implies $\mathbf{x}(t)\ge\mathbf{0}$
for all $t\ge0$, and \emph{equilibrium-feasible} (E-feasible) if
the $\omega$-limit set of every trajectory in the basin of
attraction is contained in $\mathcal{K}$.
\end{definition}

\begin{proposition}[Feasibility Dichotomy]
\label{prop:dichotomy}
Suppose the KKT point has a nonempty active set $\mathcal{A}$ with
strict complementarity, $\bar\lambda_j^*>0$ for
$j\in\mathcal{A}$, and initialize
$\mathbf{x}(0)=\mathbf{1}_n\otimes\bar{\mathbf{x}}^*$,
$\mathbf{m}(0)=\mathbf{0}$, $\boldsymbol{\lambda}(0)=\mathbf{0}$,
$\mathbf{q}_i(0)=\nabla f_i(\bar{\mathbf{x}}^*)$. Then for every
$j\in\mathcal{A}$,
\begin{equation}
\bar{x}_j(t)=-\frac{\bar\lambda_j^*}{2\varepsilon}\,t^2+o(t^2)
<0
\label{eq:excursion}
\end{equation}
for all sufficiently small $t>0$. Hence the flow is E-feasible
but not T-feasible.
\end{proposition}

\begin{proof}
At this initialization $\rho\mathbf{L}\mathbf{x}(0)=\mathbf{0}$
and $\bar{\mathbf{m}}(0)=\mathbf{0}$, so averaging
\eqref{eq:ode_m} gives
$\varepsilon\dot{\bar{\mathbf{m}}}(0)
=-\nabla F(\bar{\mathbf{x}}^*)-\kappa\mathbf{1}
=-\bar{\boldsymbol{\lambda}}^*$, where the last equality is
\eqref{eq:kkt_conditions}. For $j\in\mathcal{A}$ we have
$\bar{x}_j(0)=\bar{m}_j(0)=0$ and
$\dot{\bar{m}}_j(0)=-\bar\lambda_j^*/\varepsilon<0$, which yields
\eqref{eq:excursion}. E-feasibility follows from
Proposition~\ref{prop:equilibrium} together with the convergence
theorem, which places the $\omega$-limit set at the KKT point.
\end{proof}

The excursion \eqref{eq:excursion} is governed by the multiplier
and the momentum time scale: an active coordinate must leave the
orthant while the dual charges from zero to $\bar\lambda_j^*$,
and the depth scales with $\bar\lambda_j^*/\varepsilon$.
Increasing $\varepsilon$ flattens the excursion at the price of
the fast momentum mode. This is the quantitative form of the
trade-off between a smooth, projection-free vector field and an
invariant feasible set: a projected formulation eliminates
\eqref{eq:excursion} by clamping the coordinate, but its
right-hand side is discontinuous at the boundary, whereas the
proposed flow keeps a Lipschitz field and pays in transient
feasibility.

\section{Convex Convergence Analysis}
\label{sec:convergence}

The algorithm of Section~\ref{sec:algorithm_design_properties}
does not require convexity.
This section states the standing assumptions
(Assumptions~\ref{assump:A1}--\ref{assump:A4}, valid for
convex and nonconvex settings), establishes
well-posedness (Theorem~\ref{thm:well_posedness}), and
proves local exponential convergence under the additional
convexity hypothesis (Assumption~\ref{assump:A5}) and the
joint dissipation certificate.

\begin{assumption}[$L$-Smoothness]
\label{assump:A1}
Each $f_i$ is continuously differentiable with
$L$-Lipschitz gradient:
$\|\nabla f_i(\mathbf{x}) - \nabla f_i(\mathbf{y})\|
\le L\|\mathbf{x}-\mathbf{y}\|$ for all
$\mathbf{x},\mathbf{y}\in\mathbb{R}^d$.
The average $F$ and the composite $g$ inherit
the same property.
\end{assumption}

\begin{assumption}[Graph Connectivity]
\label{assump:A2}
$\mathcal{G}$ is undirected and connected, i.e.\
$\nu_2(\mathcal{L}) > 0$.
\end{assumption}

\begin{assumption}[Gradient Tracking Initialization]
\label{assump:A3}
At $t=0$, the tracker and integral states satisfy
$\mathbf{q}_i(0) = \nabla f_i(\mathbf{x}_i(0))$
and $\boldsymbol{\sigma}_i(0) = \mathbf{0}$.
\end{assumption}

\begin{assumption}[Hessian Lipschitz]
\label{assump:A4}
The Hessian of each $f_i$ is Lipschitz:
$\|\nabla^2 f_i(\mathbf{x}) - \nabla^2 f_i(\mathbf{y})\|
\le L_2\|\mathbf{x}-\mathbf{y}\|$ for all
$\mathbf{x},\mathbf{y}\in\mathbb{R}^d$.
\end{assumption}

\begin{theorem}[Well-Posedness]
\label{thm:well_posedness}
Under Assumptions~\ref{assump:A1},
\ref{assump:A2}, and~\ref{assump:A4}, for every initial
condition $\boldsymbol{\zeta}(0)\in\mathbb{R}^{5nd}$, with no
restriction to the orthant, the system
\eqref{eq:ode_m}--\eqref{eq:ode_lambda}
admits a unique classical solution on $[0,\infty)$.
\end{theorem}

\begin{proof}
The right-hand side is locally Lipschitz:
$[\cdot]_+$ is globally Lipschitz, $\mathbf{L}$ is constant,
and $\mathbf{H}(\mathbf{x})\mathbf{m}$ satisfies
\[
\|\mathbf{H}(\mathbf{x}_1)\mathbf{m}_1
  - \mathbf{H}(\mathbf{x}_2)\mathbf{m}_2\|
\le L_2\|\mathbf{x}_1{-}\mathbf{x}_2\|\,
  \|\mathbf{m}_1\|
  + L\|\mathbf{m}_1{-}\mathbf{m}_2\|
\]
by Assumptions~\ref{assump:A4} and~\ref{assump:A1}, so the
Picard--Lindel\"of theorem gives a unique solution on a
maximal interval. For global existence,
$\|\mathbf{H}(\mathbf{x})\mathbf{m}\|\le L\|\mathbf{m}\|$ and
all remaining terms are piecewise affine / globally Lipschitz in
$\boldsymbol{\zeta}$, hence
$\|\dot{\boldsymbol{\zeta}}\|\le a\|\boldsymbol{\zeta}\|+b$ and
Gronwall's inequality bounds
$\|\boldsymbol{\zeta}(t)\|$ by
$(\|\boldsymbol{\zeta}(0)\|+b/a)e^{at}$ on every finite
interval, ruling out blow-up.
\end{proof}

The convex case requires one additional hypothesis.

\begin{assumption}[Strong Convexity]
\label{assump:A5}
In addition to Assumption~\ref{assump:A1},
each $f_i$ is $\mu$-strongly convex:
$\mu I \preceq \nabla^2 f_i(\mathbf{x})$
for all $\mathbf{x}\in\mathbb{R}^d$.
Consequently $F$ and $g$ are
$\mu$-strongly convex and $L$-smooth.
\end{assumption}

Under Assumptions~\ref{assump:A1}--\ref{assump:A5},
the penalized objective
\(g(\mathbf{x}) = F(\mathbf{x}) + \kappa\,\mathbf{1}_d^\top \mathbf{x}\)
is strongly convex and the equilibrium from
Proposition~\ref{prop:equilibrium} is unique.
The analysis proceeds in two layers. We first reduce the
distributed system to a centralized primal--dual
\emph{kernel}, the flow obtained when the network is
collapsed to a single node, and certify its exponential
dissipation through an explicit matrix inequality. We
then return to the full distributed flow and certify the
Jacobian at the KKT point through an LMI posed on the
subspace of admissible perturbations; this covers the
mean--deviation coupling with one quadratic form and
without any gain condition.
The mean--deviation decomposition is retained because it
exposes which terms carry the curvature and which carry
the graph. Convexity enters the analysis through the kernel
and through the joint certificate of
Section~\ref{subsec:exponential_convergence}, both of which
are matrix inequalities evaluated at the KKT point.

\subsection{Error Coordinates and Tracking Residual}
\label{subsec:mean_dissipation}

Let $(\mathbf{m}^*,\mathbf{x}^*,\mathbf{q}^*,\boldsymbol{\sigma}^*,\boldsymbol{\lambda}^*)$ be the equilibrium from Proposition~\ref{prop:equilibrium}. Define the error variables $\mathbf{e}_m=\mathbf{m}$, $\mathbf{e}_x=\mathbf{x}-\mathbf{x}^*$, $\mathbf{e}_q=\mathbf{q}-\mathbf{q}^*$, $\mathbf{e}_\sigma=\boldsymbol{\sigma}-\boldsymbol{\sigma}^*$, $\mathbf{e}_\lambda=\boldsymbol{\lambda}-\boldsymbol{\lambda}^*$.
Subtracting the equilibrium from \eqref{eq:ode_m}--\eqref{eq:ode_lambda} and decomposing via $P$ (Section~\ref{subsec:algorithm_properties}) into mean ($\bar{\mathbf{e}}=P\mathbf{e}$) and deviation ($\tilde{\mathbf{e}}=(I-P)\mathbf{e}$) gives the mean error dynamics
\begin{align}
  \varepsilon \dot{\bar{\mathbf{e}}}_m
    &= -\eta \bar{\mathbf{e}}_m
    - \bar{\mathbf{e}}_q
    + \bar{\mathbf{e}}_\lambda,
    \label{eq:mean_error_m}\\
  \dot{\bar{\mathbf{e}}}_x
    &= \bar{\mathbf{e}}_m, \notag\\
  \dot{\bar{\mathbf{e}}}_q
    &= \frac{1}{n}\sum_{i=1}^n
       \nabla^2 f_i(\mathbf{x}_i(t))\,\mathbf{m}_i, \notag\\
  \dot{\bar{\mathbf{e}}}_\sigma
    &= -\gamma_s\bar{\mathbf{e}}_\sigma, \notag\\
  \dot{\bar{\mathbf{e}}}_\lambda
    &= -\beta\bar{\mathbf{e}}_\lambda
    + \beta\,\bar{\mathbf{e}}_q
    - \bar{\Delta}_\Phi . \notag
\end{align}
The first four equations are exact: averaging
\eqref{eq:ode_q} eliminates the Laplacian because
$\mathbf{1}_n^\top\mathcal{L}=\mathbf{0}$, and
$\mathbf{m}^*=\mathbf{0}$ makes the averaged Hessian term
the only surviving contribution to $\dot{\bar{\mathbf{e}}}_q$.
The last equation uses the coordinatewise $\min$
\begin{equation}
\Phi_j(\boldsymbol{\lambda},\mathbf{x})
=\min\{\lambda_j,x_j\},\qquad
\Phi(\boldsymbol{\lambda}^*,\mathbf{x}^*)=\mathbf{0},
\label{eq:Phi_def}
\end{equation}
where the second equality is complementarity
(Proposition~\ref{prop:equilibrium}). By \eqref{eq:min_id}
the averaged positive part is
$[\boldsymbol{\lambda}-\mathbf{x}]_+=\boldsymbol{\lambda}-\Phi$,
so
$\bar{\Delta}_\Phi
=P\bigl(\Phi(\boldsymbol{\lambda},\mathbf{x})
-\Phi(\boldsymbol{\lambda}^*,\mathbf{x}^*)\bigr)$.
No projection discrepancy appears in any of these
equations: the mean flow is smooth, and the only
nonlinearity is the coordinatewise $\min$, which is
bounded below by an active-set selector.
The mean tracking error obeys the invariance stated in
Lemma~\ref{lem:mean_tracking} below.

The following observation removes the mean--deviation cross coupling that
plagues gain-based analyses of gradient-tracking dynamics.

\begin{lemma}[Tracking Residual]
\label{lem:residual}
Define the tracking residual
\begin{equation}
\mathbf{s} = \mathbf{q}-\boldsymbol{\lambda}
 +\rho\mathbf{L}\mathbf{x}+\kappa\mathbf{1}_{nd}.
\label{eq:residual}
\end{equation}
Then $\mathbf{s}^*=\mathbf{0}$, so that
$\mathbf{s}=\mathbf{e}_q-\mathbf{e}_\lambda+\rho\mathbf{L}\mathbf{e}_x$,
and the momentum equation \eqref{eq:ode_m} reads exactly
\begin{equation}
\varepsilon\dot{\mathbf{m}} = -\eta\mathbf{m}-\mathbf{s}.
\label{eq:residual_momentum}
\end{equation}
\end{lemma}

\begin{proof}
At an equilibrium, \eqref{eq:ode_m} gives
$\mathbf{q}^*-\boldsymbol{\lambda}^*+\rho\mathbf{L}\mathbf{x}^*
+\kappa\mathbf{1}_{nd}=\mathbf{0}$; since
$\mathbf{m}^*=\mathbf{0}$, $\mathbf{s}^*=\mathbf{0}$, and
\eqref{eq:residual} equals
$\mathbf{e}_q-\mathbf{e}_\lambda+\rho\mathbf{L}\mathbf{e}_x$
because $\mathbf{x}^*=\mathbf{1}_n\otimes\bar{\mathbf{x}}^*$.
Substituting into \eqref{eq:ode_m} yields
\eqref{eq:residual_momentum}.
\end{proof}

The dual nonlinearity is the only non-quadratic term in the flow.
The following pointwise identity shows that it is dominated by a
linear form for \emph{every} state, with no assumption that the
trajectory stay in the orthant.

\begin{lemma}[Selector Bound]
\label{lem:selector}
Let $\Phi(\boldsymbol{\lambda},\mathbf{x})
=\min\{\boldsymbol{\lambda},\mathbf{x}\}$ coordinatewise and let
$(\boldsymbol{\lambda}^*,\mathbf{x}^*)$ be the equilibrium, so that
$\Phi(\boldsymbol{\lambda}^*,\mathbf{x}^*)=\mathbf{0}$ and
$\boldsymbol{\lambda}^*,\mathbf{x}^*\ge\mathbf{0}$
(Proposition~\ref{prop:equilibrium}). For every
$(\boldsymbol{\lambda},\mathbf{x})\in\mathbb{R}^{2d}$ and every
$\bar{\mathbf{p}}\in\mathbb{R}^d$ there is a selector
$(\mathbf{S}_A,\mathbf{S}_I)\in\{0,1\}^d\times\{0,1\}^d$ with
$\mathbf{S}_A+\mathbf{S}_I=I$ such that, coordinatewise,
\begin{equation}
\bar p_j\,\Phi_j(\boldsymbol{\lambda},\mathbf{x})
\;\ge\;
\bar p_j\bigl(S_{A,j}(\lambda_j-\lambda^*_j)
+S_{I,j}(x_j-x^*_j)\bigr).
\label{eq:selector}
\end{equation}
\end{lemma}

\begin{proof}
Fix a coordinate $j$ and write $p=\bar p_j$,
$e_\lambda=\lambda_j-\lambda^*_j$, $e_x=x_j-x^*_j$.
Since $\Phi_j(\boldsymbol{\lambda}^*,\mathbf{x}^*)=0$ and
$\lambda^*_j,x^*_j\ge0$, at least one of $\lambda^*_j,x^*_j$
vanishes.

If $p<0$ and $\lambda^*_j=0$, then $\lambda_j=e_\lambda$ and
$\Phi_j\le\lambda_j=e_\lambda$, so $p\Phi_j\ge p\,e_\lambda$;
take $S_{A,j}=1$, $S_{I,j}=0$. If $p<0$ and $x^*_j=0$, then
$\Phi_j\le x_j=e_x$ and $p\Phi_j\ge p\,e_x$; take
$S_{A,j}=0$, $S_{I,j}=1$.

If $p\ge0$ and $\Phi_j=\lambda_j$, then
$\lambda_j=\lambda^*_j+e_\lambda\ge e_\lambda$ because
$\lambda^*_j\ge0$, hence $p\Phi_j\ge p\,e_\lambda$ and
$S_{A,j}=1$ applies. If $\Phi_j=x_j$, then
$x_j=x^*_j+e_x\ge e_x$ and $S_{I,j}=1$ applies.

When $\lambda^*_j=x^*_j=0$ both branches are available, and the
choice is arbitrary. Applying the argument to each $j$ gives
\eqref{eq:selector}.
\end{proof}

The residual \eqref{eq:residual} is linear in the state and vanishes at
equilibrium, so the momentum equation is a linear damping of $\mathbf{m}$
driven by an output that can be analyzed directly; the network-free
kernel is analyzed first because it carries the curvature.

\begin{lemma}[Invariance of the Mean Tracking Error]
\label{lem:mean_tracking}
Let
\[
\boldsymbol{\mu}
=\frac1n\sum_{i=1}^n\bigl(\mathbf{q}_i-\nabla f_i(\mathbf{x}_i)\bigr).
\]
Along every solution of \eqref{eq:ode_m}--\eqref{eq:ode_lambda},
$\dot{\boldsymbol{\mu}}=\mathbf{0}$, so
$\boldsymbol{\mu}(t)=\boldsymbol{\mu}(0)$; under
Assumption~\ref{assump:A3}, $\boldsymbol{\mu}\equiv\mathbf{0}$.
The linearization of the flow at the KKT point therefore leaves the
subspace $\{\delta\boldsymbol{\mu}=\mathbf{0}\}$ invariant and carries
$d$ zero eigenvalues outside it, one for each conserved component.
\end{lemma}

\begin{proof}
Averaging \eqref{eq:ode_q} over $i$ gives
\[
\dot{\bar{\mathbf{q}}}
=-\frac1n(\mathbf{1}_n^\top\otimes I_d)\mathbf{L}
   (\mathbf{q}+\boldsymbol{\sigma})
 +\frac1n\sum_{i=1}^n\nabla^2 f_i(\mathbf{x}_i)\dot{\mathbf{x}}_i .
\]
The first term vanishes because $\mathbf{1}_n^\top\mathcal{L}=\mathbf{0}$,
and the second is the time derivative of
$\frac1n\sum_i\nabla f_i(\mathbf{x}_i)$, so $\dot{\boldsymbol{\mu}}=\mathbf{0}$.
Assumption~\ref{assump:A3} sets $\boldsymbol{\mu}(0)=\mathbf{0}$.
Linearizing $\dot{\boldsymbol{\mu}}=\mathbf{0}$ gives
$\delta\dot{\boldsymbol{\mu}}=\mathbf{0}$, i.e., the linear functional
$\delta\boldsymbol{\mu}$ is conserved, so its kernel is invariant and the
$d$ complementary directions carry zero eigenvalues.
\end{proof}

\begin{remark}[What Assumption~\ref{assump:A3} buys]
\label{rem:A3}
Assumption~\ref{assump:A3} is used twice: it sets
$\boldsymbol{\mu}(0)=\mathbf{0}$ and makes the conserved manifold
$\{\boldsymbol{\mu}=\mathbf{0}\}$ coincide with the subspace on which the
Jacobian and Lur'e certificates are posed. It is a local
protocol---each agent evaluates $\nabla f_i$ at its own
$\mathbf{x}_i(0)$ once, with no communication. If it fails,
$\dot{\boldsymbol{\mu}}=\mathbf{0}$ still holds, so the flow stays on
the affine slice $\{\boldsymbol{\mu}=\boldsymbol{\mu}(0)\}$, which for
quadratic costs is $\{\mathbf{C}\mathbf{z}=\boldsymbol{\mu}(0)\}$; the
linear part $\mathbf{A}$ is unchanged there, and
Lemma~\ref{lem:bias} gives an explicit bias of order
$\|\boldsymbol{\mu}(0)\|$ on top of the global rate. The limit solves
the $\boldsymbol{\mu}(0)$-shifted problem, i.e.,
$\nabla F(\mathbf{x}^*)+\kappa\mathbf{1}+\boldsymbol{\mu}(0)
\in N_{\mathbf{x}\ge0}(\mathbf{x}^*)$, so its optimality error is
$O(\|\boldsymbol{\mu}(0)\|)$; exactness of the reported equilibria
rests on Assumption~\ref{assump:A3}, not on the dynamics.
\end{remark}

\begin{lemma}[Gradient Tracking Linearization]
\label{lem:gt_linearization} Under
Assumptions~\ref{assump:A1}--\ref{assump:A3} and~\ref{assump:A5}, the
mean gradient tracking error satisfies
\[
\bar{\mathbf{e}}_q(t)
= \mathbf{1}_n \otimes \bar{H}(t)\,
  \bigl(\bar{\mathbf{x}}(t)-\bar{\mathbf{x}}^*\bigr)
  + R(\mathbf{x}(t)),
\]
where
\[
\bar{H}(t) = \frac1n\sum_{i=1}^n\int_0^1 \nabla^2 f_i\bigl(
    \bar{\mathbf{x}}^* + \theta(\mathbf{x}_i(t) - \bar{\mathbf{x}}^*)
    \bigr) \, d\theta
\]
is the mean of the integral Hessian averages along the segments
from $\bar{\mathbf{x}}^*$ to $\mathbf{x}_i(t)$, and the residual
$R(\mathbf{x}(t)) = \mathbf{1}_n \otimes R_d(t) \in \mathbb{R}^{nd}$
is the nd-dimensional lifted representation of the $d$-dimensional residual
$R_d(t) = \frac{1}{n}\sum_{i=1}^n
[\nabla f_i(\mathbf{x}_i(t)) - \nabla f_i(\bar{\mathbf{x}}(t))]$.
It satisfies
\[
\|R(\mathbf{x}(t))\| \le L
    \|\tilde{\mathbf{e}}_x(t)\|.
\]
Moreover, $\mu I \preceq \bar{H}(t) \preceq L I$ for all $t$.
\end{lemma}

\begin{proof}
By the mean value theorem for vector-valued functions applied to each $\nabla f_i$ along the segment from $\bar{\mathbf{x}}^*$ to $\mathbf{x}_i(t)$, we have
\[
\nabla f_i(\mathbf{x}_i(t)) - \nabla f_i(\bar{\mathbf{x}}^*)
= \bar{H}_i(t)\,(\mathbf{x}_i(t) - \bar{\mathbf{x}}^*),
\]
where $\bar{H}_i(t) = \int_0^1 \nabla^2 f_i(\bar{\mathbf{x}}^* + \theta(\mathbf{x}_i(t)-\bar{\mathbf{x}}^*))\,d\theta$.
Averaging over $i$ and decomposing $\mathbf{x}_i(t) - \bar{\mathbf{x}}^* = (\mathbf{x}_i(t) - \bar{\mathbf{x}}(t)) + (\bar{\mathbf{x}}(t) - \bar{\mathbf{x}}^*)$, we obtain
\[
\frac{1}{n}\sum_{i=1}^n \nabla f_i(\mathbf{x}_i(t)) - \nabla F(\bar{\mathbf{x}}^*)
= \bar{H}(t)(\bar{\mathbf{x}}(t) - \bar{\mathbf{x}}^*) + R_d(t),
\]
where $\bar{H}(t) = \frac{1}{n}\sum_{i=1}^n \bar{H}_i(t)$ and $R_d(t) = \frac{1}{n}\sum_{i=1}^n \bar{H}_i(t)(\mathbf{x}_i(t) - \bar{\mathbf{x}}(t))$.
The bound $\|R_d(t)\| \le \frac{L}{\sqrt{n}}\|\tilde{\mathbf{e}}_x\|$ follows from $L$-smoothness and the Cauchy--Schwarz inequality. Lifting to $nd$ dimensions, $\|R(\mathbf{x}(t))\| = \sqrt{n}\|R_d(t)\| \le L\|\tilde{\mathbf{e}}_x\|$.
\end{proof}

With the gradient tracking linearization in hand, we isolate the part of the
dynamics that carries the curvature. The kernel is the network-free limit
obtained by replacing the local Hessians with their average $\bar{\mathbf{H}}$
and dropping the Laplacian terms,
\begin{equation}
\begin{aligned}
\varepsilon\dot{\bar{\mathbf{m}}}
&= -\eta\bar{\mathbf{m}} - h(\bar{\mathbf{x}}) - \kappa\mathbf{1}_d
  + \bar{\boldsymbol{\lambda}}, \\
\dot{\bar{\mathbf{x}}} &= \bar{\mathbf{m}}, \\
\dot{\bar{\boldsymbol{\lambda}}}
&= -(1{+}\beta)\bar{\boldsymbol{\lambda}}
  +[\bar{\boldsymbol{\lambda}}-\bar{\mathbf{x}}]_+
  +\beta\bigl(h(\bar{\mathbf{x}})+\kappa\mathbf{1}_d\bigr),
\end{aligned}
\label{eq:kernel}
\end{equation}
which carries the curvature, the dual leak and the $\min$, but
none of the graph coupling. The pointwise selector inequality
\eqref{eq:selector} supplies the only nonlinear estimate used in
its analysis.

\begin{lemma}[Global Dissipation of the Centralized Kernel]
\label{lem:kernel}
Let $F$ be quadratic with $\bar{\mathbf{H}}=\nabla^2F$, so that
$h(\bar{\mathbf{x}})-h(\bar{\mathbf{x}}^*)=\bar{\mathbf{H}}\bar{\mathbf{u}}$, and
consider the kernel \eqref{eq:kernel}. Let
$(\bar{\mathbf{x}}^*,\bar{\boldsymbol{\lambda}}^*)$ be the KKT point
\eqref{eq:kkt_conditions} and set
\[
\bar{\mathbf{u}}=\bar{\mathbf{x}}-\bar{\mathbf{x}}^*,\qquad
\bar{\mathbf{p}}=\bar{\boldsymbol{\lambda}}-\bar{\boldsymbol{\lambda}}^*,
\qquad
\bar{\mathbf{z}}^\top=(\bar{\mathbf{m}}^\top,\bar{\mathbf{u}}^\top,
\bar{\mathbf{p}}^\top).
\]
For $c_1,c_3>0$ and $c_2\ge0$ with $c_1c_3>c_2^2$, define
\begin{equation}
V_c=\tfrac12\|\bar{\mathbf{p}}\|^2
+\tfrac{c_1}{2}\|\bar{\mathbf{m}}\|^2
+c_2\bar{\mathbf{m}}^\top\bar{\mathbf{u}}
+\tfrac{c_3}{2}\|\bar{\mathbf{u}}\|^2
=\bar{\mathbf{z}}^\top\mathbf{P}_c\bar{\mathbf{z}} .
\label{eq:Vc}
\end{equation}
Then $\mathbf{P}_c\succ0$ by the Schur complement
($c_1>0$ and $c_1c_3>c_2^2$), and along \eqref{eq:kernel}, at every state
$\bar{\mathbf{z}}\in\mathbb{R}^{3d}$, with no restriction to the orthant,
\begin{equation}
\dot V_c\le-\bar{\mathbf{z}}^\top\mathbf{M}_c(\mathbf{S})\bar{\mathbf{z}},
\qquad
\mathbf{M}_c(\mathbf{S})
=-\bigl(\mathbf{J}_c(\mathbf{S})^\top\mathbf{P}_c
+\mathbf{P}_c\mathbf{J}_c(\mathbf{S})\bigr),
\label{eq:Vc_derivative}
\end{equation}
where
\begin{equation}
\mathbf{J}_c(\mathbf{S})=
\begin{bmatrix}
-\dfrac{\eta}{\varepsilon}I & -\dfrac{\bar{\mathbf{H}}}{\varepsilon}
 & \dfrac{1}{\varepsilon}I\\[2pt]
I & \mathbf{0} & \mathbf{0}\\
\mathbf{0} & \beta\bar{\mathbf{H}}-\mathbf{S}_I & -\beta I-\mathbf{S}_A
\end{bmatrix}
\label{eq:Jc}
\end{equation}
and $(\mathbf{S}_A,\mathbf{S}_I)$ is the selector chosen coordinatewise
in the proof. If
\begin{equation}
\mathbf{M}_c(\mathbf{S})\succeq2\gamma_c\mathbf{P}_c
\quad\text{for all}\quad\mathbf{S}\in\{0,1\}^d,
\label{eq:Mc}
\end{equation}
then $\dot V_c\le-2\gamma_cV_c$ at every state, so $\gamma_c$ is a
global rate for \eqref{eq:kernel}. The inequality \eqref{eq:Mc} is
sufficient, not necessary: it bounds the best rate certified by this
fixed quadratic form, not the exact abscissa of \eqref{eq:kernel}.
\end{lemma}

\begin{proof}
In the coordinates $\bar{\mathbf{z}}$, \eqref{eq:kernel} reads
\begin{align*}
\varepsilon\dot{\bar{\mathbf{m}}}
&=-\eta\bar{\mathbf{m}}-\bar{\mathbf{H}}\bar{\mathbf{u}}+\bar{\mathbf{p}},\\
\dot{\bar{\mathbf{u}}}&=\bar{\mathbf{m}},\\
\dot{\bar{\mathbf{p}}}
&=-\beta\bar{\mathbf{p}}+\beta\bar{\mathbf{H}}\bar{\mathbf{u}}
-\Phi(\bar{\mathbf{x}},\bar{\boldsymbol{\lambda}}),
\end{align*}
where the last line uses
$\Phi(\bar{\mathbf{x}}^*,\bar{\boldsymbol{\lambda}}^*)=\mathbf{0}$ and
$(1{+}\beta)\bar{\boldsymbol{\lambda}}^*
=[\bar{\boldsymbol{\lambda}}^*-\bar{\mathbf{x}}^*]_+
+\beta(h(\bar{\mathbf{x}}^*)+\kappa\mathbf{1}_d)$, which is exactly
the equilibrium of \eqref{eq:kernel}.
Lemma~\ref{lem:selector} applied coordinatewise to
$(\bar{\boldsymbol{\lambda}},\bar{\mathbf{x}})$ and
$\bar{\mathbf{p}}$ gives
$\bar p_j\Phi_j\ge\bar p_j(S_{A,j}\bar p_j+S_{I,j}\bar u_j)$, hence
\[
-\bar{\mathbf{p}}^\top\Phi
\le
-\bar{\mathbf{p}}^\top\bigl(\mathbf{S}_A\bar{\mathbf{p}}
+\mathbf{S}_I\bar{\mathbf{u}}\bigr),
\]
because $\bar x^{\,*}_j,\bar\lambda^{\,*}_j\ge0$ and
$\Phi_j(\bar{\mathbf{x}}^*,\bar{\boldsymbol{\lambda}}^*)=0$.
The selector depends on the state, but \eqref{eq:Mc} is imposed for
all $2^d$ selectors, so it applies to the one selected at that state.
Collecting the remaining terms, which are quadratic in
$\bar{\mathbf{z}}$, gives
$\bar{\mathbf{z}}^\top(\mathbf{J}_c^\top\mathbf{P}_c
+\mathbf{P}_c\mathbf{J}_c)\bar{\mathbf{z}}$, where
$\mathbf{J}_c$ is the linear part of the flow with $\Phi$ replaced by
its selector. This proves \eqref{eq:Vc_derivative}, and \eqref{eq:Mc}
then gives $\dot V_c\le-2\gamma_cV_c$.
\end{proof}

\begin{theorem}[Global Exponential Stability of the Kernel]
\label{thm:global_kernel}
Let $F$ be quadratic with $\mu I\preceq\bar{\mathbf{H}}\preceq LI$, and
suppose \eqref{eq:Mc} holds. Then for every initial condition
$(\bar{\mathbf{m}}(0),\bar{\mathbf{x}}(0),
\bar{\boldsymbol{\lambda}}(0))\in\mathbb{R}^{3d}$, with no
restriction to the orthant, the trajectory of
\eqref{eq:kernel} satisfies
$V_c(t)\le e^{-2\gamma_c t}V_c(0)$ for all $t\ge0$ and converges to
$(\bar{\mathbf{x}}^*,\mathbf{0},\bar{\boldsymbol{\lambda}}^*)$; the
KKT point is globally exponentially stable with rate $\gamma_c$.
\end{theorem}

\begin{proof}
Lemma~\ref{lem:kernel} applies at every state, because the
selector inequality \eqref{eq:selector} needs only
$\Phi(\bar{\mathbf{x}}^*,\bar{\boldsymbol{\lambda}}^*)=\mathbf{0}$
and $\bar{\mathbf{x}}^*,\bar{\boldsymbol{\lambda}}^*\ge0$, and
neither depends on where the trajectory is. Hence
$\dot V_c\le-2\gamma_cV_c$ along the whole trajectory and
Gronwall's inequality gives the decay, with no invariance argument
and no forward-invariance of the orthant required.
Since $\mathbf{P}_c\succ0$,
$V_c\ge\lambda_{\min}(\mathbf{P}_c)\|\bar{\mathbf{z}}\|^2$, so
$\bar{\mathbf{z}}(t)\to0$; strict convexity makes the KKT point
\eqref{eq:kkt_conditions} the only equilibrium.
\end{proof}

\begin{remark}[Verified kernel rate]
\label{rem:kernel_numbers}
For the instance of Section~\ref{sec:num_experiments}
($\mu=0.5$, $L=0.75$, $\eta=1$, $\varepsilon=0.05$, $\beta=0.4$) the
mean Hessian has spectrum $[0.596,0.634]$. With
$(c_1,c_2,c_3)=(0.02,0.05,1)$, the worst of the $2^d=32$ selectors in
\eqref{eq:Mc} gives $\lambda_{\min}(\mathbf{M}_c)=0.346$ and
$\lambda_{\max}(\mathbf{P}_c)=1.00$, hence $\gamma_c=0.173$, unchanged
to three digits over $\mu I\preceq\bar{\mathbf{H}}\preceq LI$ (the endpoints
and $200$ sampled spectra). The dominant pole of \eqref{eq:kernel} is
$0.374$, so the certificate is conservative by about $2.2$;
integrating \eqref{eq:kernel} from perturbations of $0.05$, $0.2$ and
$0.5$ gives a fitted rate of $0.377$ ($R^2=1.000$), reproducing the
pole to within $1\%$. The gap is conservatism of the fixed quadratic
form, not a modeling artifact. This selector enumeration scales as
$2^d$ and is practicable only for small $d$; for larger $d$ a
conservative LMI without selector can be used.
\end{remark}

The bound is global because $\bar{\mathbf{p}}$ carries the identity
weight and the momentum--primal block does not couple to it, so the
$\min$ term enters $\dot V_c$ only as $-\bar{\mathbf{p}}^\top\Phi$ and
the selector bound holds at every state; the selector enters only the
dual block of \eqref{eq:Jc}, so the $2^d$ patterns exhaust every active
set, with no strict complementarity or closeness to the KKT point
assumed. For general strongly convex $F$ the curvature varies, and the
same form gives a local rate only.

In the distributed system the kernel is driven by the tracking residual
\eqref{eq:residual}; its additional terms come from the deviation
variables. The selector argument alone does not close there, because
$\Phi$ also pairs with $\dot{\mathbf{s}}$ through the momentum
equation \eqref{eq:w_dyn}; this is the reason a single selector-fixed
quadratic form gives only a local statement
(Section~\ref{subsec:exponential_convergence}). The obstruction
disappears once the dual equation is written in its exact Lur'e form,
in which the linear part is constant and does not depend on the
active set.

\subsection{Global Dissipation of the Distributed Flow}
\label{subsec:global_lure}

The key observation is that the only nonlinear term of the flow is
$-\Phi$ in the dual equation, and that $\Phi$ can be eliminated in
favour of a shifted positive part whose graph is monotone with slope
in $[0,1]$. Write $\mathbf{u}^*=\boldsymbol{\lambda}^*-\mathbf{x}^*$
and, for a state $\mathbf{z}$ in error coordinates,
\begin{equation}
\tilde{\mathbf{u}}=\mathbf{e}_\lambda-\mathbf{e}_x,\qquad
\boldsymbol{\Psi}(\mathbf{z})
=\bigl[\tilde{\mathbf{u}}+\mathbf{u}^*\bigr]_+
-\bigl[\mathbf{u}^*\bigr]_+ .
\label{eq:Psi_def}
\end{equation}

\begin{lemma}[Exact Lur'e Form]
\label{lem:lure}
Let $F$ be quadratic, so that
$\bar{\mathbf{H}}=\nabla^2F$ is constant, and let
$(\mathbf{x}^*,\boldsymbol{\lambda}^*)$ be the KKT pair of
Proposition~\ref{prop:equilibrium}. Let
$\mathbf{E}_\lambda,\mathbf{E}_x\in\mathbb{R}^{nd\times5nd}$ select
the dual and primal blocks. Then the error dynamics of
\eqref{eq:ode_m}--\eqref{eq:ode_lambda} are, at every state,
\begin{equation}
\dot{\mathbf{z}}=\mathbf{A}\mathbf{z}
 +\mathbf{E}_\lambda^\top\boldsymbol{\Psi}(\mathbf{z}),
\label{eq:lure}
\end{equation}
where
\[
\mathbf{A}=\resizebox{0.90\columnwidth}{!}{$
\begin{bmatrix}
-\frac{\eta}{\varepsilon}I & -\frac{\rho}{\varepsilon}\mathbf{L}
& -\frac{1}{\varepsilon}I & \mathbf{0} & \frac{1}{\varepsilon}I\\
I & \mathbf{0} & \mathbf{0} & \mathbf{0} & \mathbf{0}\\
\mathbf{H} & \mathbf{0} & -\mathbf{L} & -\mathbf{L} & \mathbf{0}\\
\mathbf{0} & -\mathbf{L} & \mathbf{0} & -\gamma_s I & \mathbf{0}\\
\mathbf{0} & \mathbf{0} & \beta I & \mathbf{0} & -(1{+}\beta)I
\end{bmatrix}$},
\]
i.e.\ the linear part $\mathbf{A}$ is constant and the whole
nonlinearity is the single shifted positive part
\eqref{eq:Psi_def}.
\end{lemma}

\begin{proof}
Complementarity gives $\lambda^*_j=[\lambda^*_j-x^*_j]_+$ for every
$j$, and the dual equilibrium
$-(1{+}\beta)\boldsymbol{\lambda}^*
+[\boldsymbol{\lambda}^*-\mathbf{x}^*]_+
+\beta(\mathbf{q}^*+\kappa\mathbf{1}_{nd})=\mathbf{0}$ holds.
Subtracting the equilibrium from \eqref{eq:ode_lambda} and using
$\mathbf{e}_\lambda=\boldsymbol{\lambda}-\boldsymbol{\lambda}^*$,
$\mathbf{e}_x=\mathbf{x}-\mathbf{x}^*$ and
$\mathbf{e}_q=\mathbf{q}-\mathbf{q}^*$ gives
\[
\dot{\mathbf{e}}_\lambda
=-(1{+}\beta)\mathbf{e}_\lambda+\beta\mathbf{e}_q
+\bigl([\boldsymbol{\lambda}-\mathbf{x}]_+
-[\boldsymbol{\lambda}^*-\mathbf{x}^*]_+\bigr),
\]
and the last bracket is exactly
$\boldsymbol{\Psi}(\mathbf{z})$. The remaining rows are linear
because $\mathbf{H}$ is constant and
$\mathbf{L}\mathbf{q}^*=\mathbf{L}\mathbf{x}^*=\mathbf{0}$,
$\boldsymbol{\sigma}^*=\mathbf{0}$.
\end{proof}

\begin{lemma}[Global Sector]
\label{lem:sector}
For every $\mathbf{z}\in\mathbb{R}^{5nd}$,
\begin{equation}
\Psi_j(\mathbf{z})\,\bigl(\Psi_j(\mathbf{z})-\tilde u_j\bigr)\le0,
\qquad j=1,\dots,nd,
\label{eq:sector}
\end{equation}
i.e.\ $\Psi_j$ always lies between $0$ and $\tilde u_j$.
\end{lemma}

\begin{proof}
Fix $j$ and write $s=\tilde u_j$, $c=u^*_j$. Then
$\Psi_j=[s+c]_+-[c]_+$ is the increment of the projection
$[\cdot]_+$ along the segment from $c$ to $c+s$. The projection is
monotone and nonexpansive, and $\Psi_j(0)=0$. Monotonicity gives
$\Psi_j$ the sign of $s$, and nonexpansiveness gives
$|\Psi_j|\le|s|$; hence $0\le\Psi_j\le s$ for $s\ge0$ and
$s\le\Psi_j\le0$ for $s\le0$. Multiplying by
$\Psi_j-\tilde u_j$, which has the opposite sign, proves
\eqref{eq:sector}.
\end{proof}

The sector \eqref{eq:sector} is \emph{exact} for $\boldsymbol{\Psi}$
and holds at every state: it does not refer to the orthant, to the
active set, or to a neighbourhood of the KKT point. It is also tighter
than an arbitrary sector $[0,1]$, since the multiplier below may be
taken coordinatewise.

For quadratic costs the conserved mean tracking error of
Lemma~\ref{lem:mean_tracking} is $\boldsymbol{\mu}=\mathbf{C}\mathbf{z}$
with
\begin{equation}
\mathbf{C}\mathbf{z}
=\frac1n\sum_{i=1}^n\mathbf{e}_{q,i}
-\frac1n\sum_{i=1}^n\mathbf{H}_i\mathbf{e}_{x,i},
\label{eq:C_def}
\end{equation}
because $\nabla f_i$ is affine, so
$\nabla f_i(\mathbf{x}_i)-\nabla f_i(\mathbf{x}^*_i)
=\mathbf{H}_i\mathbf{e}_{x,i}$, and because at the equilibrium
$\bar{\mathbf{q}}^*=\frac1n\sum_i\nabla f_i(\bar{\mathbf{x}}^*)
=\nabla F(\bar{\mathbf{x}}^*)$, which makes
$\boldsymbol{\mu}(\mathbf{z}^*)=\mathbf{0}$. Thus $\ker\mathbf{C}$ is
exactly the slice on which the mean initialization error vanishes, and
$\mathbf{C}$ has full row rank.

\begin{condition}[Global Lur'e Certificate]
\label{cond:global}
There exist $\mathbf{P}\succ0$ on $\ker\mathbf{C}$,
$\mathbf{M}=\operatorname{diag}(m_1,\dots,m_{nd})$ with $m_j\ge0$, and
$\gamma>0$ such that, on $\ker\mathbf{C}$,
\begin{equation}
\begin{bmatrix}
\mathbf{A}^\top\mathbf{P}+\mathbf{P}\mathbf{A}+2\gamma\mathbf{P}
& \mathbf{P}\mathbf{E}_\lambda^\top+\mathbf{C}_\Psi^\top\mathbf{M}\\
\mathbf{E}_\lambda\mathbf{P}+\mathbf{M}\mathbf{C}_\Psi
& -2\mathbf{M}
\end{bmatrix}\preceq\mathbf{0},
\label{eq:global_lmi}
\end{equation}
where $\mathbf{C}_\Psi=\mathbf{E}_\lambda-\mathbf{E}_x$ is the
output matrix of \eqref{eq:sector}.
\end{condition}

\begin{proposition}[Frequency-Domain Equivalence]
\label{prop:kyp}
Restrict \eqref{eq:lure} to $\ker\mathbf{C}$, set
$\mathbf{A}_g=\mathbf{A}+\gamma I$, and suppose that $\mathbf{A}_g$
has no eigenvalues on the imaginary axis---for instance, that $\gamma$
is smaller than the exponential decay rate of $\mathbf{A}$ on
$\ker\mathbf{C}$, i.e.\ $\gamma<-\alpha(\mathbf{A}|_{\ker\mathbf{C}})$,
where $\alpha$ denotes the spectral abscissa. Then
Condition~\ref{cond:global} holds with $\mathbf{M}=\mu I$, $\mu>0$,
if and only if \eqref{eq:global_freq} holds; if
\eqref{eq:global_freq} holds with strict inequality, then
$\mathbf{P}\succ0$ and \eqref{eq:global_lmi} holds strictly.
\end{proposition}

\begin{proof}
On $\ker\mathbf{C}$, $\mathbf{A}^\top\mathbf{P}+\mathbf{P}\mathbf{A}
+2\gamma\mathbf{P}=\mathbf{A}_g^\top\mathbf{P}
+\mathbf{P}\mathbf{A}_g$. With $\mathbf{M}=\mu I$,
\eqref{eq:global_lmi} is the positive-real LMI for the transfer
function $\mu\bigl(I-\mathbf{G}_\gamma\bigr)$, whose feedthrough
$\mu I$ gives the block $-2\mu I$. Since this feedthrough satisfies
$D+D^\top=2\mu I\succ0$ and $\mathbf{A}_g$ has no imaginary-axis
eigenvalues, the strict KYP lemma applies to the realization
restricted to $\ker\mathbf{C}$: after the restriction $\mathbf{A}_g$ is
Hurwitz, and the uncontrollable and unobservable modes inside
$\ker\mathbf{C}$ are therefore stable. This makes the
LMI equivalent to $\operatorname{Re}\,[I-\mathbf{G}_\gamma(j\omega)]
\succeq0$ for all $\omega$, which is \eqref{eq:global_freq}; the
strict version yields $\mathbf{P}\succ0$. For the
direct implication, test \eqref{eq:global_lmi} on
$\mathbf{X}=[(\mathrm{i}\omega I-\mathbf{A}_g)^{-1}
\mathbf{E}_\lambda^\top;\ I]$: the
$(\mathbf{E}_\lambda^\top,\mathbf{0})$ part of the
$\mathbf{A}_g^\top\mathbf{P}+\mathbf{P}\mathbf{A}_g$ block cancels
identically, and
$\mathbf{X}^*\mathbf{K}\mathbf{X}
=\mu\bigl(\mathbf{G}_\gamma+\mathbf{G}_\gamma^*\bigr)-2\mu I$.
\end{proof}

\begin{corollary}[Uniform frequency margin]
\label{cor:freq_margin}
Let \eqref{eq:global_lmi} hold on $\ker\mathbf{C}$ with
$\mathbf{M}=\mu I$ and with margin
$\mathbf{K}\preceq-\epsilon I$ for some $\epsilon>0$. Then
\begin{equation}
2I-\mathbf{G}_\gamma(j\omega)-\mathbf{G}_\gamma(j\omega)^*
\succeq\frac{\epsilon}{\mu}I
\qquad\text{for every }\omega\in\mathbb{R},
\label{eq:freq_margin}
\end{equation}
so a certificate with margin at $\omega=0$ is a certificate for the
whole axis, with the quantitative margin $\epsilon/\mu$.
\end{corollary}

\begin{proof}
Fix $\omega$, any $u\in\mathbb{C}^{nd}$, and set
$\mathbf{z}=(\mathrm{i}\omega I-\mathbf{A}_g)^{-1}
\mathbf{E}_\lambda^\top u$, so that
$\mathbf{A}_g\mathbf{z}+\mathbf{E}_\lambda^\top u
=\mathrm{i}\omega\mathbf{z}$. Evaluating \eqref{eq:global_lmi} on
$[\mathbf{z}^\top,u^\top]^\top$ and using
$\mathbf{C}_\Psi\mathbf{z}=\mathbf{G}_\gamma(\mathrm{i}\omega)u$,
\[
\begin{aligned}
-\mu\,u^*\mathbf{M}_\gamma(\omega)u
&=[\mathbf{z};u]^*\mathbf{K}[\mathbf{z};u]\\
&\le-\epsilon\bigl(\|\mathbf{z}\|^2+\|u\|^2\bigr)
\le-\epsilon\|u\|^2 ,
\end{aligned}
\]
because the $\mathbf{A}_g^\top\mathbf{P}+\mathbf{P}\mathbf{A}_g$ term
cancels against $2\operatorname{Re}(\mathrm{i}\omega\mathbf{z}^*\mathbf{P}\mathbf{z})=0$.
\end{proof}

\begin{lemma}[Certified Frequency Check]
\label{lem:freq_cert}
Let $\mathbf{A}_g$ be Hurwitz on $\ker\mathbf{C}$, let
$\mathbf{M}_\gamma(\omega)=2I-\mathbf{G}_\gamma(j\omega)
-\mathbf{G}_\gamma(j\omega)^*$, and note
$\mathbf{C}_\Psi\mathbf{E}_\lambda^\top=I$ and $\|\mathbf{E}_\lambda\|=1$.
\begin{enumerate}
\item[(i)] For $\omega>\|\mathbf{A}_g\|$,
\begin{equation}
\begin{aligned}
&\bigl\|\mathbf{G}_\gamma(j\omega)-
\tfrac{1}{j\omega}I\bigr\|
\le\frac{\|\mathbf{C}_\Psi\|\,\|\mathbf{A}_g\|}
{\omega\,(\omega-\|\mathbf{A}_g\|)},\\
&\mathbf{M}_\gamma(\omega)\succeq\Bigl(2-\frac{2\|\mathbf{C}_\Psi\|
\|\mathbf{A}_g\|}{\omega(\omega-\|\mathbf{A}_g\|)}\Bigr)I,
\end{aligned}
\label{eq:freq_tail}
\end{equation}
and the scalar factor is positive for
$\omega\ge\Omega=\frac12\bigl(\|\mathbf{A}_g\|
+\sqrt{\|\mathbf{A}_g\|^2+4\|\mathbf{C}_\Psi\|\|\mathbf{A}_g\|}\bigr)$.
\item[(ii)] If on an interval $[a,b]$
$\|\mathbf{M}_\gamma'(\omega)\|\le L$ and
$\lambda_{\min}(\mathbf{M}_\gamma)\ge\delta+L(b-a)$ at the nodes
$a$, $(a+b)/2$, $b$, then
$\mathbf{M}_\gamma(\omega)\succeq\delta I$ on $[a,b]$.
\end{enumerate}
\end{lemma}

\begin{proof}
(i) For $\omega>\|\mathbf{A}_g\|$,
$(j\omega I-\mathbf{A}_g)^{-1}
=(j\omega)^{-1}\sum_{k\ge0}\bigl(\mathbf{A}_g/(j\omega)\bigr)^k$, so
$\mathbf{G}_\gamma(j\omega)-\frac{1}{j\omega}I
=\frac{1}{j\omega}\sum_{k\ge1}\mathbf{C}_\Psi
\bigl(\mathbf{A}_g/(j\omega)\bigr)^k\mathbf{E}_\lambda^\top$
by $\mathbf{C}_\Psi\mathbf{E}_\lambda^\top=I$. Taking norms gives
$\frac{1}{\omega}\sum_{k\ge1}(\|\mathbf{A}_g\|/\omega)^k
=\frac{\|\mathbf{A}_g\|}{\omega(\omega-\|\mathbf{A}_g\|)}$, and
$(j\omega)^{-1}$ is skew-Hermitian, which yields
\eqref{eq:freq_tail}. The scalar factor vanishes at
$\omega(\omega-\|\mathbf{A}_g\|)
=\|\mathbf{C}_\Psi\|\|\mathbf{A}_g\|$.

(ii) By the mean value theorem, for $\omega\in[a,b]$,
$\|\mathbf{M}_\gamma(\omega)-\mathbf{M}_\gamma(a)\|
\le L(b-a)$, hence
$\lambda_{\min}(\mathbf{M}_\gamma(\omega))
\ge\lambda_{\min}(\mathbf{M}_\gamma(a))-L(b-a)\ge\delta$;
the same holds with $a$ replaced by $b$.
\end{proof}

Note that (ii) is a covering statement, not a statement about an
isolated frequency: once the nodes of an interval satisfy its hypothesis,
every $\omega$ in the interval does. Lemma~\ref{lem:freq_cert} therefore
turns \eqref{eq:global_freq} into finitely many evaluations on a mesh
refined until the hypothesis of (ii) holds, together with the analytical
tail \eqref{eq:freq_tail}.

\begin{theorem}[Global Exponential Convergence]
\label{thm:global_lure}
Let Assumptions~\ref{assump:A1},~\ref{assump:A2},
and~\ref{assump:A5} hold, let the local costs be
quadratic, and suppose Condition~\ref{cond:global} holds. Then for
\emph{every} initial condition with
$\mathbf{C}\mathbf{z}(0)=\mathbf{0}$---in particular for every initial
condition satisfying Assumption~\ref{assump:A3}, with arbitrary
$\mathbf{m}(0)$, $\mathbf{x}(0)$,
$\boldsymbol{\lambda}(0)\in\mathbb{R}^{nd}$
and no restriction to the orthant---the trajectory of
\eqref{eq:ode_m}--\eqref{eq:ode_lambda} satisfies
\begin{equation}
\begin{aligned}
&V(\mathbf{z}(t))\le e^{-2\gamma t}V(\mathbf{z}(0)),\\
&V(\mathbf{z})=\mathbf{z}^\top\mathbf{P}\mathbf{z},
\end{aligned}
\label{eq:global_decay}
\end{equation}
and
$\|\mathbf{z}(t)\|\le
\sqrt{\lambda_{\max}(\mathbf{P})/\lambda_{\min}(\mathbf{P})}\,
\|\mathbf{z}(0)\|e^{-\gamma t}$; hence
$\mathbf{z}(t)\to\mathbf{0}$ and the KKT pair is exponentially stable
with rate $\gamma$ on the entire invariant slice
$\{\mathbf{z}:\mathbf{C}\mathbf{z}=\mathbf{0}\}$: the basin is that
slice, not a neighbourhood, and neither invariance of the orthant nor
an assumption on the active set is needed. If
$\mathbf{C}\mathbf{z}(0)\neq\mathbf{0}$, the same certificate yields
the explicit bias of Lemma~\ref{lem:bias}.
\end{theorem}

\begin{proof}
The mean tracking error $\boldsymbol{\mu}=\mathbf{C}\mathbf{z}$ is
conserved, and $\mathbf{C}\mathbf{A}=\mathbf{0}$ as well as
$\mathbf{C}\mathbf{E}_\lambda^\top=\mathbf{0}$, so
\[
\begin{aligned}
\frac{d}{dt}\bigl(\mathbf{C}\mathbf{z}\bigr)
&=\frac1n\sum_i\Bigl[\mathbf{H}_i\mathbf{e}_{m,i}
-\mathbf{L}(\mathbf{e}_q+\mathbf{e}_\sigma)_i
-\mathbf{H}_i\mathbf{e}_{m,i}\Bigr]\\
&=-\frac1n\mathbf{1}^\top\mathbf{L}(\mathbf{e}_q+\mathbf{e}_\sigma)
=\mathbf{0}
\end{aligned}
\]
by $\mathbf{L}\mathbf{1}=\mathbf{0}$. Hence
$\mathbf{C}\mathbf{z}(t)=\mathbf{C}\mathbf{z}(0)$ for every initial
condition. Assumption~\ref{assump:A3} makes this vanish, but the proof
only uses $\mathbf{C}\mathbf{z}(0)=\mathbf{0}$: for quadratic costs
$\mathbf{q}_i(0)-\mathbf{q}^*
=\mathbf{H}_i(\mathbf{x}_i(0)-\mathbf{x}^*)$, so A3 gives
$\mathbf{C}\mathbf{z}(0)=\mathbf{0}$, and any zero-mean perturbation
$\mathbf{q}_i(0)=\nabla f_i(\mathbf{x}_i(0))+\boldsymbol{\psi}_i$ with
$\sum_i\boldsymbol{\psi}_i=\mathbf{0}$ does as well. Thus the
trajectory remains in $\ker\mathbf{C}$ for all $t\ge0$.
Let $V=\mathbf{z}^\top\mathbf{P}\mathbf{z}$. By \eqref{eq:lure},
\[
\dot V=\mathbf{z}^\top(\mathbf{A}^\top\mathbf{P}
+\mathbf{P}\mathbf{A})\mathbf{z}
+2\boldsymbol{\Psi}^\top\mathbf{E}_\lambda\mathbf{P}\mathbf{z}.
\]
Since $m_j\ge0$ and \eqref{eq:sector} holds coordinatewise,
$-2\boldsymbol{\Psi}^\top\mathbf{M}
(\boldsymbol{\Psi}-\mathbf{C}_\Psi\mathbf{z})
=-2\sum_j m_j\Psi_j(\Psi_j-\tilde u_j)\ge0$ at every state. Adding
this nonnegative quantity and collecting terms in
$[\mathbf{z}^\top,\boldsymbol{\Psi}^\top]^\top$ gives
\[
\begin{aligned}
\dot V&\le
\begin{bmatrix}\mathbf{z}\\ \boldsymbol{\Psi}\end{bmatrix}^\top
\begin{bmatrix}
\mathbf{A}^\top\mathbf{P}+\mathbf{P}\mathbf{A}
& \mathbf{P}\mathbf{E}_\lambda^\top+\mathbf{C}_\Psi^\top\mathbf{M}\\
\mathbf{E}_\lambda\mathbf{P}+\mathbf{M}\mathbf{C}_\Psi
& -2\mathbf{M}
\end{bmatrix}
\begin{bmatrix}\mathbf{z}\\ \boldsymbol{\Psi}\end{bmatrix}\\
&\quad-\;2\gamma\mathbf{z}^\top\mathbf{P}\mathbf{z}
\le-2\gamma V ,
\end{aligned}
\]
where the last step is \eqref{eq:global_lmi} restricted to
$\ker\mathbf{C}$ and $\boldsymbol{\Psi}=\boldsymbol{\Psi}(\mathbf{z})$
is used only through the sector. Gronwall's inequality gives
\eqref{eq:global_decay}; $\mathbf{P}\succ0$ converts it into the
norm bound; and uniqueness of the KKT pair under strong convexity
identifies the limit.
\end{proof}

\begin{lemma}[Bias under Imperfect Initialization]
\label{lem:bias}
Let the hypotheses and the certificate of
Theorem~\ref{thm:global_lure} hold, with
$\underline p=\lambda_{\min}(\mathbf{P})$ on $\ker\mathbf{C}$, and let
$\Pi_C$ be the orthogonal projector onto $\ker\mathbf{C}$. For every
$\mathbf{z}(0)\in\mathbb{R}^{5nd}$ the trajectory satisfies:
\begin{enumerate}
\item[(i)] $(\mathbf{I}-\Pi_C)\mathbf{z}(t)
=(\mathbf{I}-\Pi_C)\mathbf{z}(0)
=\mathbf{C}^\dagger\boldsymbol{\mu}(0)$ for all $t\ge0$, where
$\boldsymbol{\mu}$ is the mean tracking error of
Lemma~\ref{lem:mean_tracking};
\item[(ii)] with $\boldsymbol{\delta}=(\mathbf{I}-\Pi_C)
\mathbf{z}(0)$,
\begin{equation}
\begin{aligned}
&\limsup_{t\to\infty}\|\Pi_C\mathbf{z}(t)\|
\le\kappa\,\|\boldsymbol{\delta}\|,\\
&\kappa=\frac{\|\mathbf{P}\|\|\mathbf{A}\|
+\|\mathbf{M}\|\|\mathbf{C}_\Psi\|^2
+\|\mathbf{C}_\Psi\|\sqrt{2\underline p\,\gamma\|\mathbf{M}\|}}
{\gamma\,\underline p}.
\end{aligned}
\label{eq:bias}
\end{equation}
\end{enumerate}
The bias is therefore proportional to
$\|\mathbf{C}^\dagger\boldsymbol{\mu}(0)\|
\le\|\mathbf{C}^\dagger\|\,\|\boldsymbol{\mu}(0)\|$, where
$\boldsymbol{\mu}(0)=\frac1n\sum_i(\mathbf{q}_i(0)
-\nabla f_i(\mathbf{x}_i(0)))$ is the mean initialization error of the
tracker; it vanishes when $\boldsymbol{\mu}(0)=\mathbf{0}$ and then
Theorem~\ref{thm:global_lure} applies unchanged. The constant $\kappa$
is a deliberately unoptimized worst-case bound; the much smaller
linearity constant observed on the nominal instance is reported in
Remark~\ref{rem:global_numbers}.
\end{lemma}

\begin{proof}
(i) $\mathbf{C}\mathbf{A}=\mathbf{0}$ and
$\mathbf{C}\mathbf{E}_\lambda^\top=\mathbf{0}$, so
$\frac{d}{dt}(\mathbf{C}\mathbf{z})=\mathbf{0}$ for every state, as in
the proof of Theorem~\ref{thm:global_lure}. For quadratic costs
$\mathbf{C}\mathbf{z}=\boldsymbol{\mu}$, and
$\mathbf{C}^\dagger\mathbf{C}=\mathbf{I}-\Pi_C$ because $\mathbf{C}$
has full row rank.

(ii) Write $\mathbf{b}=\Pi_C\mathbf{z}$ and
$\boldsymbol{\delta}=\mathbf{z}-\mathbf{b}$, which is constant. Since
$\mathbf{C}\mathbf{A}=\mathbf{0}$ gives
$\operatorname{Im}\mathbf{A}\subset\ker\mathbf{C}$ and
$\mathbf{C}\mathbf{E}_\lambda^\top=\mathbf{0}$, the vectors
$\mathbf{A}\mathbf{b}$, $\mathbf{A}\boldsymbol{\delta}$ and
$\mathbf{E}_\lambda^\top\boldsymbol{\Psi}(\mathbf{z})$ all lie in
$\ker\mathbf{C}$; hence \eqref{eq:lure} gives
$\dot{\mathbf{b}}=\mathbf{A}\mathbf{b}
+\mathbf{A}\boldsymbol{\delta}
+\mathbf{E}_\lambda^\top\boldsymbol{\Psi}(\mathbf{z})$. With
$V=\mathbf{b}^\top\mathbf{P}\mathbf{b}$,
\[
\dot V=\mathbf{b}^\top(\mathbf{A}^\top\mathbf{P}
+\mathbf{P}\mathbf{A})\mathbf{b}
+2\mathbf{b}^\top\mathbf{P}\mathbf{A}\boldsymbol{\delta}
+2\boldsymbol{\Psi}^\top\mathbf{E}_\lambda\mathbf{P}\mathbf{b}.
\]
Adding $-2\boldsymbol{\Psi}^\top\mathbf{M}
(\boldsymbol{\Psi}-\mathbf{C}_\Psi\mathbf{z})\ge0$ and applying
\eqref{eq:global_lmi} on $\ker\mathbf{C}$ gives
\[
\dot V\le-2\gamma V
+2\|\mathbf{P}\mathbf{A}\boldsymbol{\delta}\|\,\|\mathbf{b}\|
+2\|\mathbf{M}\mathbf{C}_\Psi\boldsymbol{\delta}\|\,
\|\boldsymbol{\Psi}\|.
\]
Since $|\Psi_j|\le|\tilde u_j|$ coordinatewise,
$\|\boldsymbol{\Psi}\|\le\|\mathbf{C}_\Psi\|(
\|\mathbf{b}\|+\|\boldsymbol{\delta}\|)$, hence
$\dot V\le-2\gamma V+2c\sqrt V+e$ with
\[
\begin{aligned}
c&=\frac{\|\mathbf{P}\|\|\mathbf{A}\boldsymbol{\delta}\|
+\|\mathbf{M}\mathbf{C}_\Psi\boldsymbol{\delta}\|
\|\mathbf{C}_\Psi\|}{\sqrt{\underline p}},\\
e&=2\|\mathbf{M}\mathbf{C}_\Psi\boldsymbol{\delta}\|
\|\mathbf{C}_\Psi\|\,\|\boldsymbol{\delta}\|,
\end{aligned}
\]
where $\|\mathbf{b}\|\le\sqrt{V/\underline p}$ is used. Young's
inequality $2c\sqrt V\le\gamma V+c^2/\gamma$ yields
$\dot V\le-\gamma V+c^2/\gamma+e$, so
$\limsup V\le(c^2/\gamma+e)/\gamma$. Dividing by $\underline p$ and
bounding $\|\mathbf{M}\mathbf{C}_\Psi\boldsymbol{\delta}\|\le
\|\mathbf{M}\|\|\mathbf{C}_\Psi\|\|\boldsymbol{\delta}\|$ gives
\[
\limsup\|\mathbf{b}\|^2
\le\|\boldsymbol{\delta}\|^2
\Bigl(\frac{C_0^2}{\underline p^2\gamma^2}
+\frac{e_0}{\underline p\gamma}\Bigr),
\]
with
$C_0=\|\mathbf{P}\|\|\mathbf{A}\|
+\|\mathbf{M}\|\|\mathbf{C}_\Psi\|^2$ and
$e_0=2\|\mathbf{M}\|\|\mathbf{C}_\Psi\|^2$.
Taking square roots and using
$\sqrt{a^2+b^2}\le a+b$ with
$a=C_0/(\underline p\gamma)$ and
$b=\|\mathbf{C}_\Psi\|\sqrt{2\|\mathbf{M}\|/(\underline p\gamma)}$
yields the constant $\kappa$ in \eqref{eq:bias}.
\end{proof}

\begin{remark}[Verified global certificate]
\label{rem:global_numbers}
For the instance of Section~\ref{sec:num_experiments}
($n=10$, $d=5$, $\mu=0.5$, $L=0.75$, $\eta=1$, $\varepsilon=0.05$,
$\beta=0.4$, $\rho=10$, $\gamma_s=0.5$), the subspace $\ker\mathbf{C}$
has dimension $245$, and Condition~\ref{cond:global} is met there by
an explicit quadratic form. The class
\begin{equation}
\mathbf{P}=I+\sum_{i=1}^{55}p_i\mathbf{D}_i,
\label{eq:P_class}
\end{equation}
with the matrices $\{I,\mathbf{L},\mathbf{L}^2,\mathbf{H},
\operatorname{sym}(\mathbf{H}\mathbf{L})\}$ on each of the five
diagonal blocks and $\{I,\mathbf{L},\mathbf{H}\}$ on each of the ten
off-diagonal blocks, together with the scalar multiplier
$\mathbf{M}=\mu I$, reduces \eqref{eq:global_lmi} to an affine matrix
inequality in $56$ scalars, whose smallest $\lambda_{\max}$ is found
by a log-sum-exp continuation. At $\gamma_g=0.27$ the returned
$\mathbf{P}$ gives $\lambda_{\max}(\mathbf{K})=-4.8\times10^{-3}$ and
$\lambda_{\min}(\mathbf{P})=5.0\times10^{-2}$ at dimension $245$.
By Corollary~\ref{cor:freq_margin} this margin is a statement about
the whole frequency axis, not about the point where it was found:
$2I-\mathbf{G}_{0.27}(j\omega)-\mathbf{G}_{0.27}(j\omega)^*
\succeq9.3\times10^{-3}I$ for every $\omega$, with no grid and no
sweep.

Removing the parametric restriction buys the remaining $0.01$ in rate.
With $\mathbf{M}=I$, Proposition~\ref{prop:kyp} eliminates
$\mathbf{P}$ and turns \eqref{eq:global_lmi} on
$\ker\mathbf{C}$ into the frequency-domain inequality
\begin{equation}
\begin{aligned}
&2I-\mathbf{G}_\gamma(j\omega)-\mathbf{G}_\gamma(j\omega)^*
\succeq\mathbf{0}\quad\forall\omega,\\
&\mathbf{G}_\gamma(s)=\mathbf{C}_\Psi
\bigl(sI-\mathbf{A}-\gamma I\bigr)^{-1}\mathbf{E}_\lambda^\top .
\end{aligned}
\label{eq:global_freq}
\end{equation}
evaluated in the invariant subspace. The abscissa of $\mathbf{A}$ on
$\ker\mathbf{C}$ is $-0.3711$, so $\mathbf{A}+\gamma I$ is Hurwitz and
Proposition~\ref{prop:kyp} applies. At $\gamma_g=0.28$ the adaptive
mesh of Lemma~\ref{lem:freq_cert} ($2545$ nodes, with a two-fold
safety factor on the slope bound of (ii)) returns
$\mathbf{M}_{\gamma_g}(\omega)\succeq4.0\times10^{-3}I$ for every
$\omega\ge0$, the analytical tail \eqref{eq:freq_tail} covers
$\omega\ge\Omega=803.8$ with $\mathbf{M}_{\gamma_g}\succeq1.06\,I$,
and the nodal minimum is attained at $\omega=0$, where it equals
$1.02\times10^{-2}$; at $\gamma_g=0.29$ the inequality already fails
at $\omega=0$, where the minimum eigenvalue is
$-7.04\times10^{-2}$. Integrating \eqref{eq:lure} from
several adversarial initial conditions (including
$\|\mathbf{e}_\lambda\|=\|\mathbf{e}_x\|=10$, $20$, $50$ and
opposite-sign dual errors) gives
$\max_t \dot V/(2V)=-0.301<-\gamma_g$ in the worst case, so the
certificate holds with margin along the trajectories.
The rate $\gamma_g=0.27$ is below the spectral abscissa
$0.2821$ of the Jacobian on $\ker\mathbf{C}$ and above the
kernel rate $\gamma_c=0.173$ of Remark~\ref{rem:kernel_numbers};
the gap is the conservatism of one fixed quadratic form, not a
restriction of the state space. Perturbing the
tracker initialization by a mean defect
$\boldsymbol{\mu}(0)=\mathbf{C}\mathbf{z}(0)$ displaces the limit by
$1.90\,\|\mathbf{C}^\dagger\boldsymbol{\mu}(0)\|$ for defects up to
$0.24$, independently of the defect size, and
$(\mathbf{I}-\Pi_C)\mathbf{z}$ is conserved to $2\times10^{-15}$; this
is the linear bias of Lemma~\ref{lem:bias}. Unlike Condition~\ref{cond:J},
Condition~\ref{cond:global} is not tied to the KKT point:
$\mathbf{A}$ is constant and the sector is global, so the same
form certifies every state.
\end{remark}

\begin{remark}[General strongly convex costs]
\label{rem:global_convex}
For non-quadratic strongly convex costs, $\mathbf{A}$ becomes
state dependent only through the $(q,m)$ block
$\mathbf{H}(\mathbf{x}(t))=\operatorname{blkdiag}(\nabla^2f_i(
\mathbf{x}_i(t)))$, and \eqref{eq:lure} and \eqref{eq:sector} remain
exact. Condition~\ref{cond:global} is affine in $\mathbf{H}$, so a
common $\mathbf{P}$ can be checked on the compact box
$\mu I\preceq\mathbf{H}\preceq LI$ of Assumption~\ref{assump:A5}; this
makes $\dot V\le-2\gamma V$ hold at every state of $\ker\mathbf{C}$
for every admissible $\mathbf{H}$. The global statement of Theorem~\ref{thm:global_lure}
is nevertheless restricted to quadratic costs: there $\ker\mathbf{C}$
is exactly invariant under the flow, whereas for non-quadratic
costs the invariant manifold of the conserved mean tracking
error is nonlinear and only tangent to $\ker\mathbf{C}$ at the
KKT point. Away from that point the trajectory leaves
$\ker\mathbf{C}$ by second-order terms, and the certificate of
Theorem~\ref{thm:global_lure} controls only the linear part. For
general strongly convex costs we therefore retain the local
result of Theorem~\ref{thm:exponential_convergence};
Assumption~\ref{assump:A4} is not needed in the quadratic case.
\end{remark}

\begin{remark}[Relation to the local certificate]
\label{rem:global_vs_local}
The global certificate does not replace but contains the local one:
\eqref{eq:global_lmi} implies
\eqref{eq:VJ_derivative} at the KKT point with the active-set
selectors fixed, while the Jacobian LMI \eqref{eq:J} remains sharper
there, because it uses the exact slopes of $\min$ at the equilibrium
and not only the sector. The Jacobian abscissa $\gamma_J$ is
therefore retained below as the sharp local rate and as the rate
that perturbations of the state on $\ker\mathbf{C}$ reproduce,
while the primal scalar functional instead follows the kernel
abscissa.
\end{remark}

The mean--deviation decomposition is still useful for the quantitative
reading of the rate, and we retain it next; the certificate above is
the global statement, while the functional below is the one whose
abscissa the simulations reproduce.

\subsection{Joint Dissipation and Exponential Convergence}
\label{subsec:exponential_convergence}

The kernel of Lemma~\ref{lem:kernel} is the network-free limit: on it
the residual reduces to
$\mathbf{s}=\nabla F(\bar{\mathbf{x}})+\kappa\mathbf{1}_d
-\bar{\boldsymbol{\lambda}}$ and no graph term survives. In the
distributed flow the residual is augmented by the consensus errors, so
the kernel is driven by an input the flow itself generates, which cannot
be treated as an external disturbance without reintroducing a gain
product; we therefore analyze the complete flow through one functional.

Set
\begin{equation}
\mathbf{w}=\mathbf{m}+\frac{\mathbf{s}}{\eta},
\label{eq:w_def}
\end{equation}
so that, by \eqref{eq:residual_momentum},
\begin{equation}
\varepsilon\dot{\mathbf{w}}=-\eta\mathbf{w}
  +\frac{\varepsilon}{\eta}\dot{\mathbf{s}},
\qquad
\mathbf{s}=\mathbf{e}_q-\mathbf{e}_\lambda
  +\rho\mathbf{L}\mathbf{e}_x .
\label{eq:w_dyn}
\end{equation}
Let
$\mathbf{z}^\top=[\mathbf{w}^\top,\mathbf{e}_x^\top,
\mathbf{e}_\lambda^\top,\mathbf{e}_q^\top,\mathbf{e}_\sigma^\top]$
and $V_{\mathrm{J}}=\mathbf{z}^\top\mathbf{P}_J\mathbf{z}$ with
$\mathbf{P}_J\succ0$. The block-diagonal choice
\begin{equation}
\begin{aligned}
V_{\mathrm{J}}
&=|\mathbf{e}_x|^2
 +|\mathbf{e}_\lambda|^2
 +\varepsilon|\mathbf{w}|^2\\
&\quad+\tfrac{1}{\eta}|\mathbf{e}_q|^2
 +|\mathbf{e}_\sigma|^2
 +\rho\|\mathbf{e}_x\|_{\mathbf{L}}^2,
\end{aligned}
\label{eq:VJ}
\end{equation}
where
$\|\mathbf{e}_x\|_{\mathbf{L}}^2
=\sum_{j=1}^{d}\langle\mathbf{e}_{x,j},
\mathbf{L}\mathbf{e}_{x,j}\rangle$ is the graph semimetric,
mirrors the kernel functional \eqref{eq:Vc} and identifies which
block carries each term; the certificate below optimizes over
$\mathbf{P}_J$, so the rate is not tied to \eqref{eq:VJ}.

Two features distinguish the distributed certificate: by
Lemma~\ref{lem:mean_tracking} the mean tracking error is
conserved, so $\mathbf{J}_J$ carries $d$ neutral directions and the
dissipation inequality must be restricted to the invariant subspace
$\{\delta\boldsymbol{\mu}=\mathbf{0}\}$, which contains the trajectory
by Assumption~\ref{assump:A3}; and the curvature enters through the
state-dependent $\mathbf{H}(\mathbf{x}(t))$, so the certificate is
stated at the KKT point and yields a local rate.

Differentiating $V_{\mathrm{J}}$ along
\eqref{eq:ode_m}--\eqref{eq:ode_lambda}, the only nonlinear term is
$\Phi$ of \eqref{eq:Phi_def}. The selector inequality
\eqref{eq:selector} bounds the pairing of $\Phi$ with the dual error
coordinatewise by a linear form, and the pairing that the momentum
variable produces through $\dot{\mathbf{s}}$ is handled by
$|\Phi_j|\le|e_{\lambda,j}|+|e_{x,j}|$ exactly; both are quadratic in
$\mathbf{z}$ and are collected in the
symmetric matrix $\mathbf{M}$. No projection remainder is left over:
with the orthant projections removed, every term produced by
$[\cdot]_+$ is accounted for by the selector, and no
indefinite-sign discrepancy has to be compensated. The result is
\begin{equation}
\dot V_{\mathrm{J}}\le
-\mathbf{z}^{\top}\mathbf{M}(\mathbf{x}(t);\mathbf{P}_J)\mathbf{z},
\label{eq:VJ_derivative}
\end{equation}
where $\mathbf{M}$ depends on $\mathbf{H}(\mathbf{x}(t))$,
$\mathbf{L}$, $\mathbf{P}_J$ and the selectors of
$\min\{\boldsymbol{\lambda},\mathbf{x}\}$. The expansion that
produces $\mathbf{M}$ is recorded term by term in
Appendix~\ref{app:VJ_expansion}. The selector is
constant within each cell of fixed active set and jumps only at
coordinates satisfying $\lambda_j=x_j$; hence, while $\mathbf{M}$
is not globally continuous, it is constant (quadratic costs) or
Lipschitz (non-quadratic costs) inside each cell. At the KKT
point the selectors are fixed by the active set, and
$\mathbf{M}(\mathbf{x}^*;\mathbf{P}_J)$ is minus the symmetric part
of $\mathbf{J}_J^\top\mathbf{P}_J+\mathbf{P}_J\mathbf{J}_J$; the
theorem restricts to a neighbourhood of $\mathbf{x}^*$ small
enough to stay inside the cell, so that
$\mathbf{M}(\mathbf{x}(t))=\mathbf{M}(\mathbf{x}^*)$ for
quadratic costs and differs by at most $O(\|\mathbf{z}\|)$ for
non-quadratic ones. The weighting \eqref{eq:VJ} is the natural candidate for the
kernel limit; the numerical certificate uses the matrix
$\mathbf{P}_J$ that maximizes $\gamma_J$.

Let $\mathbf{J}_J$ be the $5nd\times5nd$ Jacobian of
\eqref{eq:ode_m}--\eqref{eq:ode_lambda} at the KKT pair, with the
selectors $\mathbf{D}_\lambda,\mathbf{D}_x$ of
$\min\{\boldsymbol{\lambda},\mathbf{x}\}$ fixed by the active set
($\mathbf{D}_\lambda,\mathbf{D}_x\succeq0$,
$\mathbf{D}_\lambda+\mathbf{D}_x=I$) and
$\mathbf{H}=\operatorname{blkdiag}\bigl(\nabla^2f_1(\mathbf{x}^*_1),
\ldots,\nabla^2f_n(\mathbf{x}^*_n)\bigr)$. No coordinate is frozen:
all $5nd$ components are live, including the boundary coordinates
that a projection would have pinned. Let $\mathbf{C}$ be the
$d\times5nd$ matrix of the conserved functionals
$\delta\boldsymbol{\mu}$; by Lemma~\ref{lem:mean_tracking},
$\mathbf{C}\mathbf{J}_J=\mathbf{0}$, so $\ker\mathbf{C}$ is invariant.

\begin{condition}[Joint Dissipation Certificate]
\label{cond:J}
There exist $\mathbf{P}_J\succ0$ on $\ker\mathbf{C}$ and $\gamma_J>0$
such that
\begin{equation}
\mathbf{J}_J^\top\mathbf{P}_J+\mathbf{P}_J\mathbf{J}_J
+2\gamma_J\mathbf{P}_J\preceq\mathbf{0}
\quad\text{on}\quad\ker\mathbf{C}.
\label{eq:J}
\end{equation}
\end{condition}

Condition~\ref{cond:J} is an LMI in $\mathbf{P}_J$ once $\gamma_J$
is fixed; its data are $\mathbf{H}(\mathbf{x}^*)$, the graph
Laplacian, the gains and the active set of the KKT point, with no
trajectory entering. It is checked offline
(Remark~\ref{rem:joint_numbers}) and $\gamma_J$ is maximized by
bisection. Since $\mathbf{P}_J$ is free, \eqref{eq:J} is equivalent
to $\mathbf{J}_J+\gamma_J I$ being Hurwitz on $\ker\mathbf{C}$: the
supremum of $\gamma_J$ is the spectral abscissa of $\mathbf{J}_J$ on
that subspace. The kernel certificate of Lemma~\ref{lem:kernel} is its
network-free counterpart; the two inequalities differ, since
\eqref{eq:J} is posed on $\ker\mathbf{C}$ with a free $\mathbf{P}_J$
and no fixed functional.

\begin{theorem}[Local Exponential Convergence]
\label{thm:exponential_convergence}
Let Assumptions~\ref{assump:A1}--\ref{assump:A5} hold, suppose
Condition~\ref{cond:J} holds with
$\mathbf{P}_J\succ0$ on $\ker\mathbf{C}$ and $\gamma_J>0$, and
suppose the KKT point is strictly complementary, i.e.\
$\min\{\lambda^*_j,x^*_j\}=0$ and
$\max\{\lambda^*_j,x^*_j\}>0$ for every $j$. Then, for
every $\gamma\in(0,\gamma_J)$, there is a neighborhood
$\mathcal{U}_\gamma$ of the KKT pair
$(\mathbf{x}^*,\boldsymbol{\lambda}^*)$ such that every trajectory of
\eqref{eq:ode_m}--\eqref{eq:ode_lambda} starting in
$\mathcal{U}_\gamma$ satisfies, for all $t\ge0$,
\begin{equation}
\|\mathbf{e}_x(t)\|^2+\|\mathbf{e}_\lambda(t)\|^2
\;\le\;\frac{V_{\mathrm{J}}(0)}{\underline{c}_J}
\,e^{-2\gamma t},
\label{eq:thm_exp}
\end{equation}
where $V_{\mathrm{J}}=\mathbf{z}^\top\mathbf{P}_J\mathbf{z}$
with the certificate matrix $\mathbf{P}_J$ (not the block-diagonal
candidate of \eqref{eq:VJ}), $\mathbf{z}\in\ker\mathbf{C}$, and
$\underline{c}_J$ is the
smallest eigenvalue of $\mathbf{P}_J$ on $\ker\mathbf{C}$; the KKT
pair is therefore locally exponentially stable with rate
$\gamma_J$; more precisely, every rate $\gamma<\gamma_J$ is attained,
and $\gamma_J$ is the supremum. The inequality \eqref{eq:J} is
sufficient, not necessary, for this conclusion: with a free
$\mathbf{P}_J$ it returns the spectral abscissa of $\mathbf{J}_J$ on
$\ker\mathbf{C}$.
\end{theorem}

\begin{proof}
Let $V_{\mathrm{J}}=\mathbf{z}^\top\mathbf{P}_J\mathbf{z}$. Along a
solution, \eqref{eq:VJ_derivative} gives
$\dot V_{\mathrm{J}}\le-\mathbf{z}^\top
\mathbf{M}(\mathbf{x}(t);\mathbf{P}_J)\mathbf{z}$.
At $\mathbf{x}^*$, \eqref{eq:J} gives
$\mathbf{M}(\mathbf{x}^*;\mathbf{P}_J)\succeq2\gamma_J\mathbf{P}_J$ on
$\ker\mathbf{C}$. Let $\delta_{\mathrm{cell}}$ be the distance from
$(\mathbf{x}^*,\boldsymbol{\lambda}^*)$ to the nearest switch
$\{\lambda_j=x_j\}$; under strict complementarity it equals
$\min_j\max(\lambda^*_j,x^*_j)/\sqrt{2}$, which is positive by the
strict-complementarity hypothesis. Inside the
ball $\|\mathbf{z}\|<\delta_{\mathrm{cell}}$ the selector is
identical to the equilibrium selector, so
$\mathbf{M}(\mathbf{x}(t);\mathbf{P}_J)$ depends on
$\mathbf{H}(\mathbf{x}(t))$ only. For quadratic costs
$\mathbf{H}$ is constant and hence
$\mathbf{M}(\mathbf{x}(t))=\mathbf{M}(\mathbf{x}^*)$ exactly;
for non-quadratic costs Assumption~\ref{assump:A4} gives
$\|\mathbf{H}(\mathbf{x}(t))-\mathbf{H}(\mathbf{x}^*)\|
\le L_2\|\mathbf{e}_x\|\le L_2\|\mathbf{z}\|$,
so
$\mathbf{M}(\mathbf{x}(t))=\mathbf{M}(\mathbf{x}^*)+O(\|\mathbf{z}\|)$
inside the cell. In both cases
$\mathbf{M}(\mathbf{x}(t))\succeq2\gamma\mathbf{P}_J$ whenever
$\|\mathbf{z}\|\le\delta_\gamma$ and $\gamma<\gamma_J$, with the
explicit radius $\delta_\gamma
=\min(\delta_{\mathrm{cell}}, 2(\gamma_J-\gamma)\underline{c}_J/\ell_M)$,
where $\ell_M$ is the Lipschitz constant of $\mathbf{M}$ for
non-quadratic costs
($\ell_M\le2\|\mathbf{P}_J\|\kappa_H L_2$, with $\kappa_H$ the
graph- and gain-dependent constant of
Appendix~\ref{app:VJ_expansion}); for
quadratic costs $\ell_M=0$ and the radius is
$\delta_{\mathrm{cell}}$.
Then $\dot V_{\mathrm{J}}\le-2\gamma V_{\mathrm{J}}$ on
$\{\|\mathbf{z}\|\le\delta_\gamma\}$, and Gronwall's inequality gives
\eqref{eq:thm_exp} as long as the trajectory stays there, which holds
for $\|\mathbf{z}(0)\|$ small since $V_{\mathrm{J}}$ dominates
$\|\mathbf{z}\|^2$ through $\underline{c}_J$.
\end{proof}

\begin{remark}[Verified joint rate]
\label{rem:joint_numbers}
For the instance of Section~\ref{sec:num_experiments}
($n=10$, $d=5$, $\mu=0.5$, $L=0.75$, $\eta=1$,
$\varepsilon=0.05$, $\beta=0.4$, $\rho=10$) Condition~\ref{cond:J}
is feasible up to $\gamma_J=0.2821$, which coincides with the spectral
abscissa of $\mathbf{J}_J$ on $\ker\mathbf{C}$ (of dimension $245$
here), because a free $\mathbf{P}_J$ makes \eqref{eq:J} equivalent to
$\mathbf{J}_J+\gamma_J I$ being Hurwitz on that subspace. Perturbations
confined to $\ker\mathbf{C}$ relax with a fitted state-norm rate of
$0.278$--$0.284$ over several amplitudes, reproducing the abscissa.
The relative primal error used in Fig.~\ref{fig:convex_convergence}
decays faster, and the modal decomposition says why. The six slowest
eigenmodes of $\mathbf{J}_J$ on $\ker\mathbf{C}$ lie in
$[-0.2868,-0.2821]$; the four slowest carry about $1\%$ of their mass
in the primal block and $71$--$72\%$ in the dual block, the remaining
two less than $0.3\%$ in the primal block. The slowest eigenvector with
a substantial primal component ($64\%$) has rate $0.3741$, the kernel
abscissa. The primal functional therefore follows the kernel rate over
the simulated horizon---$0.375$ fitted on $[7.2,45]$ of the same
curve---and the joint abscissa is invisible in it,
though sharp for
the state norm, because a $1\%$ primal projection needs a longer
horizon to overtake. The kernel abscissa
is $0.374$ and the fixed quadratic form certifies $\gamma_c=0.173$
there. Removing the projections changes the linearization: the boundary
coordinates are live instead of frozen, and the abscissa on
$\ker\mathbf{C}$ is $0.2821$ rather than the $0.3741$ of the projected
formulation. The rate remains exponential, and the projection-free
flow has an exact KKT limit from every initial condition with
$\mathbf{C}\mathbf{z}(0)=\mathbf{0}$; otherwise the limit is the
shifted KKT point of Lemma~\ref{lem:bias}.
\end{remark}

\begin{remark}[How the two layers separate]
\label{rem:separation}
Table~\ref{tab:sweep} sweeps the data. The rate increases
monotonically in $\beta$ over $[0.1,0.5]$
($0.0995$, $0.1708$, $0.2325$, $0.2815$, $0.2920$): for
small $\beta$ it essentially tracks $\beta$, so the dual leak is the
bottleneck, and it saturates near $0.30$ once the primal and consensus
blocks take over. It also increases in $\rho$ ($0.1584$ at $\rho=1$,
$0.2150$ at $\rho=3$, $0.2815$ at $\rho=10$, $0.3013$ at $\rho=30$);
at $\rho=1$ consensus no longer separates from the kernel. The rate is
insensitive to $\varepsilon$ to four digits over $[0.01,0.2]$
($0.2815$ in all four cases), so the momentum time scale affects
conditioning rather than the rate. It is nearly flat in $L/\mu$
($0.2804$, $0.2815$, $0.2828$, $0.2835$ for $1.2$, $1.5$, $2.0$,
$3.0$), so curvature is not the binding factor at this scale. For
$n\le10$ the rate is set by the optimization block ($0.2842$ at $n=5$,
$0.2815$ at $n=10$); for $n\ge20$ it follows $\nu_2$
($0.0773\approx0.79\,\nu_2$ at $n=20$ and
$0.0098\approx0.40\,\nu_2$ at $n=40$), so sparse graphs are capped by
consensus.
\end{remark}

\emph{Scope of the global claim.} The kernel argument above uses a
fixed $\mathbf{P}_c$ together with an active-set selector. In the
distributed flow the selector realizing \eqref{eq:selector} changes
with the state, the Hessian enters through
$\mathbf{H}(\mathbf{x}(t))$, and the graph stiffness
$\rho\mathbf{L}/\varepsilon$ scales the blocks differently, so the
selector-fixed form is stated at the KKT point and yields the local
result with the explicit radius $\delta_\gamma$. A global
functional is obtained not by enlarging that selector argument but
by changing coordinates: Section~\ref{subsec:global_lure} removes
the state-dependent selector altogether, writing the dual
nonlinearity as a shifted positive part with an exact sector that
holds at every state. The resulting LMI
(Condition~\ref{cond:global}) certifies global exponential
stability for quadratic costs, and on the nominal instance an explicit
quadratic form certifies a rate within five percent of the Jacobian
abscissa. The Jacobian LMI is
retained alongside it because it is sharp at the KKT point and its
abscissa is the rate that the full-trajectory tail fits reproduce,
to within $0.6$--$1.3\%$ in $\varepsilon$ and $L/\mu$ and to within
$11\%$ in $\beta$.

\section{Numerical Experiments}
\label{sec:num_experiments}

All simulations use RK4 integration with
step size $\Delta t = 2\times10^{-3}$
($\Delta t = 10^{-3}$ for $d=30$).
No projection or clipping is applied at any step: the
integrator advances the ODE \eqref{eq:ode_m}--\eqref{eq:ode_lambda}
verbatim, so $\mathbf{x}$ and $\boldsymbol{\lambda}$ may be negative
during transients and enter the orthant only through the equilibrium
equations (Proposition~\ref{prop:equilibrium}).
The baselines (DGD, DIGing and the PNN without tracking) do use the
orthant projection, since they require a feasible iterate.
Unless otherwise stated, the algorithm parameters
are $\eta = 1.0$, $\varepsilon = 0.05$,
$\beta = 0.4$, $\rho = 10$, $\gamma_s = 0.5$ (the
$\boldsymbol{\sigma}$ damping gain, not one of the rates
$\gamma_c$, $\gamma_J$, $\gamma_g$), and $\kappa = 0.05$.
Each agent draws its initial state independently,
$\mathbf{x}_i(0)\sim U[0,2]^d$, with no shared initialization; the
remaining states start at
$\mathbf{m}_i(0)=\boldsymbol{\sigma}_i(0)=\mathbf{0}$,
$\mathbf{q}_i(0)=\nabla f_i(\mathbf{x}_i(0))$ and
$\boldsymbol{\lambda}_i(0)=[\kappa\mathbf{1}-\mathbf{x}_i(0)]_+$;
the positive part here is an initialization only, not a projection
step of the dynamics.
Assumption~\ref{assump:A3} therefore holds by construction and,
for this quadratic instance, $\mathbf{C}\mathbf{z}(0)=\mathbf{0}$,
so every run lies on the slice
of Theorem~\ref{thm:global_lure}.
Every configuration is repeated over three independent (instance,
initialization) seed pairs and reported as a median with its
$[\min,\max]$ range; the ablation figure shows these medians as curves
with $[\min,\max]$ bands, the other figures show one representative seed.
The step size is integrator-safe on the scalability instance (the ratio
to the RK4 real-axis boundary is $0.57$), and halving $\Delta t$
changes terminal values only in the ninth significant digit; the
coupling sweep is checked separately below. On the $\varepsilon$ sweep the
matrix-norm ratio is conservative---it reaches $8.01$ at
$\varepsilon=0.01$---while the largest eigenvalue modulus gives
$\Delta t\max|\lambda(\mathbf{J})|\le0.20$ and halving $\Delta t$ twice
leaves the fitted rate unchanged to six digits. The decay rate is fitted on the state norm
$\|\boldsymbol{\zeta}(t)-\boldsymbol{\zeta}^*\|$ over the last
quarter $[0.75T,T]$ of the horizon, where the slow mode dominates
and the fit is clean ($R^2\ge0.999$ in every run); scanning
$[0.70T,T]$, $[0.75T,T]$, $[0.80T,T]$ across seven configurations
gives differences of at most $4\times10^{-3}$ (at most
$1\times10^{-3}$ in six of the seven rows), all inside the seed
ranges, so the reported window is representative, not tuned. A pre-registered
fit on $[0.12T,0.75T]$ (after the initial peak) overestimates the
rate because it still contains the fast tracker residual; it is
kept only as a finite-horizon diagnostic. Primal and consensus
rates are reported as references only. The network is a ring graph
with unit edge weights ($w_{ij}=1$), spectral gap
$\nu_2=2(1-\cos(2\pi/n))$, the sparsest connected topology and
therefore the least favorable case for consensus.
Supplementary material provides the Python implementation that
reproduces every figure.

We consider the distributed optimization
problem
\[
\min_{\mathbf{x} \ge 0}
\frac{1}{n}\sum_{i=1}^n f_i(\mathbf{x})
+ \kappa\|\mathbf{x}\|_1,
\]
where each agent $i$ holds a private cost
function $f_i: \mathbb{R}^d \to \mathbb{R}$.
We test a convex quadratic instance
$f_i(\mathbf{x})
= \frac{1}{2}(\mathbf{x} - \mathbf{c}_i)^\top
H_i (\mathbf{x} - \mathbf{c}_i)$ with
$H_i \succ 0$, $\mu = 0.5$, and $L = 0.75$,
which models a sparse nonnegative regression
problem with heterogeneous local data.

\subsection{Convex Results}
\label{subsec:exp_convex}

We first validate the rates predicted by the theory, then study the
sparsity of the limit, and then scalability, the coupling strength,
the gradient-tracking ablation, the comparison with other neurodynamic
solvers and the block ablations, and the effect of a mean initialization
defect.

\subsubsection{Exponential Convergence and Rates}

\begin{figure}[ht]
\centering
\includegraphics[width=\columnwidth]{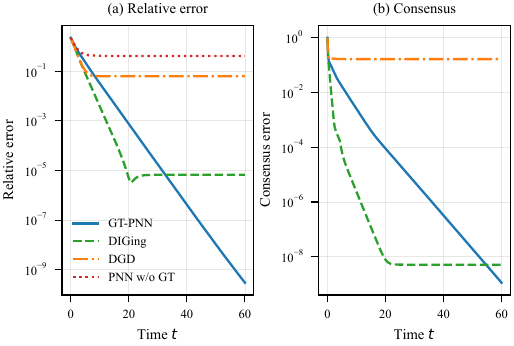}

\caption{Exponential convergence on convex quadratic
($n{=}10$, $d{=}5$, $L/\mu{=}1.5$, per-agent IC):
relative error (a) and consensus error (b) for projection-free
GT-PNN, continuous-time DIGing, projected DGD, and PNN without
gradient tracking. Both GT-PNN errors decrease from $t=0$; its
state-norm rate matches the certified Jacobian abscissa to within
$1\%$ on the tail window, and the primal scalar decays at the
faster kernel abscissa.}
\label{fig:convex_convergence}
\end{figure}

Fig.~\ref{fig:convex_convergence} verifies the exponential
convergence of the projection-free GT-PNN on the convex quadratic
($n = 10$, $d = 5$, $L/\mu = 1.5$). The functional that is directly
comparable with the certificates is the norm of the error state on
$\ker\mathbf{C}$: after a short transient caused by the tracker
initialization---$\mathbf{q}_i(0)=\nabla f_i(\mathbf{x}_i(0))$ is far
from $\bar{\mathbf{q}}^*$---it decays exponentially and reaches
$3.7\times10^{-8}$ at $t=60$. Two estimates are recorded. On the
pre-registered window $[7.2,45]$ the fit gives $0.3051$ ($R^2=0.998$;
median $0.3051$, range $[0.2949,0.3109]$); this window still contains
the residual of the fast tracker mode, and the estimate is $8\%$ above
the Jacobian abscissa. On the last quarter $[45,60]$ the fit gives
$0.2832$ ($R^2=1.000$; median $0.2832$, range $[0.2822,0.2835]$),
i.e.\ $0.4\%$ above $\gamma_J$. The three instance certificates are
$0.2821$, $0.2815$ and $0.2815$ (median $0.2815$); on each instance
the LMI value coincides with the spectral abscissa to seven digits.
Both estimates lie above the global rate $\gamma_g=0.27$, as a lower
bound requires. The primal
relative error
$\|\bar{\mathbf{x}}(t)-\bar{\mathbf{x}}^*\|
/\|\bar{\mathbf{x}}^*\|$ decays faster---$0.375$ fitted on the
pre-registered window---because the slowest eigenmodes carry about
one percent of their mass in the primal block
(Remark~\ref{rem:joint_numbers}); after $60$ time units it is
$3.0\times10^{-10}$ and the consensus error $1.1\times10^{-9}$, the
latter decreasing monotonically from $1.0$--$1.3$ at $t=0$.
Both errors decay exponentially, which is the statement of
Theorem~\ref{thm:exponential_convergence}.
No hump appears at the start: with per-agent initialization the
consensus error is already $O(1)$ at $t=0$ and decreases from there,
whereas a shared initialization starts on the consensus manifold and
produces a spurious rise.
The limit is
$\bar{\mathbf{x}}^*
=(0.6959,\,0.1858,\,0,\,0.3463,\,0)$ with
$\bar{\boldsymbol{\lambda}}^*
=(0,\,0,\,0.0443,\,0,\,0.3130)$: coordinates $3$ and $5$ are pinned
to the boundary by a positive multiplier and the other three are
interior, so the sparsity pattern is produced by the dual and not by
a projection.
The certificate is local, with the explicit radius
$\delta_\gamma$ of Theorem~\ref{thm:exponential_convergence}; the
initial distance of the figure need not lie inside $\delta_\gamma$,
which is why the global certificate is used here as well. On the
last quarter the modal decomposition of
Remark~\ref{rem:joint_numbers} applies: the state norm relaxes at
$\gamma_J$, while the primal scalar functional follows the kernel
abscissa $0.3741$. On the same instance the global Lur'e certificate
of Condition~\ref{cond:global} holds with $\gamma_g=0.27$
(Remark~\ref{rem:global_numbers}), so the decay is certified from
every initial condition with $\mathbf{C}\mathbf{z}(0)=\mathbf{0}$---in
particular every initialization satisfying
Assumption~\ref{assump:A3}---and not only inside $\delta_\gamma$.
Because no coordinate is clipped, primal and dual entries may cross
zero during the transient; they return and settle on the boundary,
so feasibility is recovered at the equilibrium and not enforced along
the way.
As a baseline we use continuous-time DIGing with
the orthant projection: it tracks gradients
but has no dual dynamics.
It reaches $6.8\times10^{-6}$ (median, range
$[2.5\times10^{-6},\,3.1\times10^{-5}]$), whereas projected
DGD and the PNN without gradient tracking stagnate
at $6.3\times10^{-2}$ (range
$[5.1\times10^{-2},\,2.8\times10^{-1}]$) and
$4.5\times10^{-1}$ (range
$[4.1\times10^{-1},\,1.4]$), since neither cancels the
heterogeneity of the local costs.

\begin{figure}[ht]
\centering
\includegraphics[width=\columnwidth]{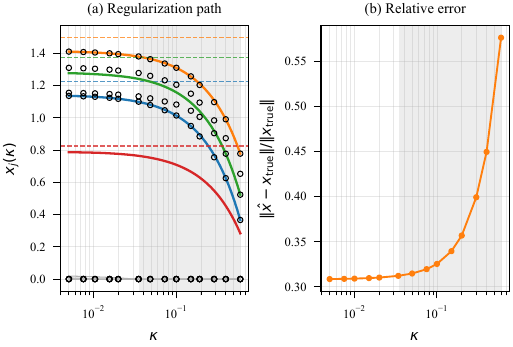}
\caption{Sparse regularization path on the nonnegative sparse
regression instance ($n{=}8$, $d{=}8$, $m_i{=}12$, four nonzeros).
(a)~Coordinates of the limit versus $\kappa$: solid curves are the
coordinates in $\operatorname{supp}(\mathbf{x}_{\rm true})$, dashed
lines their true values, gray curves the remaining coordinates, and
circles the limits obtained by integrating
\eqref{eq:ode_m}--\eqref{eq:ode_lambda}
at the fourteen anchor values $\kappa=0.005,\dots,0.6$. (b)~Relative
distance to $\mathbf{x}_{\rm true}$; the shaded region is where the
centralized path recovers the support of $\mathbf{x}_{\rm true}$
($\kappa\gtrsim0.035$).}
\label{fig:sparse_path}
\end{figure}

\subsubsection{Sparsity, Complementarity, and the Path}
\label{subsubsec:exp_sparse}

Sparsity is a property of the equilibrium, so it is tested directly on a
nonnegative sparse regression instance of the same form,
$f_i(\mathbf{x})=\tfrac12\|\mathbf{A}_i\mathbf{x}-\mathbf{b}_i\|^2
+\tfrac{\mu}{2}\|\mathbf{x}\|^2$ with
$\mathbf{b}_i=\mathbf{A}_i\mathbf{x}_{\rm true}
+\sigma\boldsymbol{\xi}_i$, $n=8$, $d=8$, $m_i=12$ samples per agent,
and $\mathbf{x}_{\rm true}$ supported on four of the eight coordinates.
Fig.~\ref{fig:sparse_path} traces the regularization path at fourteen
values of $\kappa$ from $0.005$ to $0.6$, each integrated from three
(instance, initialization) seed pairs. Across the $42$ runs the limit
of the flow coincides with the centralized NN-LASSO solution: the
relative distance
$\|\bar{\mathbf{x}}-\mathbf{x}_{\rm cen}\|/\|\mathbf{x}_{\rm cen}\|$
lies between $9.6\times10^{-11}$ and $2.2\times10^{-7}$, the KKT
residual between $1.2\times10^{-13}$ and $6.7\times10^{-13}$, the
complementarity residual between $1.5\times10^{-15}$ and
$1.9\times10^{-14}$, and the consensus error between
$1.4\times10^{-14}$ and $4.3\times10^{-14}$. The support also coincides
for $\kappa\ge0.034$; at the four smaller anchors two of the three
seeds recover it and one retains a single extra coordinate, exactly as
the centralized path does. The markers in
Fig.~\ref{fig:sparse_path}(a) are these integrated limits, the curves
are the centralized path.

Exact recovery of the support of $\mathbf{x}_{\rm true}$ holds from
$\kappa\approx0.035$ in the convex program---the extra coordinate falls
from $7.4\times10^{-3}$ at $\kappa=0.02$ to $8.1\times10^{-7}$ at
$\kappa=0.035$---and from $\kappa=0.034$ among the fourteen integrated
anchors. Below it the $\ell_1$ path retains one extra coordinate (at
$\kappa=0.02$ the centralized solution has
$x_{{\rm cen},7}=7.4\times10^{-3}$), so the discrepancy is a property
of the convex program and not of the dynamics, which reproduces the
program at every $\kappa$. This is structural rather than numerical:
by Corollary~\ref{cor:oblique} the path is a straight line in the
$\bar{\mathbf{H}}$-metric, so the support changes only at the
breakpoints where a coordinate of $\mathbf{z}_\kappa$ crosses the
orthant boundary and is constant between two consecutive
breakpoints. The measured threshold $\kappa\approx0.035$ is one of
these breakpoints of the centralized path, not a property of the
integration.

At $\kappa=0.1$ the zeros are produced by the dual. The limit has
exactly four active coordinates, $\{1,3,6,8\}$, all of order one (the
smallest is $0.709$), while $\{2,4,5,7\}$ are pinned to zero by positive
multipliers; the complementarity residual is
$\max_j|x_j\lambda_j|\in[1.6,6.0]\times10^{-15}$, the KKT residual
$3.1$--$4.4\times10^{-13}$ and the consensus error
$2.0$--$4.2\times10^{-14}$. The distance to the centralized solution is
at most $1.9\times10^{-8}$, so it sits on the support coordinates and
not on a weak threshold.

\subsubsection{Scalability}

\begin{figure}[ht]
\centering
\includegraphics[width=\columnwidth]{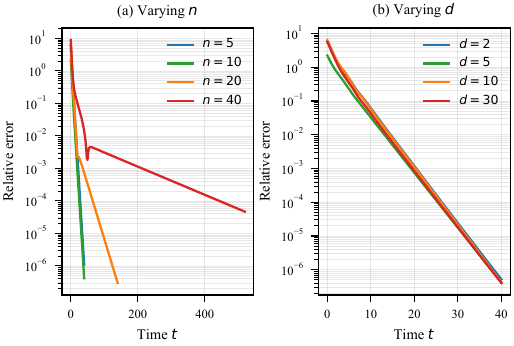}

\caption{Scalability on the convex quadratic ($L/\mu{=}1.5$, ring).
(a)~Relative error versus time for $n\in\{5,10,20,40\}$:
the decay slows as $\nu_2$ shrinks.
(b)~Relative error versus time for $d\in\{2,5,10,30\}$: the rate is
essentially independent of $d$.}
\label{fig:scalability}
\end{figure}

Fig.~\ref{fig:scalability} shows how
convergence scales with network size $n$
and dimension $d$.
For $n \le 10$ the tail state-norm rate is $0.286$ at $n=10$ and
$0.339$ (range $[0.312,0.356]$) at $n=5$; the certified rate is
$0.282$--$0.284$, i.e., the spectral abscissa of the Jacobian on
$\ker\mathbf{C}$ reported in Remark~\ref{rem:joint_numbers}. At
$n=5$ the slow mode carries little mass over the recorded horizon,
so the estimate stays $19\%$ above the abscissa; at $n=10$ it is
$1.5\%$ above.
As $n$ grows, the spectral gap of the ring
graph shrinks ($\nu_2 = 1.38, 0.382, 0.098,
0.0246$ for $n = 5, 10, 20, 40$), and
consensus becomes the binding block of the
certificate: the certified rate drops to
$0.0773$ ($n = 20$, $0.79\,\nu_2$) and $0.0098$
($n = 40$, $0.40\,\nu_2$), the measured tail rates being
$0.0780$ (range $[0.0775,0.0792]$) and $0.0099$
(range $[0.0099,0.0100]$), i.e.\ $0.9\%$ and $1.0\%$ above the
certificate.
This is the two-layer structure in action: the
optimization rate is set by the kernel, which follows
the curvature and the dual leak, while on sparse graphs
the consensus block caps the certified rate at a
fraction of $\nu_2$ and becomes the bottleneck.
For $n = 40$ the relative error
decays monotonically from $8.9$ to
$7.3\times10^{-5}$ (median, range
$[4.5\times10^{-5},1.8\times10^{-4}]$) over $t \in [0,520]$; this
terminal value is not a step-size artifact, since
$\Delta t=2\times10^{-3}$ gives $\|\mathbf{J}\|_2\Delta t=1.60$
against the RK4 boundary $2.785$, and halving $\Delta t$ changes the
terminal value in the ninth significant digit. The earlier
$n=20,40$ runs used $\Delta t=5\times10^{-3}$, for which
$\|\mathbf{J}\|_2\Delta t=4.00$ exceeds that boundary, so their
terminal values carried an integrator error and are superseded.
Increasing $d$ leaves the rate essentially
unchanged (Fig.~\ref{fig:scalability}(b)):
the tail rates are $0.284$, $0.286$, $0.284$ and $0.286$ for
$d=2,5,10,30$; the first three lie within $1.5\%$ of their
certificates, and the $d=30$ value, for which no certificate matrix
is formed, matches that range. The final relative errors are of order
$10^{-7}$, because the curvature range is fixed by
$\mu$ and $L$ independently of $d$. For $d=2$ one seed gives
$\bar{\mathbf{x}}^*=\mathbf{0}$, for which the relative error is not
defined; that run is reported through its state-norm rate ($0.283$)
only.

\subsubsection{Parameter Sensitivity and Ablation}

\begin{figure}[ht]
\centering
\includegraphics[width=\columnwidth]{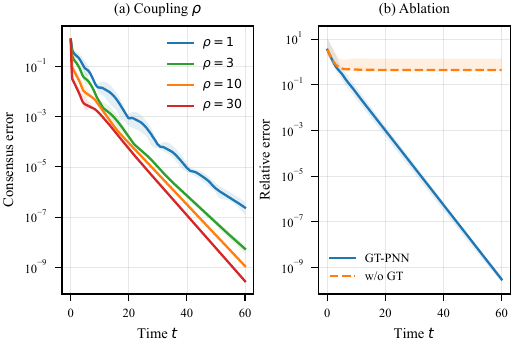}

\caption{Coupling strength and gradient-tracking ablation
($n{=}10$, $d{=}5$, convex quadratic, $L/\mu{=}1.5$, per-agent
initial conditions $x_i(0)\sim U[0,2]^d$, independent across agents).
Median over three instance seeds and two initialization seeds;
shaded bands: $[\min,\max]$.
(a)~Consensus error for different $\rho$: it starts at $O(1)$ and
decreases with $\rho$, while the transient to a $10^{-2}$ primal
error is nearly unchanged.
(b)~Without gradient tracking the error stagnates.}
\label{fig:param_ablation}
\end{figure}

Fig.~\ref{fig:param_ablation} studies the
effect of the coupling strength and the
contribution of gradient tracking.

\emph{Coupling strength \(\rho\).}
Fig.~\ref{fig:param_ablation}(a) shows the
consensus error for $\rho \in \{1,3,10,30\}$.
The consensus residual at $t=60$ decreases
monotonically with $\rho$, from
$2.4\times10^{-7}$ at $\rho=1$ to
$2.8\times10^{-10}$ at $\rho=30$.
The terminal relative errors are
$1.3\times10^{-8}$ at $\rho=1$ and between
$5.6\times10^{-10}$ and $3.0\times10^{-10}$ for $\rho\ge3$; apart from
that first row, the relative-error curves stay within a factor of two, and
the time to reach a relative error of
$10^{-2}$ is governed by the primal block, which changes little
with $\rho$ ($13.8$--$13.9$ across the four values).
The certified rate, however, is not $\rho$-independent: it rises
from $0.1584$ at $\rho=1$ to $0.2815$ at $\rho=10$ and $0.3013$
at $\rho=30$ (Table~\ref{tab:sweep}). The measured tail rate
rises from $0.166$ to $0.303$, i.e.\ by $0.137$ against a
certified increment of $0.143$; on the medians measured and certified
rates therefore agree to $4.5\%$, $3.6\%$, $0.6\%$ and $0.5\%$ at
$\rho=1,3,10,30$, and the certified increment overstates the measured
one by $1.04$. (The pre-registered finite window gives
$0.203$--$0.320$ and a factor $1.2$; the difference is the tracker
transient retained by that window.) The primal relative error
also rises ($0.367$ to $0.377$) but saturates at the kernel abscissa
$0.374$, which is why the primal functional is not used as the
certificate metric. At $\rho=30$ the matrix-norm product
$\|\mathbf{J}\|_2\Delta t=4.80$ exceeds the RK4 real-axis boundary
$2.785$, but the norm test is conservative there: the largest eigenvalue
modulus gives $\Delta t\max|\lambda(\mathbf{J})|=0.098$. Halving
$\Delta t$ changes the measured rate by less than
$0.1\%$ and the trajectory by $2\times10^{-6}$ in relative terms,
so the values at that end of the sweep are step-converged rather than
integrator-limited.

\emph{Gradient tracking ablation.}
Fig.~\ref{fig:param_ablation}(b) compares
GT-PNN with a variant without gradient
tracking, where $\mathbf{q}_i$ is replaced
by the local gradient $\nabla f_i(\mathbf{x}_i)$.
Without tracking, the equilibrium satisfies
only a perturbed KKT system, and the
relative error stagnates at
$4.5\times10^{-1}$ (median over six runs, range
$[4.1\times10^{-1},1.4]$).
With tracking, it reaches
$3.0\times10^{-10}$ (range
$[2.4\times10^{-10},4.1\times10^{-10}]$), confirming that tracking is
essential for exact optimality under
heterogeneous local costs.

Table~\ref{tab:sweep} pairs, for the same instance, the certified rate
$\gamma_J$ of Condition~\ref{cond:J} with the measured state-norm rate.
The certified entries are offline LMI evaluations involving no
simulation; the measured entries are tail fits of the state norm on
full trajectories over six runs (three instance seeds and two
initializations), with $R^2\ge0.998$, and they lie above the certificates
in every column, as a lower bound requires. For $\varepsilon$,
$L/\mu$ and $\beta$ the two agree to $0.6\%$, $1.3\%$ and $11\%$;
the $\beta$ row is the most conservative, the gap shrinking from
$11\%$ at $\beta=0.3$ to $0.6\%$ at $\beta=0.4$ as the dual leak
ceases to be the limiting block.
The trends agree: the rate rises with $\beta$ and
saturates near $0.30$, rises with $\rho$, is insensitive to
$\varepsilon$ to four digits, is nearly flat in $L/\mu$, and is
capped by $\nu_2$ once the graph is sparse enough. The interpretation is given in
Remark~\ref{rem:separation}.

\begin{table}[ht]
\centering
\caption{Rate versus problem data ($n{=}10$, $d{=}5$, $\mu{=}0.5$ unless
stated). \emph{Cert.}: Jacobian abscissa $\gamma_J$ of
Condition~\ref{cond:J} (offline LMI at the KKT point). \emph{Measured}:
exponential fit of $\|\boldsymbol{\zeta}(t)-\boldsymbol{\zeta}^*\|$
on $[0.75T,T]$ ($R^2\ge0.998$), median over six runs (three instance
seeds and two initializations). All rows are full-trajectory tail
fits; the $n$ row gives certificates only, and its measured rates are
in the text.}
\label{tab:sweep}
\begin{tabular}{lcccc}
\hline
$\beta$ & $0.1$ & $0.2$ & $0.3$ & $0.4$ \\
$\gamma_J$ (cert.) & $0.0995$ & $0.1708$ & $0.2325$ & $0.2815$ \\
measured & $0.1053$ & $0.1809$ & $0.2583$ & $0.2832$ \\
\hline
$\rho$ & $1$ & $3$ & $10$ & $30$ \\
$\gamma_J$ (cert.) & $0.1584$ & $0.2150$ & $0.2815$ & $0.3013$ \\
measured & $0.166$ & $0.223$ & $0.283$ & $0.303$ \\
\hline
$\varepsilon$ & $0.01$ & $0.02$ & $0.05$ & $0.20$ \\
$\gamma_J$ (cert.) & $0.2815$ & $0.2815$ & $0.2815$ & $0.2815$ \\
measured & $0.2833$ & $0.2833$ & $0.2832$ & $0.2832$ \\
\hline
$L/\mu$ & $1.2$ & $1.5$ & $2.0$ & $3.0$ \\
$\gamma_J$ (cert.) & $0.2804$ & $0.2815$ & $0.2828$ & $0.2835$ \\
measured & $0.2841$ & $0.2832$ & $0.2848$ & $0.2855$ \\
\hline
$n$ & $5$ & $10$ & $20$ & $40$ \\
$\nu_2$ & $1.382$ & $0.382$ & $0.098$ & $0.0246$ \\
$\gamma_J$ & $0.2842$ & $0.2815$ & $0.0773$ & $0.0098$ \\
\hline
\end{tabular}
\end{table}

\emph{Local heterogeneity.} The certificate is evaluated at the KKT
point; its dependence on the individual Hessian blocks is weak. To
test it, the local costs are deformed as
$\mathbf{H}_i=\bar{\mathbf{H}}+\mu_s(\mathbf{H}_i^0-\bar{\mathbf{H}})$
and
$\mathbf{c}_i=\bar{\mathbf{c}}+\mu_s(\mathbf{c}_i^0-\bar{\mathbf{c}})$,
with the spectrum clipped to $[\mu,L]$, for a mismatch parameter
$\mu_s\in\{0,0.5,1,2\}$ on the same three seed pairs. At $\mu_s=0$
all agents share one local cost; at $\mu_s=2$ the deviation is
doubled. The pre-registered finite window gives $0.349$ (range
$[0.325,0.371]$) at $\mu_s=0$ and $0.294$ (range $[0.292,0.296]$) at
$\mu_s=2$, a drop of $16\%$; the tail fit gives $0.289$ (range
$[0.285,0.306]$) and $0.283$ (range $[0.282,0.284]$), a change of
$1.9\%$, smaller than the $6.1\%$ seed-to-seed spread of the
$\mu_s=0$ row. Heterogeneity therefore
costs transient time rather than asymptotic rate, consistent with the
certificate moving only from $0.2823$ to $0.2814$. The flow remains
exact at the KKT point in every case, with terminal relative errors
between
$2.0\times10^{-10}$ and $8.6\times10^{-10}$.

\subsubsection{Comparison with Other Neurodynamic Solvers and Block Ablations}
\label{subsec:neuro_compare}

\begin{figure}[ht]
\centering
\includegraphics[width=\columnwidth]{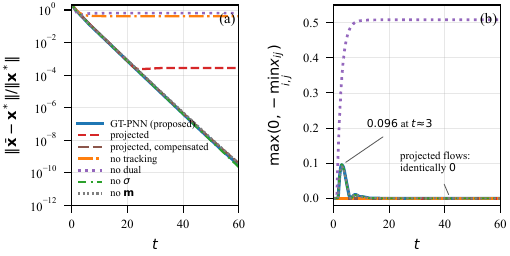}

\caption{The seven flows of Table~\ref{tab:neuro} on the nominal instance.
(a)~Relative primal error: the flows without tracking ($\mathrm{N}_3$) or
without a dual ($\mathrm{N}_4$) stagnate, and the clipped flow without
compensation ($\mathrm{N}_1$) stalls near $10^{-4}$; the
compensated clipped flow ($\mathrm{N}_2$) coincides with the proposed one.
(b)~Violation of the orthant: the clipped flows are feasible by
construction, the projection-free flows return to the orthant after a
transient, and the flow without a dual stays outside it.}
\label{fig:neuro_compare}
\end{figure}

The continuous-time literature disposes of the orthant in two standard
ways: projecting the primal iterate after each integration step, and
replacing the tracker by the local gradient. Both are applied here to
the same template, so that the comparison isolates the design choice
from the tuning. Table~\ref{tab:neuro} and Fig.~\ref{fig:neuro_compare}
report seven flows integrated with the same RK4 step, the same
$\Delta t$, the same output grid and the same per-agent initialization as
the main runs, over the same fifteen (five instances and three
initializations) pairs:
Table~\ref{tab:neuro} gives their medians, Fig.~\ref{fig:neuro_compare}
one representative run.
Besides the proposed flow we integrate (i)~$\mathrm{N}_1$, the same flow
with $\mathbf{x}$ and $\boldsymbol{\lambda}$ clipped to the orthant after
every step; (ii)~$\mathrm{N}_2$, the same clipped flow in which the
tracker is re-aligned with the clipped iterate,
$\mathbf{q}_i\mapsto\mathbf{q}_i+\nabla f_i(\mathbf{x}_i^{+})
-\nabla f_i(\mathbf{x}_i^{-})$, whenever the clip acts, which is the
correction used by projected DIGing; (iii)~$\mathrm{N}_3$, a clipped flow
without gradient tracking, i.e.\ the variant of
Fig.~\ref{fig:param_ablation}(b); and (iv)~$\mathrm{N}_4$, a projection-free
tracking flow with the dual block removed. The last two rows are block
ablations of the proposed flow, $\boldsymbol{\sigma}\equiv\mathbf{0}$ and
$\mathbf{m}\equiv\mathbf{0}$.

\begin{table}[ht]
\centering
\caption{Solvers and block ablations on the convex quadratic, with the
same instance, integrator, step and initialization as the main runs;
median over the fifteen (instance, initialization) pairs.
\emph{rel}: relative primal error at $t=60$; \emph{viol}:
$\max(0,-\min_{i,j}x_{ij})$ over the horizon; \emph{KKT} and
\emph{comp}: stationarity and complementarity residuals at $t=60$;
$t_{\mathrm{act}}$: last time from which the zero pattern of
$\bar{\mathbf{x}}$ matches that of $\mathbf{x}^*$; \emph{supp}: runs
whose active set equals the centralized one (the five instances have
$3$, $3$, $1$, $1$ and $4$ nonzeros). Dashes: no dual state.}
\label{tab:neuro}
\scriptsize
\setlength{\tabcolsep}{2.6pt}
\begin{tabular}{lcccccc}
\hline
flow & \emph{rel} & \emph{viol} & \emph{KKT} & \emph{comp} &
$t_{\mathrm{act}}$ & \emph{supp}\\
\hline
GT-PNN & $3.0{\times}10^{-10}$ & $8.8{\times}10^{-2}$ &
$3.7{\times}10^{-11}$ & $1.3{\times}10^{-9}$ & $17.22$ & $15/15$\\
$\mathrm{N}_1$ & $6.4{\times}10^{-4}$ & $0$ & $4.8{\times}10^{-3}$ &
$8.7{\times}10^{-10}$ & $16.52$ & $15/15$\\
$\mathrm{N}_2$ & $6.3{\times}10^{-10}$ & $0$ & $1.4{\times}10^{-9}$ &
$5.8{\times}10^{-10}$ & $16.12$ & $15/15$\\
$\mathrm{N}_3$ & $1.2$ & $0$ & $2.1{\times}10^{-3}$ &
$2.9{\times}10^{-1}$ & --- & $0/15$\\
$\mathrm{N}_4$ & $1.9$ & $5.1{\times}10^{-1}$ & --- &
--- & --- & $0/15$\\
$-\boldsymbol{\sigma}$ & $2.8{\times}10^{-10}$ & $8.7{\times}10^{-2}$ &
$4.1{\times}10^{-11}$ & $1.9{\times}10^{-10}$ & $17.22$ & $15/15$\\
$-\mathbf{m}$ & $4.1{\times}10^{-10}$ & $8.4{\times}10^{-2}$ &
$5.5{\times}10^{-11}$ & $1.3{\times}10^{-9}$ & $16.70$ & $15/15$\\
\hline
\end{tabular}
\end{table}

\emph{Tracking and the dual are load-bearing.} $\mathrm{N}_3$ stalls at
$1.2$ in relative error (median; up to $1.7$ on the five-instance grid,
whose instances are harder than the three used in
Fig.~\ref{fig:param_ablation}(b), where the same variant reads
$4.5\times10^{-1}$), leaves the consensus error
at $1.2\times10^{-1}$, and its complementarity residual is
$2.9\times10^{-1}$ although its stationarity residual is already
$2.1\times10^{-3}$: stationarity alone is not a certificate, because the
multiplier absorbs the gradient mismatch, and it keeps all five
coordinates active in every run. $\mathrm{N}_4$ is the sharper test: it
preserves the tracker (consensus $4.8\times10^{-12}$, and
\eqref{eq:conservation_law} to $6.5\times10^{-14}$ on the nominal run) and
still never
becomes feasible, the violation sitting at $5.1\times10^{-1}$ and not
decaying, so its limit is the consensus minimizer of the penalized
\emph{unconstrained} problem, with all five coordinates active. In this
template nonnegativity is therefore neither a consequence of consensus
nor of the primal dynamics: it is produced at the equilibrium by the dual
through the one-sided threshold of Proposition~\ref{prop:onesided}, and
$\mathrm{N}_4$ shows that nothing else in the template reproduces it.

\emph{Projection.} $\mathrm{N}_1$ is feasible at every instant, yet its
relative error freezes at $6.4\times10^{-4}$ (range
$[1.2\times10^{-5},1.4\times10^{-3}]$) and its stationarity residual at
$4.8\times10^{-3}$, six and eight orders above the proposed flow. The
cause is visible in the conservation law: clipping is agent-wise, so
\eqref{eq:conservation_law} is violated by $5.3\times10^{-3}$ on the
nominal run, and the stationarity residual equals that violation to all
printed digits---the limit is a KKT point of a problem whose perturbation
is precisely the clip. Compensating the tracker, $\mathrm{N}_2$, restores
the identity and recovers exactness ($6.3\times10^{-10}$ relative error,
$1.4\times10^{-9}$ stationarity residual, and the centralized active set
in all fifteen runs), at the price of an event-triggered correction at
every step in which the clip acts. That correction turns the right-hand
side into a hybrid field: the excursion \eqref{eq:excursion} disappears,
but the smooth autonomous structure on which
Theorem~\ref{thm:global_lure} rests is lost. The projection-free design
therefore does not buy accuracy, which $\mathrm{N}_2$ also attains; it
buys a single smooth vector field, and it pays in transient feasibility:
the median excursion is $8.8\times10^{-2}$ ($9.6\times10^{-2}$ at
$t\approx3$ on the nominal run), it starts quadratically by
Proposition~\ref{prop:dichotomy}, and it is removed by the dual, the
violation at $t=60$ being $1.4\times10^{-9}$.

\emph{Cost.} A single RK4 stage of every flow in
Table~\ref{tab:neuro} costs exactly one local evaluation per agent,
either $\nabla f_i(\mathbf{x}_i)$ or $\mathbf{H}_i\dot{\mathbf{x}}_i$,
both of which are one $d\times d$ matvec on this quadratic instance,
so tracking is free in this
accounting and the measured cost is $4.00$ evaluations per step for
GT-PNN, $\mathrm{N}_3$ and $\mathrm{N}_4$. The clip itself is an $O(d)$
operation, but it acts on $97\%$ of the steps---the constraint is active
at the equilibrium in the sense that the inactive coordinates are pushed
against the boundary every step---so the compensation of $\mathrm{N}_2$
costs $4.97$ evaluations per step, a $24\%$ steady-state overhead rather
than a transient one. Measured in work to reach a relative error of
$10^{-6}$, the proposed flow needs $7.7\times10^{4}$ evaluations per
agent, $\mathrm{N}_2$ needs $9.7\times10^{4}$, and $\mathrm{N}_1$ never
arrives: its floor is a bias of order $\Delta t$, measured as $5.4$,
$2.8$, $1.4$ and $0.71$ times $10^{-4}$ at $\Delta t=4,2,1,0.5$ times
$10^{-3}$ (empirical order $0.96$--$0.99$), so $10^{-6}$ would require
$\Delta t\approx7\times10^{-6}$, i.e.\ some $4\times10^{2}$ times the work
of the proposed flow, and $10^{-9}$ is unreachable at any step. For the
two exact flows the error at $t=60$ is independent of $\Delta t$
($3.0$ and $3.9$ times $10^{-10}$ at all four step sizes), that is, it is
set by the horizon and the rate, not by the integrator, so the largest
stable step already delivers the accuracy.

\emph{Block ablations.} Dropping $\boldsymbol{\sigma}$ leaves the limit
unchanged ($2.8\times10^{-10}$, KKT $4.1\times10^{-11}$) and dropping
$\mathbf{m}$ leaves $4.1\times10^{-10}$ and $5.5\times10^{-11}$; both
identify the active set at $17.22$ and $16.70$ against $17.22$ for the
full flow, and both reproduce the centralized support in all fifteen
runs at the same work ($7.8\times10^{4}$ evaluations to $10^{-6}$ against
$7.7\times10^{4}$).
This is what Proposition~\ref{prop:gradient_tracking} predicts for
$\boldsymbol{\sigma}$, which is absent from the conservation law, and
what the separation of time scales predicts for $\mathbf{m}$, which
shapes only the fast mode. Their cost is transient: on the nominal run
the times to a relative error of $10^{-6}$ are $37.62$ and $38.20$
against $37.66$ for the full flow.

\subsubsection{Bias under a Mean Initialization Defect}

Table~\ref{tab:bias} tests Lemma~\ref{lem:bias} on the nominal instance,
whose state dimension is $250$ and whose $\ker\mathbf{C}$ has dimension
$245$. A mean defect $\boldsymbol{\mu}(0)$ is injected into $\mathbf{q}(0)$,
so that the trajectory starts outside $\ker\mathbf{C}$, and the limit is
compared with the reference run that satisfies
Assumption~\ref{assump:A3}. The displacement is linear in the defect
with the constant $1.90$ over mean defects up to $0.24$, i.e.\ up to
$\|\boldsymbol{\delta}\|=0.64$, and leaves the linear regime at $0.8$;
the component outside $\ker\mathbf{C}$ is conserved to
$3\times10^{-13}$. Lemma~\ref{lem:bias} guarantees a bound
proportional to $\|\mathbf{C}^\dagger\boldsymbol{\mu}(0)\|$; the
constant $1.90$ is the measured ratio, not the conservative $\kappa$
of \eqref{eq:bias}.

\begin{table}[ht]
\centering
\caption{Displacement of the limit under a mean tracker defect on the
nominal instance ($n=10$, $d=5$, state dimension $250$,
$\dim\ker\mathbf{C}=245$, nominal parameters, $T=120$,
$\Delta t=2\times10^{-3}$). Columns: defect, $\ell_2$ norm of
$\boldsymbol{\delta}=\mathbf{C}^\dagger\boldsymbol{\mu}(0)$,
displacement over $\|\boldsymbol{\delta}\|$, and the violation of
$(\mathbf{I}-\Pi_C)\mathbf{z}(t)=(\mathbf{I}-\Pi_C)\mathbf{z}(0)$.}
\label{tab:bias}
\begin{tabular}{cccc}
\hline
$\|\boldsymbol{\mu}(0)\|$ & $\|\boldsymbol{\delta}\|$ &
$d(\infty)/\|\boldsymbol{\delta}\|$ & conservation \\
\hline
$0.012$ & $3.2\times10^{-2}$ & $1.903$ & $7\times10^{-15}$ \\
$0.059$ & $1.6\times10^{-1}$ & $1.903$ & $3\times10^{-14}$ \\
$0.236$ & $6.4\times10^{-1}$ & $1.903$ & $1\times10^{-13}$ \\
$0.943$ & $2.5$ & $1.465$ & $3\times10^{-13}$ \\
\hline
\end{tabular}
\end{table}

\section{Conclusion}
\label{sec:conclusion}

This paper introduced a continuous-time
projection-free gradient-tracking neurodynamic
system for distributed nonnegative sparse
optimization, embedding gradient tracking
into the PNN with a perturbed dual dynamics
that ensures exact KKT satisfaction under
heterogeneous local costs.
Neither the primal nor the dual is projected:
both are free, may be negative during transients, and
nonnegativity, feasibility and complementary slackness
are produced by the equilibrium equations, where
$\min\{\boldsymbol{\lambda}^*,\mathbf{x}^*\}
=-\beta\rho\mathbf{L}\mathbf{x}^*$ and an extremal argument
force consensus with $\bar{\mathbf{x}}^*\ge\mathbf{0}$.
The vector field is locally Lipschitz, so every initial
condition admits a unique classical solution and no
Filippov convention is required.
For convex objectives, three matrix inequalities
certify convergence with computable rates. The
centralized kernel is certified by
$\mathbf{M}_c(\mathbf{S})\succeq2\gamma_c\mathbf{P}_c$ for every
active-set selector; for quadratic costs
Theorem~\ref{thm:global_kernel} yields global
exponential stability with the computed rate
$\gamma_c=0.173$ on the nominal instance, the $\min$ term
being bounded by the pointwise selector inequality at every state
rather than linearized.
The distributed flow is certified by an LMI on the Jacobian at
the KKT point, restricted to the invariant subspace of the mean
tracking error, with no gain product and no stability assumption
on an isolated subsystem. Its stability is local, with the
explicit radius $\delta_\gamma$ of
Theorem~\ref{thm:exponential_convergence}, and the certified
rate equals the spectral abscissa of the Jacobian on
$\ker\mathbf{C}$, $0.2821$ on the nominal instance;
it is insensitive to $\varepsilon$, and on sparse graphs
$\nu_2$ becomes rate-limiting, the measured rate being
$0.79\,\nu_2$ at $n=20$ and $0.40\,\nu_2$ at $n=40$. A third LMI, posed in
the exact Lur'e coordinates of the dual error, makes the linear
part constant and the dual positive part a shifted positive
part with a global sector; it yields global exponential
stability of the distributed flow for quadratic costs from
every initial condition with
$\mathbf{C}\mathbf{z}(0)=\mathbf{0}$---in particular every
gradient-tracking initialization---with an explicit bias
when that condition is dropped, and with the rate
$\gamma_g=0.27$ on the nominal instance---within five percent
of the Jacobian abscissa and with no restriction to a
neighbourhood of the KKT point. On that instance the measured
state-norm rate is within one percent of the Jacobian abscissa on the
tail; growing local heterogeneity slows the finite-horizon
transient by about fifteen percent, while the tail rate changes by
less than two percent.
The comparison with other neurodynamic solvers isolates the two blocks
that carry the design: dropping the dual leaves the flow infeasible and
convergent to the penalized unconstrained minimizer, dropping the tracker
stagnates at $1.2$ with spurious support, whereas the two
together give $3.0\times10^{-10}$ with the exact active set in every run.
Clipping the primal either biases the limit at $6.4\times10^{-4}$ or, when
the tracker is compensated, requires an event-triggered correction that
removes the smooth autonomous structure on which the certificates rest and
costs a quarter more evaluations per step.
Extensions to stochastic gradients and
time-varying topologies, with explicit rates,
remain open.

\appendices
\section{Term-by-Term Expansion of $\dot V_{\mathrm{J}}$}
\label{app:VJ_expansion}

This appendix supplies the bookkeeping behind
\eqref{eq:VJ_derivative}. Put $\mathbf{N}=\rho\mathbf{L}$ and
$\mathbf{G}=\mathbf{I}+\mathbf{N}$, so that \eqref{eq:VJ} reads
$V_{\mathrm{J}}
=\langle\mathbf{G}\mathbf{e}_x,\mathbf{e}_x\rangle
+|\mathbf{e}_\lambda|^2+\varepsilon|\mathbf{w}|^2
+\eta^{-1}|\mathbf{e}_q|^2+|\mathbf{e}_\sigma|^2$.

\emph{Cell identity.} Let the KKT point be strictly complementary
and let $\mathcal{C}$ be the cell of constant active set that
contains it. For each coordinate exactly one of $\lambda^*_j$,
$x^*_j$ is positive; set $(\mathbf{D}_\lambda)_{jj}=1$ when
$\lambda^*_j=0$ and $(\mathbf{D}_\lambda)_{jj}=0$ otherwise, and
$\mathbf{D}_x=\mathbf{I}-\mathbf{D}_\lambda$. The branch of the
$\min$ selected at the KKT point is the coordinate that vanishes
there, and it does not change inside $\mathcal{C}$; hence,
throughout $\mathcal{C}$,
\begin{equation}
\Phi_j(\boldsymbol{\lambda},\mathbf{x})
=(\mathbf{D}_\lambda)_{jj}e_{\lambda,j}
+(\mathbf{D}_x)_{jj}e_{x,j},
\qquad j=1,\ldots,nd,
\label{eq:app_cell}
\end{equation}
so that \eqref{eq:selector} holds with equality for the selector
fixed by the active set. Thus $\Phi$ is linear and $\dot
V_{\mathrm{J}}$ is exactly quadratic in $\mathbf{z}$ on
$\mathcal{C}$; the selector inequality is used to identify the
selector, not to bound a residual.

\emph{The five blocks.} With
$\mathbf{s}=\mathbf{e}_q-\mathbf{e}_\lambda
+\mathbf{N}\mathbf{e}_x$ and
$\dot{\mathbf{x}}=\mathbf{m}
=\mathbf{w}-\mathbf{s}/\eta$, the error dynamics
\eqref{eq:ode_m}--\eqref{eq:ode_lambda} read, exactly on
$\mathcal{C}$,
\begin{align}
\dot{\mathbf{e}}_x&=\mathbf{w}-\eta^{-1}\mathbf{s},
\label{eq:app_dx}\\
\varepsilon\dot{\mathbf{w}}&=-\eta\mathbf{w}
+\frac{\varepsilon}{\eta}\dot{\mathbf{s}},
\label{eq:app_dw}\\
\dot{\mathbf{e}}_\lambda
&=-(\beta\mathbf{I}+\mathbf{D}_\lambda)\mathbf{e}_\lambda
+\beta\mathbf{e}_q-\mathbf{D}_x\mathbf{e}_x,
\label{eq:app_dl}\\
\dot{\mathbf{e}}_q&=-\mathbf{L}(\mathbf{e}_q+\mathbf{e}_\sigma)
+\mathbf{H}\bigl(\mathbf{w}-\eta^{-1}\mathbf{s}\bigr),
\label{eq:app_dq}\\
\dot{\mathbf{e}}_\sigma&=-\mathbf{L}\mathbf{e}_x
-\gamma_s\mathbf{e}_\sigma,
\label{eq:app_ds}
\end{align}
with $\mathbf{H}=\mathbf{H}(\mathbf{x}(t))$. Differentiating the
five blocks of \eqref{eq:VJ} along
\eqref{eq:app_dx}--\eqref{eq:app_ds} gives
\begin{align}
2\langle\mathbf{G}\mathbf{e}_x,\dot{\mathbf{e}}_x\rangle
&=2\langle\mathbf{G}\mathbf{e}_x,\mathbf{w}\rangle
-\frac{2}{\eta}\langle\mathbf{G}\mathbf{e}_x,\mathbf{e}_q\rangle
+\frac{2}{\eta}\langle\mathbf{G}\mathbf{e}_x,\mathbf{e}_\lambda\rangle
\notag\\
&\quad-\frac{2}{\eta}\langle\mathbf{G}\mathbf{e}_x,
\mathbf{N}\mathbf{e}_x\rangle,
\label{eq:app_Tx}\\
2\langle\mathbf{e}_\lambda,\dot{\mathbf{e}}_\lambda\rangle
&=-2\beta|\mathbf{e}_\lambda|^2
+2\beta\langle\mathbf{e}_\lambda,\mathbf{e}_q\rangle
-2\langle\mathbf{e}_\lambda,
\mathbf{D}_\lambda\mathbf{e}_\lambda\rangle
\notag\\
&\quad-2\langle\mathbf{e}_\lambda,
\mathbf{D}_x\mathbf{e}_x\rangle,
\label{eq:app_Tl}\\
\frac{2}{\eta}\langle\mathbf{e}_q,\dot{\mathbf{e}}_q\rangle
&=-\frac{2}{\eta}\langle\mathbf{e}_q,\mathbf{L}\mathbf{e}_q\rangle
-\frac{2}{\eta}\langle\mathbf{e}_q,\mathbf{L}\mathbf{e}_\sigma\rangle
+\frac{2}{\eta}\langle\mathbf{e}_q,\mathbf{H}\mathbf{w}\rangle
\notag\\
&\quad-\frac{2}{\eta^2}\langle\mathbf{e}_q,\mathbf{H}\mathbf{e}_q\rangle
+\frac{2}{\eta^2}\langle\mathbf{e}_q,\mathbf{H}\mathbf{e}_\lambda\rangle
\notag\\
&\quad-\frac{2}{\eta^2}\langle\mathbf{e}_q,
\mathbf{H}\mathbf{N}\mathbf{e}_x\rangle,
\label{eq:app_Tq}\\
2\langle\mathbf{e}_\sigma,\dot{\mathbf{e}}_\sigma\rangle
&=-2\langle\mathbf{e}_\sigma,\mathbf{L}\mathbf{e}_x\rangle
-2\gamma_s|\mathbf{e}_\sigma|^2,
\label{eq:app_Ts}
\end{align}
and the momentum equation, the second of
\eqref{eq:app_dx}--\eqref{eq:app_ds}, is driven by
$\dot{\mathbf{s}}=\dot{\mathbf{e}}_q-\dot{\mathbf{e}}_\lambda
+\mathbf{N}\dot{\mathbf{e}}_x$, so that block contributes
\begin{equation}
2\varepsilon\langle\mathbf{w},\dot{\mathbf{w}}\rangle
=-2\eta|\mathbf{w}|^2
+\frac{2\varepsilon}{\eta}\langle\mathbf{w},
\boldsymbol{\Omega}\rangle,
\qquad \boldsymbol{\Omega}=\dot{\mathbf{s}},
\label{eq:app_Tw}
\end{equation}
where, substituting \eqref{eq:app_dl}--\eqref{eq:app_dq},
\begin{equation}
\begin{aligned}
\boldsymbol{\Omega}
&=\mathbf{H}\mathbf{w}+\mathbf{N}\mathbf{w}
-\mathbf{L}\mathbf{e}_\sigma\\
&\quad-\Bigl(\mathbf{L}+\frac{1}{\eta}\mathbf{H}
+\beta\mathbf{I}+\frac{1}{\eta}\mathbf{N}\Bigr)\mathbf{e}_q\\
&\quad+\Bigl(\frac{1}{\eta}\mathbf{H}+\beta\mathbf{I}
+\mathbf{D}_\lambda+\frac{1}{\eta}\mathbf{N}\Bigr)\mathbf{e}_\lambda\\
&\quad+\Bigl(\mathbf{D}_x-\frac{1}{\eta}\mathbf{H}\mathbf{N}
-\frac{1}{\eta}\mathbf{N}^2\Bigr)\mathbf{e}_x .
\end{aligned}
\label{eq:app_Omega}
\end{equation}

\emph{Collected matrix.} The five contributions
\eqref{eq:app_Tx}--\eqref{eq:app_Ts} and \eqref{eq:app_Tw} are
exhaustive, because \eqref{eq:VJ} is block diagonal and each block
is differentiated separately; every coefficient is quadratic in
$\mathbf{z}$ once \eqref{eq:app_cell} is used. Collecting the
coefficients of the 15 independent block monomials yields
$\dot V_{\mathrm{J}}=-\mathbf{z}^{\top}\mathbf{M}\mathbf{z}$
with $\mathbf{M}=\mathbf{M}^{\top}$ and blocks
\begin{align}
\mathbf{M}_{ww}&=2\eta\mathbf{I}
-\frac{2\varepsilon}{\eta}(\mathbf{H}+\mathbf{N}),
\notag\\
\mathbf{M}_{wx}&=-\mathbf{G}
-\frac{\varepsilon}{\eta}\mathbf{D}_x
+\frac{\varepsilon}{\eta^2}\mathbf{H}\mathbf{N}
+\frac{\varepsilon}{\eta^2}\mathbf{N}^2,
\notag\\
\mathbf{M}_{w\lambda}&=-\frac{\varepsilon}{\eta}
\Bigl(\frac{1}{\eta}\mathbf{H}+\beta\mathbf{I}
+\mathbf{D}_\lambda+\frac{1}{\eta}\mathbf{N}\Bigr),
\notag\\
\mathbf{M}_{wq}&=\frac{\varepsilon}{\eta}
\Bigl(\mathbf{L}+\frac{1}{\eta}\mathbf{H}
+\beta\mathbf{I}+\frac{1}{\eta}\mathbf{N}\Bigr)
-\frac{1}{\eta}\mathbf{H},
\notag\\
\mathbf{M}_{w\sigma}&=\frac{\varepsilon}{\eta}\mathbf{L},
\qquad
\mathbf{M}_{xx}=\frac{2}{\eta}(\mathbf{N}+\mathbf{N}^2),
\notag\\
\mathbf{M}_{x\lambda}&=\mathbf{D}_x-\frac{1}{\eta}\mathbf{G},
\qquad
\mathbf{M}_{xq}=\frac{1}{\eta}\mathbf{G}
+\frac{1}{\eta^2}\mathbf{N}\mathbf{H},
\notag\\
\mathbf{M}_{x\sigma}&=\mathbf{L},
\qquad
\mathbf{M}_{\lambda\lambda}
=2(\beta\mathbf{I}+\mathbf{D}_\lambda),
\notag\\
\mathbf{M}_{\lambda q}&=-\beta\mathbf{I}
-\frac{1}{\eta^2}\mathbf{H},
\qquad
\mathbf{M}_{\lambda\sigma}=\mathbf{0},
\notag\\
\mathbf{M}_{qq}&=\frac{2}{\eta}\mathbf{L}
+\frac{2}{\eta^2}\mathbf{H},
\qquad
\mathbf{M}_{q\sigma}=\frac{1}{\eta}\mathbf{L},
\notag\\
\mathbf{M}_{\sigma\sigma}&=2\gamma_s\mathbf{I},
\label{eq:app_M}
\end{align}
where an off-diagonal block in the upper triangle is displayed and
the lower one is its transpose, e.g.\
$\mathbf{M}_{xw}=\mathbf{M}_{wx}^{\top}$.

No term has been dropped. The curvature enters through
$\mathbf{H}$ in the blocks $(w,w)$, $(w,x)$, $(w,q)$,
$(x,q)$, $(\lambda,q)$ and $(q,q)$; the graph enters through
$\mathbf{N}$ in $(w,w)$, $(w,x)$, $(w,q)$, $(w,\sigma)$,
$(x,x)$, $(x,\sigma)$, $(q,q)$ and $(q,\sigma)$; the selector
enters through $\mathbf{D}_\lambda$ and $\mathbf{D}_x$ in
$(\lambda,\lambda)$, $(x,\lambda)$ and $(w,x)$; and the damping
gains enter as $\beta$ in $(\lambda,\lambda)$, $(w,\lambda)$,
$(w,q)$ and $(\lambda,q)$, and as $\gamma_s$ only in
$(\sigma,\sigma)$.

The equalities above hold on $\mathcal{C}$. Outside
$\mathcal{C}$, \eqref{eq:app_cell} is replaced by the one-sided
bound \eqref{eq:selector}, and the pairing of $\Phi$ with
$\mathbf{w}$ is bounded by
$|\Phi_j|\le|e_{\lambda,j}|+|e_{x,j}|$ followed by a Young step;
both changes contribute nonnegative quadratic terms, so
$\mathbf{M}$ is replaced by a matrix that is smaller in the
Loewner order. Theorem~\ref{thm:exponential_convergence} is
stated on the neighbourhood in which the trajectory stays in
$\mathcal{C}$, so the exact form \eqref{eq:app_M} is the one used.

\emph{Arbitrary certificate matrix.} For the certificate the
block-diagonal weights of \eqref{eq:VJ} are replaced by the
solution $\mathbf{P}_J$ of \eqref{eq:J}. Let $\mathbf{J}(\mathbf{x})$
be the Jacobian of \eqref{eq:app_dx}--\eqref{eq:app_ds} with
respect to $\mathbf{z}$; it is affine in $\mathbf{H}(\mathbf{x})$
and constant when the costs are quadratic. Since the right-hand
side of \eqref{eq:app_dx}--\eqref{eq:app_ds} is exactly linear in
$\mathbf{z}$ on $\mathcal{C}$, one has
\begin{equation}
\dot V_{\mathrm{J}}
=\mathbf{z}^{\top}\bigl(\mathbf{J}(\mathbf{x})^{\top}\mathbf{P}_J
+\mathbf{P}_J\mathbf{J}(\mathbf{x})\bigr)\mathbf{z}
=-\mathbf{z}^{\top}\mathbf{M}(\mathbf{x};\mathbf{P}_J)\mathbf{z},
\label{eq:app_general}
\end{equation}
which is \eqref{eq:VJ_derivative} with
\begin{equation}
\mathbf{M}(\mathbf{x};\mathbf{P}_J)
=-\bigl(\mathbf{J}(\mathbf{x})^{\top}\mathbf{P}_J
+\mathbf{P}_J\mathbf{J}(\mathbf{x})\bigr).
\label{eq:app_M_general}
\end{equation}
The display \eqref{eq:app_M} is the special case
$\mathbf{P}_J=\operatorname{blkdiag}(\varepsilon\mathbf{I},
\mathbf{G},\mathbf{I},\eta^{-1}\mathbf{I},\mathbf{I})$ in the
order $(w,x,\lambda,q,\sigma)$, and it identifies which block
carries each term of the candidate. At the KKT point
\eqref{eq:app_M_general} reduces to
$-\bigl(\mathbf{J}_J^{\top}\mathbf{P}_J
+\mathbf{P}_J\mathbf{J}_J\bigr)$, so \eqref{eq:J} is exactly
$\mathbf{M}(\mathbf{x}^*;\mathbf{P}_J)\succeq2\gamma_J\mathbf{P}_J$
on $\ker\mathbf{C}$, which is the inequality used in the proof of
Theorem~\ref{thm:exponential_convergence}. Finally, for
non-quadratic costs $\mathbf{J}$ depends on the state only through
$\mathbf{H}$, whence
\begin{equation}
\begin{aligned}
\|\mathbf{M}(\mathbf{x};\mathbf{P}_J)
-\mathbf{M}(\mathbf{x}^*;\mathbf{P}_J)\|
&\le2\|\mathbf{P}_J\|\,\|\mathbf{J}(\mathbf{x})
-\mathbf{J}_J\|\\
&\le2\|\mathbf{P}_J\|\kappa_H L_2\|\mathbf{e}_x\|,
\end{aligned}
\label{eq:app_lip}
\end{equation}
where
$\kappa_H=8(1+\eta^{-2})(1+\|\mathbf{N}\|)$ collects the eight
blocks through which $\mathbf{H}$ enters
$\mathbf{J}$, namely $(q,w)$, $(q,q)$, $(q,\lambda)$, $(q,x)$,
$(w,w)$, $(w,q)$, $(w,\lambda)$ and $(w,x)$, each with a
coefficient bounded by $(1+\eta^{-2})(1+\|\mathbf{N}\|)$.
Hence the Lipschitz constant
$\ell_M\le2\|\mathbf{P}_J\|\kappa_H L_2$ used in the proof
covers all blocks of \eqref{eq:app_M}.

%% TEMP FILE FOR IEEE BIOGRAPHIES
\begin{IEEEbiography}[{\includegraphics[width=1in,height=1.25in,clip,keepaspectratio]{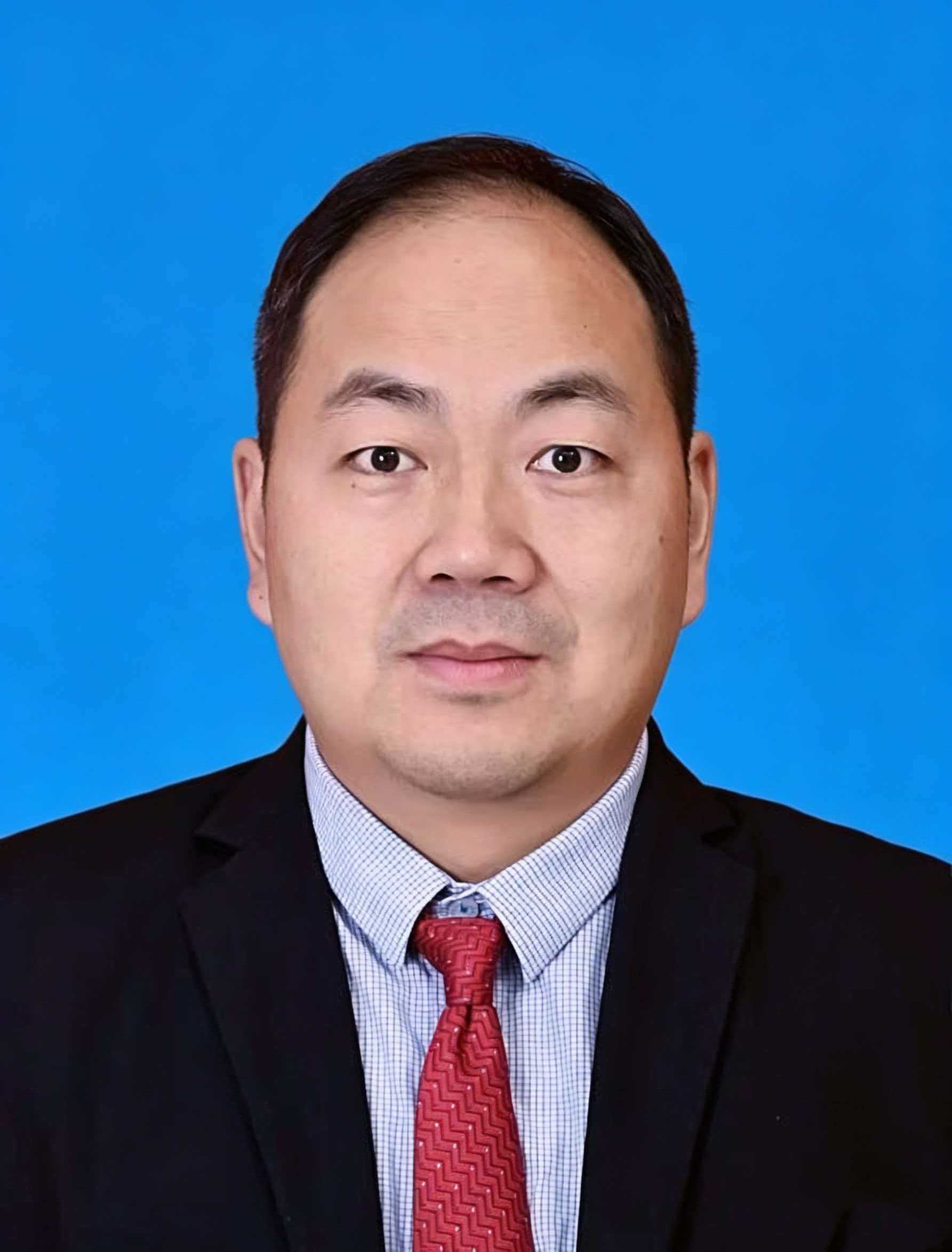}}]{Juntao Li}
Juntao Li received the B.S. and the M.S. degrees in applied mathematics from Henan Normal University, Xinxiang, China, and the Ph.D. degree in control theory and control engineering from Beihang University (BUAA), Beijing, China, in 2001, 2004, and 2010, respectively. Since 2019, he has been a professor at Henan Normal University. His research interests include multi-agent systems, data-driven modeling and control, and statistical machine learning.
\end{IEEEbiography}

\begin{IEEEbiography}[{\includegraphics[width=1in,height=1.25in,clip,keepaspectratio]{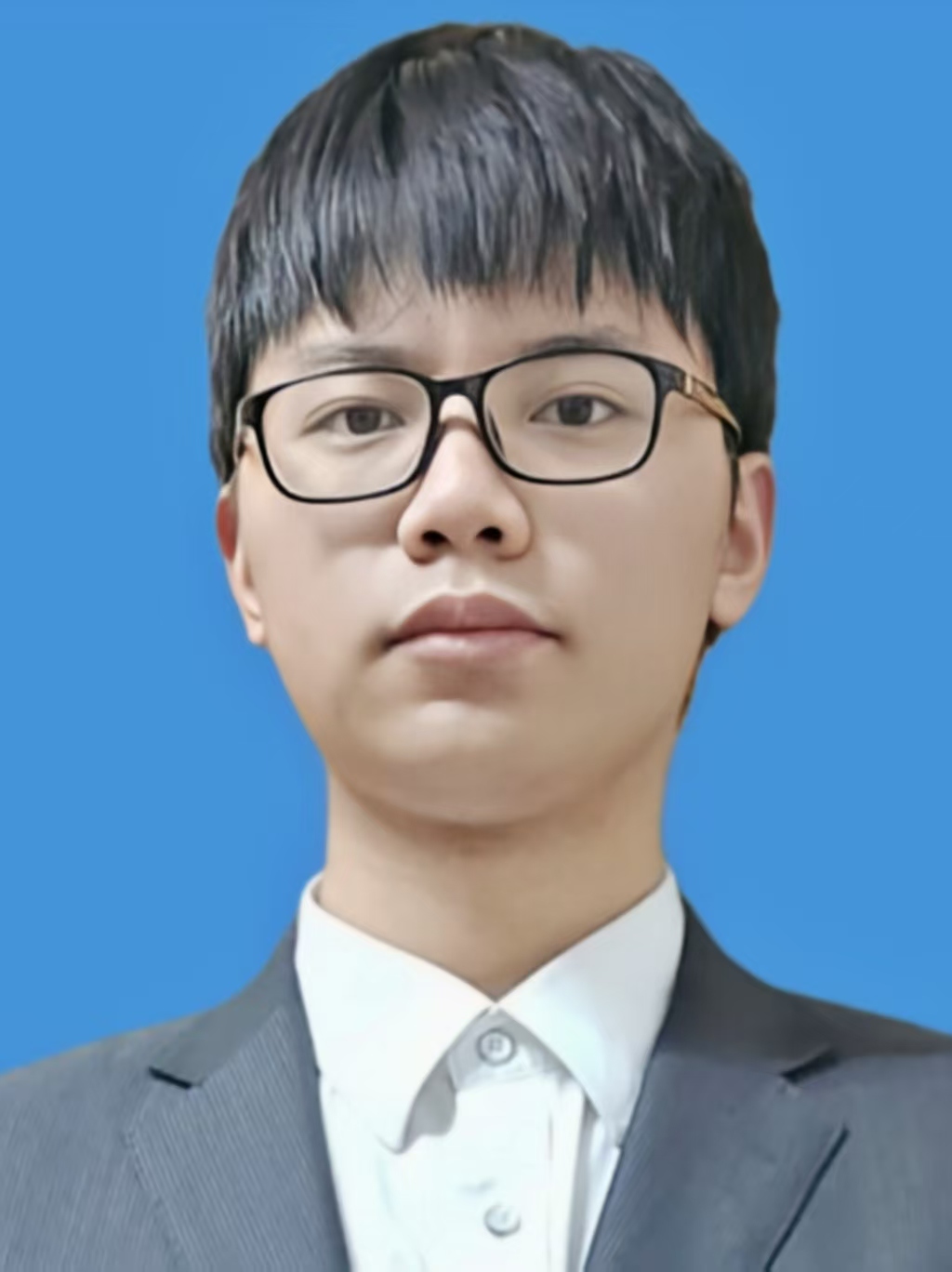}}]{Qihang Guo}
Qihang Guo received the B.S. degree in mathematics and applied mathematics from Henan Normal University, Xinxiang, China, in 2025. He is currently working toward the M.S. degree in operations research and control theory at the same university. His research interests include large-scale linear equations solving and optimization theory.
\end{IEEEbiography}

\begin{IEEEbiography}[{\includegraphics[width=1in,height=1.25in,clip,keepaspectratio]{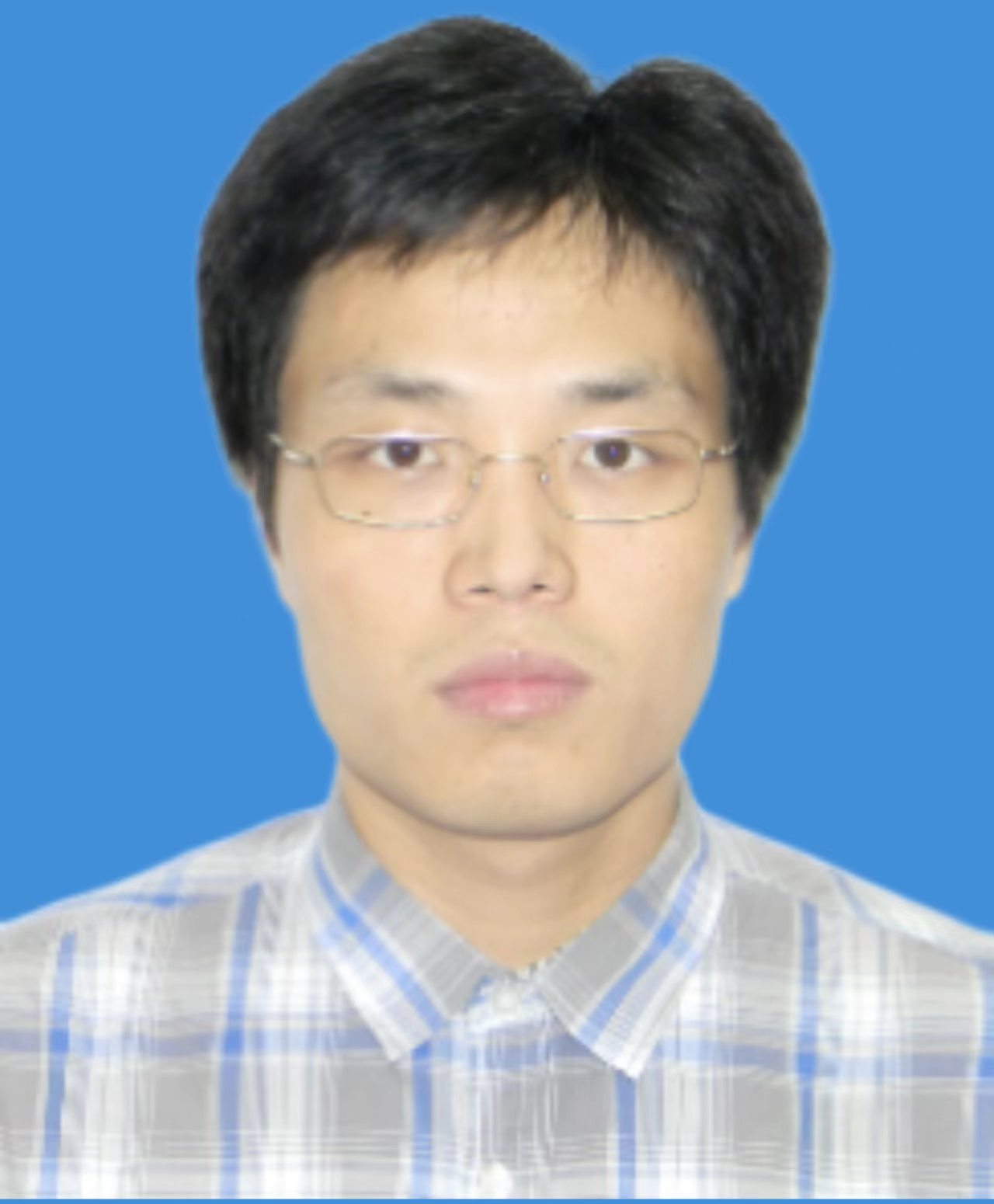}}]{Deyuan Meng}
Deyuan Meng (Senior Member, IEEE) received the B.S. degree in mathematics and applied mathematics from Ocean University of China, Qingdao, China, and the Ph.D. degree in control theory and control engineering from Beihang University (BUAA), Beijing, China. Since 2010, he has been with Beihang University (BUAA), where he is currently a Full Professor. From 2012 to 2013, he was a Visiting Scholar with the Department of Electrical Engineering and Computer Science, Colorado School of Mines, Golden, CO, USA. His research interests include iterative learning control, data-driven control, and multi-agent systems. He serves on the editorial boards of several journals and conferences, such as the International Journal of Systems Science, International Journal of Control, Automation, and Systems, and Control Theory and Technology.
\end{IEEEbiography}

\end{document}